\documentclass{amsart}

\usepackage{amssymb}
\usepackage{amsmath}
\usepackage[T1]{fontenc}

\input{macroced}
\usepackage{todonotes}
\usepackage{mathrsfs}
\usepackage{color,cite,graphicx}
\usepackage{tikz}
\definecolor{gray}{gray}{0.4}
\usepackage[all]{xy}
\usepackage{enumerate}
\usepackage{array, boldline, makecell, booktabs}
\usepackage{multirow,multicol}
\title[Complex reflection groups and K3 surfaces II]{Complex reflection groups and K3 surfaces II. \\
The groups ${\boldsymbol{G_{29}}}$, ${\boldsymbol{G_{30}}}$ and 
${\boldsymbol{G_{31}}}$}

\author{C\'edric Bonnaf\'e}
\address{IMAG, Universit\'e de Montpellier, CNRS, Montpellier, France} 

\email{cedric.bonnafe@umontpellier.fr}
\urladdr{http://imag.umontpellier.fr/~bonnafe/}
\makeatother

\author{Alessandra Sarti}
\address{	Universit\'e de Poitiers, 
Laboratoire de Math\'ematiques et Applications,
UMR 7348 CNRS, TSA 61125, 11 bd Marie et Pierre Curie,
86073 Poitiers Cedex 9, France}
\email{sarti@math.univ-poitiers.fr} 
\urladdr{http://www-math.sp2mi.univ-poitiers.fr/~sarti/}

\date{\today}

\thanks{The first author is partly supported by the ANR: 
project No. ANR-16-CE40-0010-01 (GeRepMod) and ANR-18-CE40-0024-02 (CATORE)\\
The second author is partly supported by the ANR: project No. ANR-20-CE40-0026-01 (SMAGP)}

\newtheorem{theorem}{Theorem}[section]
\newtheorem{pro}[theorem]{Proposition}
\newtheorem{lemma}[theorem]{Lemma}
\newtheorem{cor}[theorem]{Corollary}

\theoremstyle{remark}

\newtheorem{remark}[theorem]{Remark}
\DeclareMathOperator{\Pic}{Pic}

\def\petitespace{\vphantom{\DS{\frac{\DS{A}}{\DS{A}}}}}
\def\grandespace{\vphantom{\begin{pmatrix} \frac{A^1}{B_1} \\ \frac{A^1}{B_1} \end{pmatrix}}}

\def\slv{{\SSS{\mathrm{SL}}}}

\def\Deg{{\mathrm{Deg}}}
\def\Codeg{{\mathrm{Codeg}}}
\def\xyinj{\ar@{^{(}->}}

\relax 
\begin{document}

\begin{abstract}
We study some K3 surfaces obtained as minimal resolutions of quotients of subgroups of special reflection groups.
Some of these were already studied in a previous paper by W. Barth and the second author. 
We give here an easy proof that these are K3 surfaces, give equations in weighted projective space and describe their geometry.
\end{abstract}

\maketitle

\section{Introduction}

\medskip

In the first paper of this series~\cite{bonnafe-sarti-1}, 
the authors have explained how to build 
K3 surfaces from invariants of complex reflection groups 
of rank $4$ generated by reflections of order $2$. 
In this second part and the upcoming 
third part~\cite{bonnafe-sarti-3}, we complete this qualitative 
result by investigating more precisely the examples 
given by the primitive groups 
(see~\cite[\S\ref{sec:notation}]{bonnafe-sarti-1} for the definition), 
i.e. the groups $G_{28}$, $G_{29}$, $G_{30}$ and $G_{31}$ 
(as in~\cite{bonnafe-sarti-1}, 
we follow Shephard-Todd numbering for complex reflection 
groups~\cite{shephardtodd}). In particular, we investigate 
the following questions:
\begin{itemize}
\itemth{a} We show that all the K3 surfaces constructed this way 
have big Picard number compared  
to the number of moduli of the family they belong to.

\itemth{b} We compute some of the transcendental lattices of 
those K3 surfaces with Picard number 20.

\itemth{c} We give some explicit equations in 
weighted projective space.

\itemth{d} We construct explicit elliptic fibrations for all the examples of 
K3 surfaces we obtain: as shown in~\cite[Corollary~2.7]{sterk} there is only a finite number 
of elliptic fibrations for a K3 surfaces (up to automorphism) but, even though we sometimes 
contruct several non-equivalent elliptic fibrations, there is no case in which 
we pretend to have constructed all of them. For some of them, we determine 
the singular fibers. When one knows the transcendental lattice one could use the recent paper by Festi and Veniani \cite{festi-veniani}
to compute the number of elliptic fibrations (up to automorphism of the surface). 
\end{itemize}
In this second paper, 
we focus on the groups $G_{29}$, $G_{30}$ and $G_{31}$ 
while the third paper~\cite{bonnafe-sarti-3} will be devoted to 
the study of $G_{28}$. See the introduction of~\cite{bonnafe-sarti-3} 
for the reasons why $G_{28}$ deserves a particular treatment, we recall here the main points : firstly, $G_{28}'\neq G_{28}^\slv$;
secondly, there are two possible interesting degrees for the
fundamental invariants, namely $6$ and $8$; thirdly, $G_{28}$
admits an interesting outer automorphism. 
Also, we take opportunity of this work to revisit 
results from both authors who constructed highly singular 
surfaces from invariants of complex reflection 
groups~\cite{sartipencil},~\cite{bonnafesingular}. Most 
(but not all) of the singularities constructed this 
way can be obtained 
from~\cite[Corollary~\ref{coro:v-fixe-singulier}]{bonnafe-sarti-1}. 
In~\S\ref{sub:boissiere-sarti}, 
we also revisit Boissi\`ere-Sarti example of the smooth octic 
surface containing $352$ lines~\cite{boissieresarti}, 
using Springer theory~\cite[Theorem~\ref{theo:springer}]{bonnafe-sarti-1}. 

Finally, note that the first part~\cite{bonnafe-sarti-1} 
was free of computer calcutations, as the arguments 
were pretty general: in this second part, we study 
very specific examples, for which the determination 
of geometric features (singularities, transcendental 
lattices, branch locus,etc.) requires 
computer calculations. We use here the software 
{\sc Magma}~\cite{magma} (as well as some specific 
functions described in~\cite{bonnafemagma}). 

The structure of the paper is as follows. In section 3 we recall some general facts on the groups
$G_{29}$, $G_{30}$ and $G_{31}$. The section 4 is devoted to the group $G_{29}$. Here we consider the 
unique $G_{29}$-invariant polynomial of degree four which defines a quartic K3 surface in $\mathbb{P}^3(\mathbb{C})$,
this is denoted by $X_{\rm Mu}$ in \cite{bonnafe-sarti-m20}.  
We consider then the quotient of the quartic K3 surface
by the derived group $G_{29}'=G_{29}^{\rm SL}$. As remarked in \cite{bonnafe-sarti-m20} we have that $\mathbb{P}G_{29}'=M_{20}$ the Mathieu group, which acts symplectically on $X_{\rm Mu}$. It is well known that the minimal resolution is again a K3 surface and Xiao in \cite{xiao} showed that the  Picard number is 20. We give in Lemma 4.1 and Corollary 4.2 an alternative proof that uses Lehrer-Springer theory. 
In \S 4.3 we describe an elliptic fibration on this surface and thanks to that we compute the transcendental lattice. This result is new and summarized in the following theorem:
\begin{theorem}
Let $\tilde{X}^{29}$ be the minimal resolution of $X_{\rm Mu}/G_{29}'$. This is a K3 surface with Picard number 20 admitting an elliptic fibration with fibers $\tilde{E}_6+\tilde{D}_6+2\tilde{A}_2+\tilde{A}_1$ and transcendental lattice isometric to 
$$
\Tb_{\Xti^{29}}=\begin{pmatrix} 6 & 0 \\ 0 & 60 \end{pmatrix}.
$$
\end{theorem}
Observe that this surface was already studied under a different point of view by M. Sch\"utt in the 
paper \cite[Table~2]{schuett_extremal} about the construction of elliptic fibrations on extremal K3 
surfaces. Note that Sch\"utt constructs another elliptic fibration for $\Xti^{29}$: it 
would be interesting to know if there are still other elliptic fibrations.  
In section 5 we consider the group $G_{30}$ and the zero set  of the one dimensional family of invariant polynomials of degree 12. 
The group $G_{30}$ is the Coxeter group of type $H_4$. Let $\tilde{X}_{\lambda}^{30}$ denote the $K3$ surface which is the minimal resolution of the quotient of the zero set $Z(f_{2,\lambda})$ of the polynomials of degree 12 by $G_{30}^\slv$
and let  $X_{\lambda}^{30}$ denotes this singular quotient (recall that here $G_{30}'=G_{30}^\slv$). In \cite{barth-sarti}, \cite{sarti-trans} the Picard number and the transcendental lattice of the K3 surfaces were computed.  The equation of $X_\l^{30}$ and the description of the elliptic fibration is new. We show the following result :
\begin{theorem}\label{teo30}
We have the following equation
$$
X_\l^{30}=\{[y_1:y_3:y_4:j] \in \PM(1,2,3,6)~|~
j^2=r_\l(y_1,y_3,y_4)\}
$$
where $r_\l(y_1,y_3,y_4)$ is a polynomial of total degree $12$. For $\lambda$ generic the surface $\tilde{X}_{\lambda}^{30}$ has Picard number 19 and it admits an elliptic fibration with fibers $\tilde{D}_5+\tilde{A}_4+2\tilde{A}_2+3\tilde{A}_1$. The transcendental lattice as computed in \cite{sarti-trans} is
$$
\Tb_{\tilde{X}_{\lambda}^{30}}=\begin{pmatrix} 4 & 2& 0 \\ 2 & 34 & 0\\ 0& 0& -30 \end{pmatrix}.
$$
There are at least four special values of $\lambda$ for which the Picard number of the corresponding K3 surface is 20, these four values correspond to the surfaces in $Z(f_{2,\lambda})$ that have isolated $ADE$ singularities. These values of $\lambda$, the singular fibers of an elliptic fibration and the transcendental lattices are resumed in Table \ref{table3}
\end{theorem} 

Finally section 6 is devoted to the group $G_{31}$ and the one-dimensional family of invariant polynomials of degree 20. We have again $G_{31}^\slv=G_{31}'$, let $X_{\lambda}^{31}$ denote the singular quotient by $G_{31}^\slv$ and by $\tilde{X}_{\lambda}^{31}$ its minimal resolution.
The latter is a K3 surface and we show 
\begin{theorem}\label{risultatiG31}
We have the following equation
$$
X_\l^{31} = \{[y_1:y_2:y_4:j] \in \PM(2,1,2,5)~|~
j^2=q_\l(y_1,y_2,y_4)\}.
$$
where $q_\l(y_1,y_2,y_4)$ is a polynomial of total degree $10$. 
For $\lambda$ generic the surface $\tilde{X}_{\lambda}^{31}$ has Picard number 18 and it admits an elliptic fibration with singular fibers $\tilde{D}_7+3\tilde{A}_2+3\tilde{A}_1$. 
There are at least six special values of $\lambda$ for which the Picard number of the corresponding K3 surface is 19, five of these values
correspond to the singular fibers in $Z(f_{3,\lambda})$. These values of $\lambda$ and the singular fibers of an elliptic fibration are resumed in Table \ref{table4}
\end{theorem}

Finally in the Appendix we collect several useful results that allow to find  the equations of the K3 surfaces  and the elliptic fibrations. 
We remark that in the case of $G_{31}$ we described a one parameter family of K3 surfaces, we believe that this family is not isotrivial, but we could not prove it, for $G_{30}$ this was shown in \cite{barth-sarti}.

\bigskip
\noindent{\bf Acknowledgements}: We warmly thank A. Degtyarev and B. Naskr\k{e}cki for useful comments on $X_{\rm Mu}$ and X. Roulleau 
for explaining us how to use the Artin-Tate conjecture to compute the Picard number in Theorem \ref{risultatiG31} and for useful comments.
The first author is grateful to the MSRI to let him use its high performance computing 
facilities.
\bigskip

\noindent{\bf Hypothesis and notation.} 
We keep the notation introduced in~\cite{bonnafe-sarti-1}. 
We recall some of them. 
First, $V$ is a complex vector space and $W$ 
is a complex reflection group acting on $V$ of dimension $n$. If 
$v \in V \setminus \{0\}$, we denote 
by $[v] \in \PM(V)$ the line it defines (i.e. $[v]=\CM v$). 
If $S$ is a subset of $V$, we denote by $W_S$ (resp. $W(S)$) the setwise (resp. pointwise) 
stabilizer of $S$ (so that $W(S)$ is a normal subgroup of $W_S$ and 
$W_S/W(S)$ acts faithfully on $S$). 
The derived subgroup of $W$ will be denoted by $W'$, 
and we set $W^\slv=W \cap \Sb\Lb_\CM(V)$. The degrees (resp. codegrees) of $W$ 
(see~\cite[\S{3.1}]{bonnafe-sarti-1}) are denoted by $(d_1,d_2,\dots,d_n)$ 
and $(d_1^*,d_2^*,\dots,d_n^*)$ respectively. 

If $e \in \ZM_{\geqslant 1}$, we set
$$\d(e)=|\{1 \le k \le n~|~e~\text{divides}~d_k\}|\quad\text{and}\quad
\d^*(e)=|\{1 \le k \le n~|~e~\text{divides}~d_k^*\}|.$$
With this notation, we have
$$\d(e)=\max_{w \in W} \bigl(\dim V(w,\z_e)\bigr).$$
In particular, $\z_e$ is an eigenvalue of some element of $W$ 
if and only if $\d(e) \neq 0$ that is, if and only if $e$ divides 
some degree of $W$. In this case, we fix an element $w_e$ of $W$ such that 
$$\dim V(w_e,\z_e)=\d(e).$$
We set for simplification $V(e)=V(w_e,\z_e)$ and 
$W(e)=W_{V(e)}/W(V(e))$: this subquotient 
of $W$ acts faithfully on $V(e)$. 

However, 
from now on, and until the end of this paper, we assume 
that $n=\dim(V)=4$ and that $W \subset \Gb\Lb_\CM(V)$ 
is a primitive\footnote{Recall that $W$ is said {\it primitive} if there does not 
exist a decomposition $V = V_1 \oplus \cdots \oplus V_r$ with $r \ge 2$ and $V_k \neq 0$ 
such that $W$ permutes the $V_k$'s} 
complex reflection group. If $S$ is a K3 surface, 
we denote by $\Tb_S$ its transcendental lattice and by $\rhob(S)$ its 
Picard number.

If $C_1$,\dots, $C_r$ are curves on a surface $S$, 
the {\bfit intersection graph} will be represented as follows: vertices 
correspond to $C_1$,\dots, $C_r$ and are represented by circles 
(with no information if self-intersection is $-2$; otherwise, the self-intersection 
number is written inside the circle) and there is an edge between the vertices corresponding 
to $C_j$ and $C_{j'}$ if $C_j \cdot C_{j'} \neq 0$ (nothing 
more is written on the edge if $C_j \cdot C_{j'} = 1$; otherwise, 
the number $C_j \cdot C_{j'}$ is written above the edge).

The singular fibers of elliptic fibrations will be denoted as usual according to their 
intersection matrix: for instance, a singular fiber of type $\Dti_4$ is a fiber whose 
intersection matrix is the Cartan matrix of the extended Dynkin diagram of type $\Dti_4$ 
(in Kodaira's notation, it is of type $\Irm_0^*$). 
There remain some ambiguities (for types $\Ati_1$ and $\Ati_2$): we say that a singular 
fiber is of type $\Ati_1$ (resp. $\Ati_2$) if it is of type 
$\Irm_2$ (resp. $\Irm_3$) and will use Kodaira's notation (i.e. 
${\mathrm{III}}$ or ${\mathrm{IV}}$) for the other singular fibers 
whose intersection graph is of type $\Ati_1$ or $\Ati_2$.

\def\mweil{\Mrm\Wrm}
\def\triv{{\mathrm{Triv}}}

If $S$ is a K3 surface and $\ph : S \longto \PM^1(\CM)$ is an elliptic 
fibration admitting a section $\s : \PM^1(\CM) \longto S$, we denote by 
$\mweil_\s(\ph)$ (or simply $\mweil(\ph)$ if $\s$ is clear from the context) 
its Mordell-Weil group. In this case, we denote by $\triv_\s(\ph)$ (or $\triv(\ph)$) 
the {\it trivial lattice} of the fibration $\ph$, namely the lattice generated 
by the vertical divisors and the class of the image of $\s$. Then
\equat\label{eq:mordell}
\mweil_\s(\ph) \simeq \Pic(S)/\triv_\s(\ph).
\endequat

See the book \cite{SS} for more details on elliptic surfaces. 

We will often denote by $\begin{pmatrix} a & b & c \end{pmatrix}$ the $2\times 2$--matrix:
$$
\begin{pmatrix} a & b \\ b & c \end{pmatrix}.
$$

\section{Preliminaries on primitive complex 
reflection groups of rank $4$}

\medskip

Recall~\cite{shephardtodd} that there are five primitive complex reflection groups of 
rank $4$, and that they are denoted by $G_{28}$, $G_{29}$, $G_{30}$, $G_{31}$ and $G_{32}$. 
The first four are generated by reflections of order $2$ and $G_{32}$ 
is generated by reflections of order $3$. Note that $G_{28}$ (resp. $G_{30}$) 
is the Coxeter group of type $F_4$ (resp. $H_4$). 

\bigskip

When we do explicit computations, 
we use the models of the primitive 
complex reflection groups $W$ that were implemented in {\sc Magma} 
by the first author (almost copying files due to Michel~\cite{michel} and Thiel~\cite{thiel}) 
in a file 

\centerline{\tt primitive-complex-reflection-groups.m}

\noindent which can be downloaded 
in~\cite{bonnafemagma}. Most of them (but not all) are taken 
from~\cite{michel} or~\cite{thiel}. We do not pretend that these are the 
best models, but the interested reader might have a look 
at~\cite[Remark~1.3]{bonnafesingular} 
for a discussion about some of the advantages of these models. 

Representing $W$ as a subgroup of $\Gb\Lb_4(\CM)$ allows 
to identify $V$ with $\CM^4$ and we denote by $(x,y,z,t)$ the dual 
basis of the canonical basis of $\CM^4$. Therefore, 
$\CM[V]=\CM[x,y,z,t]$. 

A first advantage of the chosen models is that the group $W$ 
is implemented as a Galois stable subgroup of $\Gb\Lb_4(K)$ where $K$ 
is a finite Galois extension of $\QM$ (the fact that such a model 
always exists was proved by Marin-Michel~\cite{marin-michel}). 
This implies that we can find fundamental 
invariants in $\QM[x,y,z,t]$. For instance, 
with such a model for $W=G_{30}=\Wrm(H_4)$, the singular dodecic 
surfaces constructed by the second author~\cite{sartipencil} can be 
realized over $\QM$, as explained in~\cite[Proposition~1.1]{bonnafesingular}.

Another advantage of our models is that $W$ generally contains 
a big subgroup of monomial matrices (except for $W=G_{30}=\Wrm(H_4)$). 
This leads to expressions of fundamental invariants in terms 
of symmetric functions. For this reason, we introduce the following 
notation: if $m$ is a monomial in $x$, $y$, $z$, $t$, we denote by 
$\Sigma(m)$ the sum of the monomials obtained by permuting 
the variables. For instance,
$$\Sigma(x^4y)=x^4(y+z+t)+y^4(x+z+t)+z^4(x+y+t)+t^4(x+y+z)=\Sigma(xy^4).$$

\bigskip

\bigskip

\section{The groups $G_{29}$, $G_{30}$ and $G_{31}$}

\medskip

\boitegrise{{\bf Hypothesis.} {\it From now on, and until the 
end of this paper, we assume moreover that $W$ is one of the 
three primitive groups $G_{29}$, $G_{30}$ or $G_{31}$.}}{0.75\textwidth}

\medskip

Let us first recall in Table~\ref{table:degres-29-30-31} 
some specific data for these three groups that were contained 
in~\cite[Table~\ref{table:degres}]{bonnafe-sarti-1}. 

\begin{table}
$$\begin{array}{!{\vline width 2pt} c !{\vline width 1.3pt} c|c|c|c!{\vline width 2pt}}
\hlinewd{2pt}
 W & |W| & |W/\Zrm(W)| & |W'| &
\begin{array}{c}\Deg(W) \\ \Codeg(W) \end{array} \\
\hlinewd{2pt}
G_{29} & 7\,680 & 1\,920 & 3\,840 & 
\begin{array}{c}4, 8, 12, 20 \\ 0,8,12,16\end{array} \\
\hline
G_{30}=\Wrm(H_4) & 14\,400 & 7\,200 & 7\,200 &
\begin{array}{c}2, 12, 20, 30 \\ 0, 10, 18, 28\end{array} \\
\hline
G_{31} & 46\,080 & 11\,520 & 23\,040 &
\begin{array}{c}8, 12, 20, 24 \\ 0, 12, 16, 28 \end{array} \\
\hlinewd{2pt}
\end{array}
$$
\caption{Numerical informations for $G_{29}$, $G_{30}$ and $G_{31}$}\label{table:degres-29-30-31}
\end{table}

Note that the hypothesis implies that 
$$W'=W^\slv$$
is of index $2$ in $W$ (recall from~\cite{bonnafe-sarti-1} the notation  
$W^\slv=W \cap \Sb\Lb_\CM(V)$): we denote by $\s$ the non-trivial element of $W/W'$. According 
to~\cite[Theorem~\ref{theo:main}]{bonnafe-sarti-1}, 
the surface $\ZC(f)/W'$ is a K3 surface with ADE singularities 
(endowed with a non-symplectic automorphism given by the action of $\s$), 
provided that $f$ is a fundamental invariant of $W$ of degree 
$4$ if $W=G_{29}$, of degree $12$ of $W=G_{30}$ or of degree $20$ 
if $W=G_{31}$ such that $\ZC(f)$ has only ADE singularities 
(which is almost always the case\footnote{This is not true 
only for the one-parameter family of surfaces of degree $20$ 
built from $G_{31}$: in this family, only one surface does not have 
ADE singularities (it is in fact reducible). 
See Section~\ref{sec:g31} for details.}). 
Our aim in this paper is to study the geometry of the K3 surface 
with ADE singularities $\ZC(f)/W'$ and of its minimal resolution $\widetilde{\ZC(f)/W'}$: 
in particular, we prove that the informations given 
in Table~\ref{table:k3} are 
correct\footnote{Erratum: the singularities of the 
five singular surfaces of degree $20$ 
defined by fundamental invariants of $G_{31}$ given in 
Table~\ref{table:k3} differ from the ones given 
in~\cite[Table~4]{bonnafesingular}: in fact, there is a mistake 
in~\cite{bonnafesingular}, as can be checked with {\sc Magma} 
thanks to~\cite{bonnafemagma}, and the correct values are given 
in Table~\ref{table:k3}. See the correction statement 
at {\tt https://doi.org/10.1080/10586458.2018.1555778}.}. 

\begin{table}
$$\begin{array}{!{\vline width 2pt} c | c !{\vline width 1.3pt} c|c|c|c|c!{\vline width 2pt}}
\hlinewd{2pt}
W & 
d & 
\ZC_\sing(f)
& 
m & 
\text{singularities of $\ZC(f)/W'$} & 
\r
&
\Tb_{\widetilde{\ZC(f)/W'}} \petitespace
\\
\hlinewd{2pt}
G_{29} & 4 & \vide & 0 & D_4+2\,A_4+3\,A_2+A_1 & 20 & 
\grandespace \begin{pmatrix} 6 & 0 & 60 \end{pmatrix} \\
\hlinewd{1pt}
\multirow{5}{*}{$G_{30}$} & 
\multirow{5}{*}{12} & \vide & 1 & A_4+4\,A_2+5\,A_1 & \ge 19 & 
{\rm Theorem}\,\, \ref{teo30} \\
\cline{3-7}
&& 60\,A_1 & 0 & E_8 + 3\,A_2 + 4\,A_1 & 20 &
\grandespace \begin{pmatrix} 4 & 2 & 34 \end{pmatrix} \\
\cline{3-7}
&& 300\,A_1 & 0 & E_6+A_4+2\,A_2+4\,A_1 & 20 &
\grandespace \begin{pmatrix} 12 & 6 & 58 \end{pmatrix} \\
\cline{3-7}
&& 360\,A_1 & 0 & D_7+4\,A_2+3\,A_1 & 20 &
\grandespace \begin{pmatrix} 6 & 0 & 132 \end{pmatrix} \\
\cline{3-7}
&& 600\,A_1 & 0 & D_5+A_4+3\,A_2+3\,A_1 & 20 &
\grandespace \begin{pmatrix} 6 & 0 & 220 \end{pmatrix} \\
\hlinewd{1pt}
\multirow{6}{*}{$G_{31}$} & 
\multirow{6}{*}{20} & \vide & 1 & D_6+A_3+3\,A_2+2\,A_1 & \ge 18 &\petitespace\\
\cline{3-6}
&& 960\,A_1 & 0 & D_6+D_5+A_3+2\,A_2 & 19 &\petitespace\\
\cline{3-6}
&& 480\,A_1 & 0 & E_6+D_6+A_3+A_2+A_1 & 19 &\petitespace\\
\cline{3-6}
&& 1920\,A_1 & 0 & D_6+A_5+A_3+A_2+2\,A_1 & 19 &\petitespace\\
\cline{3-6}
&& 1440\,A_2 & 0 & D_6+D_5+3\,A_2+A_1 & 19 &\petitespace\\
\cline{3-6}
&& 640\,A_3 & 0 & D_6+2\,A_3 +2\,A_2+2\,A_1 & 19 &\petitespace\\
\hlinewd{2pt}
\end{array}$$
{(Here, 
$d=\deg(f)$, $m$ is the number of moduli of the family and 
$\r=\rhob(\widetilde{\ZC(f)/W')})$) }

{~}

\caption{K3 surfaces of the form $\widetilde{\ZC(f)/W'}$ for $W=G_{29}$, $G_{30}$ or 
$G_{31}$}\label{table:k3}
\end{table}

\medskip

\section{The group $G_{29}$}

\medskip

\boitegrise{{\bf Hypothesis.} {\it We assume in this section, 
and only in this section, that $W=G_{29}$.}}{0.75\textwidth}

\medskip

We have $G_{29}=\langle s_1,s_2,s_3,s_4\rangle$, where  
$$
\begin{array}{cc}
s_1=\begin{pmatrix}
1  &0 & 0 & 0\\
0  &1 & 0 & 0\\
0 & 0 & 1 & 0\\
0 & 0 & 0 &-1\\
\end{pmatrix}, &
s_2=\DS{\frac{1}{2}}\begin{pmatrix}
 1&  1&  i&  i\\
 1 & 1& -i& -i\\
-i&  i&  1& -1\\
-i  &i &-1&  1\\
\end{pmatrix},\\
&\\
s_3=\begin{pmatrix}
 0 &1& 0& 0\\
 1 &0 &0 &0\\
 0 &0 &1& 0\\
 0& 0 &0& 1\\
\end{pmatrix},&
s_4=\begin{pmatrix}
 1 &0& 0& 0\\
 0 &0 &1 &0\\
 0& 1 &0 &0\\
 0 &0 &0 &1\\
\end{pmatrix}. \\
\end{array}
$$
Recall that $i=\z_4$, a primitive fourth root of unity. Some of the numerical facts used below 
can be extracted from~\cite[Table~\ref{table:degres}]{bonnafe-sarti-1}.
For instance, note that $G_{29}^\slv=G_{29}'$ is of index $2$ in $G_{29}$, 
so that $G_{29}=\langle s_1 \rangle \ltimes G_{29}'$ (observe that $s_1$ is an involution). 
Note also that $\Zrm(G_{29}) \simeq \mub_4 \subset G_{29}'$, see for 
instance~\cite[$($3.2$)$]{bonnafe-sarti-1}. 
Moreover, $PG_{29}' \simeq M_{20}$ (the Mathieu group of degree $20$) 
so that we have a split exact sequence 
$$1 \longto PG_{29}' \simeq M_{20} \longto PG_{29} \longto \mub_2 
\longto 1,$$
where the last map is induced by the determinant. 

\bigskip
\def\mukai{{\mathrm{Mu}}}

\subsection{The K3 surface} 
By~\cite[Table~\ref{table:degres}]{bonnafe-sarti-1}, 
there exists a unique (up to scalar) 
homogeneous invariant polynomial $f_1$ of degree $4$: it is given by  
$$f_1=\Sigma(x^4)-6\Sigma(x^2y^2).$$
As in~\cite{bonnafe-sarti-m20}, we set $X_\mukai=\ZC(f)$ (recall 
that this surface was discovered by Mukai~\cite{mukai}).
It can easily be checked that $X_\mukai$ is a 
smooth and irreducible quartic in $\PM^3(\CM)$, so that it is a K3 surface, 
endowed with a symplectic action of $M_{20}$ and an extra non-symplectic 
automorphism of order $2$. Several properties of $X_\mukai$ are given 
in~\cite{bonnafe-sarti-m20} (transcendental lattice, automorphisms, 
polarizations: note that it is denoted by $X_{\mathrm{Mu}}$ 
in~\cite{bonnafe-sarti-m20}, as it was discovered by Mukai~\cite{mukai}). 
For instance, it is known that $X_\mukai$ 
has Picard number $20$.

Continuing with the topic of this paper, we describe here geometric properties 
of the quotient 
$X^{29}=X_\mukai/G_{29}'$: as the quotient of a K3 surface by a finite group acting 
symplectically, it is also a K3 surface with ADE singularity, 
whose minimal resolution $\Xti^{29}$ has Picard number $20$ 
(see~\cite{bonnafe-sarti-m20}). This can also be proved by examining 
the singularities of $X^{29}$, which are given below: 

\bigskip

\begin{lemma}[Xiao]\label{prop:sing-q4}
The K3 surface $X^{29}$ has singularities 
$D_4+2\,A_4+3\,A_2+A_1$.
\end{lemma}

\bigskip

\begin{proof}
See~\cite[Table~2,~last line]{xiao}. As we need concrete 
results (for instance, the coordinates of the singular points), we provide 
a proof that will provide these extra-informations. 

Since the action of $PG_{29}'$ on $X_\mukai$ is symplectic, it is sufficient 
to compute the stabilizers of points of $X_\mukai$. For this, we follow 
the discussion of~\cite[\S\ref{sub:stab}]{bonnafe-sarti-1}, from 
which we keep the notation. 
We fix $v \in V\setminus \{0\}$ such that $z=[v] \in X_\mukai$ 
and we may assume that $W_z=W_v\langle w_{e_z} \rangle$. 
Note that $4$ divides $e_z$ because $w_4=\z_4 \Id_V \in W$, 
and that $e_z$ divides one of the degrees of $W$. So $e_z \in \{4,8,12,20\}$. 
This leads to the following discussion, by using {\sc Magma}:
\begin{itemize}
\item[$\bullet$] If $e_z = 20$ then, since $\d(20)=\d^*(20)=1$, 
we have $W_v=1$ by~\cite[Theorem~\ref{theo:springer}]{bonnafe-sarti-1} and 
$\det(w_{20})=\z_{20}^{4+8+12+20-4}=1$. So the stabilizer of 
$z$ in $W$ is contained in $W'$: so the $W$-orbit of $z$ contains 
two $W'$-orbits, and the stabilizer of $z$ in $PG_{29}'$ 
is cyclic of order $5$. This leads to $2\,A_4$ singularities 
in $X^{29}$, which we denote by $a_4^\pm$

\item[$\bullet$] If $e_z=12$, then $\d(12)=1$, so 
the $W$-orbit of $z$ is completely determined, 
and a computation with {\sc Magma} shows that $|W_z|=24$. 
This shows that $W_v$ is generated by a single reflection and 
so $W_z'=\langle w_{12} \rangle \varsubsetneq W_z$. 
So the $W$-orbit of $z$ is a single $W'$-orbit, and the stabilizer 
of $z$ in $PG_{29}'$ 
is cyclic of order $3$. This leads to an $A_2$ singularity 
in $X^{29}$, which we denote by $a_2$.

\item[$\bullet$] If $e_z=8$, then $\d(8)=1$, so 
the $W$-orbit of $z$ is completely determined, 
and a computation with {\sc Magma} shows that $|W_z|=64 \neq 32=|W_z'|$. 
So the $W$-orbit of $z$ is a single $W'$-orbit, and one checks with {\sc Magma} 
that the stabilizer of $z$ in $PG_{29}'$ is the quaternionic group 
of order $8$. This leads to a $D_4$ singularity in $X^{29}$, which we denote by $d_4$.

\item[$\bullet$] Assume now that $e_z=4$. If $|W_v|=1$ or $2$, then 
$W_z'=\langle w_4 \rangle$ and so the stabilizer of $z$ in $PG_{29}'$ 
is trivial. So the image of $z$ in $X^{29}$ is smooth. 
By~\cite[Corollary~\ref{coro:v-fixe-singulier}]{bonnafe-sarti-1}, 
the group $W_v$ cannot have rank $3$, for otherwise $z$ 
would be singular in $X_\mukai$. So $W_v$ has rank $2$. There 
are three conjugacy classes of parabolic subgroups of rank $2$, 
and representatives are given by
$$W_{12}=\langle s_1,s_2\rangle,\qquad W_{13} = \langle s_1,s_3 \rangle
\qquad \text{and}\qquad W_{23}=\langle s_2,s_3 \rangle.$$
We denote by $L_{jk}$ the projective line $\PM(V^{W_{jk}})$ in $\PM(V)$. 
Since $X_\mukai$ is smooth, it follows that $L_{jk}$ meets 
$X_\mukai$ transversally~\cite[Corollary~\ref{coro:transverse}]{bonnafe-sarti-1}, and we set 
$E_{jk}=L_{jk} \cap X_\mukai$. Then $|E_{jk}|=4$ and it follows 
from~\cite[\S{sub:stab},~(c)]{bonnafe-sarti-1} that two elements of $\O_{jk}$ 
are in the same $W'$-orbit if and only if they are in the same 
$(W' \cap N_{jk})$-orbit. Now the next results can be obtained 
with {\sc Magma}:
\begin{itemize}
\item The group $W_{12}$ is of type $A_2$ and  
$|(W' \cap N_{12})/W_{12}\langle w_4 \rangle| = 2$. 
Moreover, the stabilizer of any point in $E_{12}$ 
is equal to $(W' \cap W_{12})\langle w_4 \rangle$, so its stabilizer 
in $PW'$ is cyclic of order $3$. This leads to $2\,A_2$ 
singularities in $X^{29}$, which we denote by $a_2^\pm$. 
\item The group $W_{13}$ is of type $A_1 \times A_1$ and  
$|(W' \cap N_{13})/W_{13}\langle w_4 \rangle| = 4$. 
Moreover, the stabilizer of any point in $E_{13}$ 
is equal to $(W' \cap W_{13})\langle w_4 \rangle$, so its stabilizer 
in $PW'$ is cyclic of order $2$. This leads to an $A_1$ 
singularities in $X^{29}$, which we denote by $a_1$. 
\item The set $E_{23}$ is contained in the $W$-orbit of $z_8$, so 
this case has already been treated and does not lead to new singularities 
in $X^{29}$.
\end{itemize}
\end{itemize}
The proof of the proposition is complete.
\end{proof}

\bigskip

\begin{cor}\label{cor:q4-20}
The Picard number of the K3 surface $\Xti^{29}$ is $20$.
\end{cor}

\bigskip

\begin{proof}
As $\Xti^{29}$ is algebraic, we get from 
Lemma~\ref{prop:sing-q4} that the rank of 
${\mathrm{Pic}}(\Xti^{29})$ is 
$$\ge 1 + (1 + 3 \cdot 2 + 2 \cdot 4 + 4)=20.$$
Since this rank is bounded above by $20$ for a K3 surface, 
this yields the result.
\end{proof}

\bigskip

\begin{remark}\label{rem:eq-g29}
With a suitable choice of a family $\fb$ of fundamental invariants 
and a suitable normalization of $J$, one gets  
with {\sc Magma} that
$$X^{29}=\{[x_2 : x_3 : x_4 : j] \in \PM(2,3,5,10)~|~$$
$$j^2=-64x_2^5x_4^2 + 16x_2^4x_3^4 + 32x_2^3x_3^3x_4 
+ 1800x_2^2x_3^2x_4^2 $$
$$- 432x_2x_3^6 - 5000x_2x_3x_4^3 
+ 432x_3^5x_4 + 3125x_4^4\quad\}$$
(see~\cite[Proposition~\ref{prop:zf-wsl}]{bonnafe-sarti-1}). In this model, the singular 
points are given as follows, thanks to {\sc Magma}:
$$d_4=[1:0:0:0],\qquad a_2=[0:1:0:0],\qquad a_4^\pm=[0:0:1:\pm 25 \sqrt{5}]$$
$$a_1=[\a : \b : 1 : 0]\qquad\text{and}\qquad a_2^\pm=[\a_\pm : \b_\pm : 1 : 0],$$
where 
$$\begin{cases}
(\a^5,\a\b,\b^5)=(84375/16,-25/2,3125/54),\\
(\a_\pm^5,\a_\pm\b_\pm,\b_\pm^5)=
((3987\pm 1632\sqrt{6})/16,(7\pm 2 \sqrt{6})/2,(117\pm 62\sqrt{6})/18).
\end{cases}$$
Indeed, it follows from~\cite[Lemma~\ref{lem:trivial}]{bonnafe-sarti-1} 
that $d_4$, $a_2$ and $a_4^\pm$ have coordinates of the form 
$[1:0:0:j_1]$, $[0:1:0:j_2]$ and $[0:0:1:j_3^\pm]$ respectively, and the values 
of $j_1$, $j_2$ and $j_3^\pm$ are determined by the equation of $X^{29}$. 

For the remaining points, computations with {\sc Magma} 
show that the evaluation 
of $f_4$ (the invariant of degree $20$) at points of $E_{12}$ or $E_{13}$ is different from $0$, so that the points 
$a_1$ and $a_2^\pm$ belong to the affine chart $X^{29}_{(4)}$ of $X^{29}$ defined by $x_4 \neq 0$. 
Setting $x_4=1$, 
the coordinates in the ambient space of $X^{29}_{(4)}$ are $a=x_2^5$, $b=x_2x_3$, $c=x_3^5$ 
and $j$, and $X^{29}_{(4)}$ is equal to
\begin{multline*}
X^{29}_{(4)}=\{(a,b,c,j) \!\in\! \AM^{\!4}(\CM)~|~b^5=ac~\text{and} \\
j^2=-64a+16b^4+32b^3+1800b^2-432bc-5000b+432c+3125\}.
\end{multline*}
From the second equation, we can express $a$ in terms of $b$, $c$ and $j$, and so 
\equat\label{eq:x4-g29}
X^{29}_{(4)}=\{(b,c,j) \!\in\! \AM^{\!3}(\CM)~|~ 
64b^5=c P_{29}(b,c,j)\},
\endequat
where $P_{29}(b,c,j)=16b^4+32b^3+1800b^2-432bc-5000b+432c+3125-j^2$.
The coordinates of the singular points of $X_{(4)}^{29}$ can then be computed with {\sc Magma} 
and fit with what is written above.

\medskip

Finally, recall that the action of the non-trivial element $\s$ of $W/W'$ 
is given by
\equat\label{eq:g29-sigma}
\s \cdot [x_2:x_3:x_4:j] = [x_2:x_3:x_4:-j]
\endequat
that $X/\langle \s \rangle \simeq \PM(2,3,5)$.\finl
\end{remark}

\bigskip

\def\cartan{{\mathrm{Car}}}

\subsection{Some smooth rational curves in ${\boldsymbol{X^{29}}}$}
We work in the model given by Remark~\ref{rem:eq-g29} and we denote 
by $\pi : \Xti^{29} \to X^{29}$ the natural morphism from the minimal resolution 
$\Xti^{29}$ of $X^{29}$. If $p$ is a singular point of $X^{29}$, we denote by 
$\D_1^p$,\dots, $\D_m^p$ the irreducible components of $\pi^{-1}(p)$ 
(these are smooth rational curves and $m$ is equal to the Milnor number 
of $X^{29}$ at $p$). We define
$$C_2=\{[x_2 : x_3 : x_4 : j] \in X^{29}~|~x_2=0\}$$
$$C_3=\{[x_2 : x_3 : x_4 : j] \in X^{29}~|~x_3=0\}.\leqno{\text{and}}$$
Let $\Cti_2$ and $\Cti_3$ denote the respective strict transforms 
of $C_2$ and $C_3$ in $\Xti^{29}$. 

\bigskip

\begin{pro}\label{prop:c2-c3-r-rational}
The curves $C_2$ and $C_3$ are smooth rational curves.
\end{pro}

\bigskip

\begin{proof}
First,
$$C_3 = \{[x_2 : x_4 : j] \in \PM(2,5,10)~|~j^2=-64x_2^5x_4^2 + 3125x_4^4\}.$$
But the map $\PM(2,5,10) \to \PM^2(\CM)$, $[x_2,x_4,j] \mapsto [x_2^5,x_4^2,j]$ 
is an isomorphism of varieties. Through this isomorphism, we get 
$$C_3 = \{[x_2 : x_4 : j] \in \PM^2(\CM)~|~j^2=-64x_2x_4 + 3125x_4^2\}.$$
Hence, $C_3$ is a non-degenerate conic in $\PM^2(\CM)$, i.e. $C_3$ is a smooth 
rational curve.

For $C_2$, note that 
$$C_2=\{[x_3:x_4:j] \in \PM(3,5,10)~|~j^2=432x_3^5x_4 + 3125x_4^4\}.$$
But $\PM(3,5,10) \simeq \PM(3,1,2)$ and, through this isomorphism, 
one gets
$$C_2=\{[x_3:x_4:j] \in \PM(3,1,2)~|~j^2=432x_3x_4 + 3125x_4^4\}.$$
For $k \in \{3,4\}$, we denote by $C_2^{(k)}$ the affine chart of $C_2$ 
defined by $x_k \neq 0$. Then $C_2=C_2^{(3)} \cup C_2^{(4)}$ and we only 
need to show that $C_2^{(3)}$ and $C_2^{(4)}$ are smooth rational 
affine curves. For $C_2^{(4)}$, this is obvious. For $C_2^{(3)}$, working
with the coordinates $a=x_4^3$, $b=x_4j$ and $c=j^3$, one gets 
\eqna
C_2^{(3)}&=&
\{(a,b,c) \in \AM^{\!3}(\CM)~|~b^3=ac, c=432b+3125ab~\text{and}~b^2=432a+3125a^2\} \\
&\simeq& \{(a,b) \in \AM^{\!2}(\CM)~|~b^2=432a+3125a^2\}.
\endeqna
This is clearly smooth and the result follows.
\end{proof}

\bigskip

Proposition~\ref{prop:c2-c3-r-rational} implies 
that $\Cti_2$ and $\Cti_3$ are smooth rational 
curves in $\Xti^{29}$. Adding the 19 smooth rational 
curves of the form $\D_m^p$, this gives us 21 
smooth rational curves in $\Xti^{29}$, we investigate in the next subsection if these curves are independent in the Picard group or not. 

\bigskip

\subsection{An elliptic fibration}
Any K3 surface with Picard number $20$ admits an elliptic fibration. 
We construct here an explicit one, and determine its singular fibers. 

First, let $\ph : X^{29} \setminus \{a_4^+,a_4^-\} \to \PM^1(\CM)$, 
$[x_2:x_3:x_4:j] \mapsto [x_2^3:x_3^2]$. This map is indeed well-defined 
on $X^{29} \setminus \{a_4^+,a_4^-\}$ and induces a map 
$$\pht : \Xti^{29} \setminus (\pi^{-1}(a_4^+) \cup \pi^{-1}(a_4^-)) \longto \PM^1(\CM).$$
Our elliptic fibration is obtained by extending $\pht$:

\bigskip

\begin{pro}\label{pro:elliptic-g29}
The map $\pht : \Xti^{29} \setminus (\pi^{-1}(a_4^+) \cup \pi^{-1}(a_4^-)) \longto \PM^1(\CM)$ 
extends to a morphism of algebraic varieties $\Xti^{29}\longto \PM^1(\CM)$.
\end{pro}

\bigskip

\begin{proof}
Let $\pih : \Xhat^{29} \longto X^{29}$ denote the minimal resolution of $X^{29}$ 
{\it only at the points 
$a_4^+$ and $a_4^-$.} In particular, $\Xhat^{29}$ is still singular (it has 
singularities $A_1 + 3\,A_2 + D_4$). Let 
$\phh =\ph \circ \pih : \Xhat^{29} \setminus (\pih^{-1}(a_4^+) \cup \pih^{-1}(a_4^-)) 
\longto \PM^1(\CM)$. Since the resolution $\pi : \Xti^{29} \to X^{29}$ factors through 
$\Xhat^{29}$, it is sufficient to show that $\phh$ extends to $\Xhat^{29}$. 

We will now use the results (and the notation) of Appendix~\ref{appendix-2} with $(k,l)=(2,3)$. 
Let $a_4=[0:0:1] \in \PM(2,3,5)$: it is the image of $a_4^+$ (or $a_4^-$) 
through the quotient morphism $X^{29} \longto \PM(2,3,5)$. 
Now, the map $\ph : X^{29} \setminus \{a_4^+,a_4^-\} \to \PM^1(\CM)$ is the composition of 
the quotient $X^{29} \setminus \{a_4^+,a_4^-\} \longto \PM(2,3,5) \setminus \{a_4\}$ 
and the map $\ph_{2,3} : \PM(2,3,5) \setminus \{a_4\} \longto \PM^1(\CM)$ defined 
in Appendix~\ref{appendix-2}. Therefore, $\phh$ is the composition 
$$\diagram
\Xhat^{29} \setminus (\pih^{-1}(a_4^+) \cup \pih^{-1}(a_4^-)) \rrto &&
\hat{\PM}(2,3,5) \setminus \{a_4\} 
\rrto^{\DS{\hphantom{A}\phh_{2,3}}} && \PM^1(\CM),
\enddiagram$$
where the first map is the quotient by the lift of $\s$. 
So the result follows from the fact that $\phh_{2,3}$ extends to $\hat{\PM}(2,3,5)$ 
(see~\eqref{eq:map-p1}).
\end{proof}

\bigskip

\begin{remark}\label{rem:section-g29}
Let us describe two sections of the elliptic fibration $\pht$. For this, keep 
the notation $\Xhat^{29}$, $\phh$ of the proof of 
Proposition~\ref{pro:elliptic-g29}. 
The map $\phh$ factors through
$$\diagram
\Xhat^{29} \rrto &&
\hat{\PM}(2,3,5) 
\rrto^{\DS{\hphantom{A}\phh_{2,3}}} && \PM^1(\CM),
\enddiagram$$the first map being the quotient by the action of $\s$. We denote by $\D_1$,\dots, $\D_4$ 
the smooth rational curves defined in Appendix~\ref{appendix-2} and by 
$\Delh_1^{a_4^\pm}$,\dots, $\Delh_4^{a_4^\pm}$ the smooth rational curves of 
the exceptional divisor of $\Xhat^{29}$ above $a_4^\pm$. 
Since $X^{29} \longto \PM(2,3,5)$ is unramified above $[0:0:1]$, 
we can number those last curves so that $\Delh_j^{a_4^\pm} \to \D_j$ is an isomorphism. 

Now, by Remark~\ref{rem:poids-section}, 
the map $\phh_{2,3}$ admits a section $\th : \PM^1(\CM) \longto \hat{\PM}(2,3,5)$ 
whose image is $\D_2$. This yields two sections $\theh^\pm : \PM^1(\CM) \to \Xhat^{29}$, 
whose image is $\D_2^{a_4^\pm}$. Now, $\Xti^{29}$ is obtained from $\Xhat^{29}$ 
by successive blow-ups of points not lying in $\D_2^{a_4^+} \cup \D_2^{a_4^-}$, 
so $\theh^\pm$ lifts to a section $\thet^\pm : \PM^1(\CM) \longto \Xti^{29}$, 
note that $\thet^-=\s \circ \thet^+$.\finl
\end{remark}

\bigskip

Let $[u:v] \in \PM^1(\CM)$. We denote by $X^{29}_{u,v}$ the Zariski closure 
of $\ph^{-1}([u:v])$ (endowed with its reduced structure) in $X^{29}$ 
and by $\Xti^{29}_{u,v}$ its strict transform in 
$\Xti^{29}$. Note that $\Xti_{u,v}^{29} \subset \pht^{-1}([u:v])$ and that 
\eqna
X^{29}_{u,v} &=&\Bigl(\{[x_2:x_3:x_4:j] \in X~|~vx_2^3=ux_3^2\}\Bigr)_{\mathrm{red}} \\
&=& \ph^{-1}([u:v]) \cup \{a_4^+,a_4^-\},
\endeqna
where $Y_{\mathrm{red}}$ denotes the reduced subscheme of $Y$ (this is necessary 
only if $uv=0$). 
%
%

\bigskip

\begin{cor}\label{coro:singular-fibers-g29}
The elliptic fibration $\pht$ has  
singular fibers $\Eti_6 + \Dti_6 + 2\, \Ati_2 + \Ati_1$.
\end{cor}

\bigskip

\begin{proof}
Since the map $\pht$ factorizes through the quotient $\Xti^{29}/\langle \s \rangle$, 
it follows from Proposition~\ref{prop:graphe-eclatement} that the intersection 
graph of the family of smooth rational curves 
$(\Cti_2,\Cti_3,\D_1^{a_4^+},\D_2^{a_4^+},\D_3^{a_4^+},\D_4^{a_4^+},
\D_1^{a_4^-},\D_2^{a_4^-},\D_3^{a_4^-},\D_4^{a_4^-})$ is given by
\vskip1.2cm
\equat\label{eq:intersection-g29-1}
\text
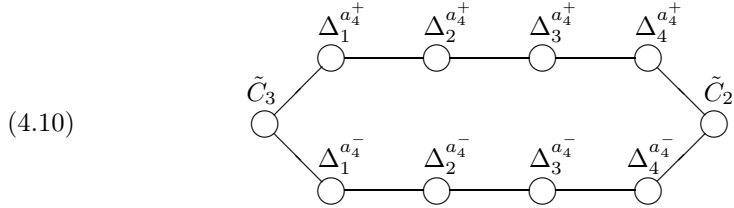
{\centerline{\begin{picture}(335,24)
\put( 85, 30){\circle{10}}\put(78,39){$\Cti_3$}
\put( 88.5355,33.5355){\line(1,1){18.1}}
\put( 88.5355,26.4645){\line(1,-1){18.1}}
\put(110, 55){\circle{10}}\put(115, 55){\line(1,0){30}}\put(105,64){$\D_1^{a_4^+}$}
\put(110,  5){\circle{10}}\put(115,  5){\line(1,0){30}}\put(105,14){$\D_1^{a_4^-}$}
\put(150, 55){\circle{10}}\put(155, 55){\line(1,0){30}}\put(145,64){$\D_2^{a_4^+}$}
\put(150,  5){\circle{10}}\put(155,  5){\line(1,0){30}}\put(145,14){$\D_2^{a_4^-}$}
\put(190, 55){\circle{10}}\put(195, 55){\line(1,0){30}}\put(185,64){$\D_3^{a_4^+}$}
\put(190,  5){\circle{10}}\put(195,  5){\line(1,0){30}}\put(185,14){$\D_3^{a_4^-}$}
\put(230, 55){\circle{10}}\put(225,64){$\D_4^{a_4^+}$}
\put(230,  5){\circle{10}}\put(222,14){$\D_4^{a_4^-}$}
\put(233.5355,51.4645){\line(1,-1){18.1}}
\put(233.5355,8.5355){\line(1,1){18.1}}
\put(255, 30){\circle{10}}\put(251,39){$\Cti_2$}
\end{picture}}}
\endequat
Moreover, Proposition~\ref{prop:graphe-eclatement} also 
shows that $\D_1^{a_4^+}$ and $\D_1^{a_4^-}$ (resp. $\D_3^{a_4^+}$, $\D_4^{a_4^+}$, 
$\D_3^{a_4^-}$ and $\D_4^{a_4^-}$) are the only rational curves among the $\D_k^{a_4^\pm}$'s 
which are contained in $\pht^{-1}([1:0])$ (resp. $\pht^{-1}([0:1])$). 

This shows that $\pht^{-1}([1:0])$ and $\pht^{-1}([0:1])$ are singular 
fibers. Let us determine their type. Note that 
$$C_3=X^{29}_{1,0}\qquad\text{and}\qquad C_2=X^{29}_{0,1}.$$
As the only singular points of $\Xti^{29}$ belonging to $C_3$ (resp. $C_2$) are 
$a_4^+$, $a_4^-$ and $d_4$ (resp. $a_4^+$, $a_4^-$ and $a_2$), this shows that 
$$\pht^{-1}([1:0]) = \Cti_3 \cup \D_1^{a_4^+} \cup \D_1^{a_4^-} 
\cup \bigl(\bigcup_{k=1}^4 \D_k^{d_4}\bigr)$$
$$\pht^{-1}([0:1]) = \Cti_2 \cup \bigl( \bigcup_{k=3}^4 (\D_k^{a_4^+} \cup \D_k^{a_4^+}) \bigr) 
\cup \bigl(\bigcup_{k=1}^2 \D_k^{a_2}\bigr).\leqno{\text{and}}$$
But $\D_k^{a_4^\pm} \cdot \D_l^{d_4} = \D_k^{a_4^\pm} \cdot \D_m^{a_2} = 0$, 
so the Kodaira-N\'eron classification of singular fibers 
forces that, with a suitable numbering of the $\D_k^{d_4}$'s and the $\D_k^{a_2}$'s, 
the intersection graphs inside $\pht^{-1}([1:0])$ and $\pht^{-1}([0:1])$ 
are respectively given by
\vskip1cm
\equat\label{eq:intersection-g29-2}
\text
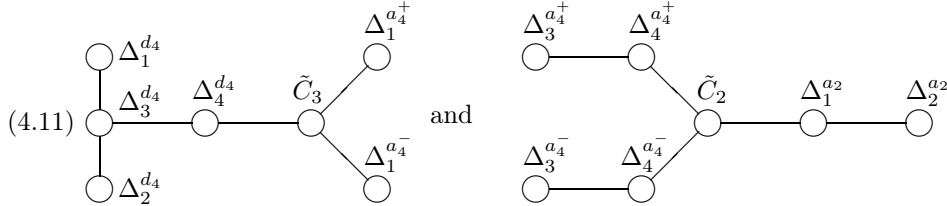
{\centerline{\begin{picture}(300,24)
\put(  5, 30){\circle{10}}\put(12,35){$\D_3^{d_4}$}
\put( 10, 30){\line(1,0){30}}
\put( 45, 30){\circle{10}}\put(40,39){$\D_4^{d_4}$}
\put( 50, 30){\line(1,0){30}}
\put(  5, 55){\circle{10}}\put(12,2){$\D_2^{d_4}$}
\put(  5,  5){\circle{10}}\put(12,54){$\D_1^{d_4}$}
\put(  5, 35){\line(0,1){15}}
\put(  5, 25){\line(0,-1){15}}
\put( 85, 30){\circle{10}}\put(78,39){$\Cti_3$}
\put( 88.5355,33.5355){\line(1,1){18.1}}
\put( 88.5355,26.4645){\line(1,-1){18.1}}
\put(110, 55){\circle{10}}
\put(105,64){$\D_1^{a_4^+}$}
\put(110,  5){\circle{10}}
\put(105,14){$\D_1^{a_4^-}$}
\put(130,30){\text{and}}
\put(170, 55){\circle{10}}\put(175, 55){\line(1,0){30}}\put(165,64){$\D_3^{a_4^+}$}
\put(170,  5){\circle{10}}\put(175,  5){\line(1,0){30}}\put(165,14){$\D_3^{a_4^-}$}
\put(210, 55){\circle{10}}\put(205,64){$\D_4^{a_4^+}$}
\put(210,  5){\circle{10}}\put(202,14){$\D_4^{a_4^-}$}
\put(213.5355,51.4645){\line(1,-1){18.1}}
\put(213.5355,8.5355){\line(1,1){18.1}}
\put(235, 30){\circle{10}}\put(231,39){$\Cti_2$}
\put(240, 30){\line(1,0){30}}
\put(275, 30){\circle{10}}\put(280, 30){\line(1,0){30}}\put(270,39){$\D_1^{a_2}$}
\put(315, 30){\circle{10}}\put(310,39){$\D_2^{a_2}$}
\end{picture}}}
\endequat

\noindent In other words, they are of type $\Dti_6$ and $\Eti_6$ respectively.

Let us now study the fibers of $\pht$ at $[\a^3 : \b^2]$ and $[\a_\pm^3 : \b_\pm^2]$, 
where $\a$, $\b$, $\a_\pm$ and $\b_\pm$ are defined in Remark~\ref{rem:eq-g29}. 
This amounts to understand the fibers of $\varphi$ passing through $a_1$, $a_2^+$ and 
$a_2^-$. Let us first determine their irreducible components (we treat only 
the cases of $a_1$ and $a_2^+$, as the case of $a_2^-$ is isomorphic 
to the case of $a_2^+$). Note that
$$\frac{\b^2}{\a^3}=\frac{\b^5}{(\a\b)^3}=\frac{4}{135}\qquad\text{and}\qquad 
\frac{\b_+^2}{\a_+^3}=\frac{\b_+^5}{(\a_+\b_+)^3}=\frac{-36+16\sqrt{6}}{45}.$$
Working inside the affine chart $X^{29}_{(4)}$, 
a {\sc Magma} computation shows that 
$a_1$ (resp. $a_2^+$) is an $A_1$-singularity of $X^{29}_{135,4}$ (resp. $X^{29}_{45,-36+16\sqrt{6}}$). 
In particular, the projective tangent cone of $X^{29}_{135,4}$ (resp. $X^{29}_{45,-36+16\sqrt{6}}$) at $a_1$ 
(resp. $a_2^+$) consists in two points, 
so $\Xhat_{135,4}^{29}$ (resp. $\Xhat_{45,-36+16\sqrt{6}}^{29}$) meets $\D_1^{a_1}$ 
(resp. $\D_1^{a_2^+} \cup \D_2^{a_2^+}$) in two points. Moreover, 
again by using {\sc Magma} computations, we can check that 
$\Xti_{135,4}^{29}$ and $\Xti_{45,-36+16\sqrt{6}}^{29}$ are irreducible of genus $0$. 
This shows that $\pht^{-1}([135:4])=\Xti_{135,4}^{29} \cup \D_1^{a_1}$ (resp. 
$\pht^{-1}([45 : -36+16\sqrt{6}])=\Xti_{45,-36+16\sqrt{6}}^{29} \cup \D_1^{a_2^+} \cup \D_2^{a_2^+}$) is 
a singular fiber of type $\Ati_1$ (resp. $\Ati_2$). 

So we have found that the elliptic fibration $\pht$ has at least $5$ singular fibers 
of respective types $\Ati_1$, $\Ati_2$, $\Ati_2$, $\Dti_6$ and $\Eti_6$. Since 
the sum of the Euler characteristic of these singular fibers is equal to $24$, the elliptic 
fibration $\pht$ has no more singular fiber.
\end{proof}

\bigskip

\subsection{Transcendental lattice}
We aim to prove that $\Pic(\Xti^{29})$ is generated by the classes of the 21 smooth 
rational curves described in the previous subsection. The intersection numbers between 
these 21 smooth rational curves have been determined in the proof of 
Corollary~\ref{coro:singular-fibers-g29} (see~(\ref{eq:intersection-g29-1}) and 
(\ref{eq:intersection-g29-2})). They are gathered in the following proposition.

\bigskip

\begin{pro}\label{prop:intersection-g29}
The intersection graph of the above 21 smooth rational curves is given by

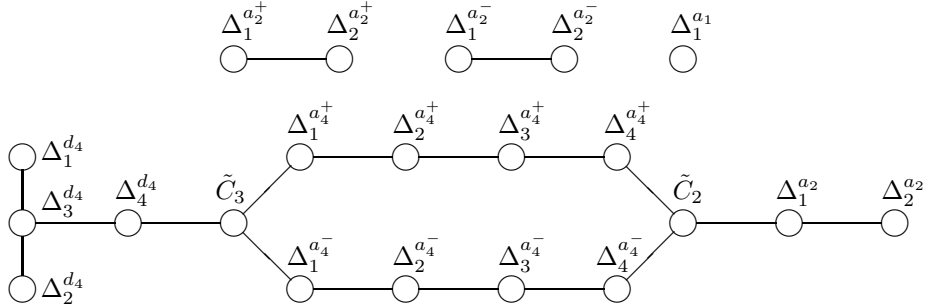
\centerline{\begin{picture}(335,120)
\put( 85, 92){\circle{10}}\put(80,102){$\D_1^{a_2^+}$}
\put( 90, 92){\line(1,0){30}}
\put(125, 92){\circle{10}}\put(120,102){$\D_2^{a_2^+}$}
\put(170, 92){\circle{10}}\put(165,102){$\D_1^{a_2^-}$}
\put(175, 92){\line(1,0){30}}
\put(210, 92){\circle{10}}\put(205,102){$\D_2^{a_2^-}$}
\put(255, 92){\circle{10}}\put(250,102){$\D_1^{a_1}$}
\put(  5, 30){\circle{10}}\put(12,35){$\D_3^{d_4}$}
\put( 10, 30){\line(1,0){30}}
\put( 45, 30){\circle{10}}\put(40,39){$\D_4^{d_4}$}
\put( 50, 30){\line(1,0){30}}
\put(  5, 55){\circle{10}}\put(12,2){$\D_2^{d_4}$}
\put(  5,  5){\circle{10}}\put(12,54){$\D_1^{d_4}$}
\put(  5, 35){\line(0,1){15}}
\put(  5, 25){\line(0,-1){15}}
\put( 85, 30){\circle{10}}\put(78,39){$\Cti_3$}
\put( 88.5355,33.5355){\line(1,1){18.1}}
\put( 88.5355,26.4645){\line(1,-1){18.1}}
\put(110, 55){\circle{10}}\put(115, 55){\line(1,0){30}}\put(105,64){$\D_1^{a_4^+}$}
\put(110,  5){\circle{10}}\put(115,  5){\line(1,0){30}}\put(105,14){$\D_1^{a_4^-}$}
\put(150, 55){\circle{10}}\put(155, 55){\line(1,0){30}}\put(145,64){$\D_2^{a_4^+}$}
\put(150,  5){\circle{10}}\put(155,  5){\line(1,0){30}}\put(145,14){$\D_2^{a_4^-}$}
\put(190, 55){\circle{10}}\put(195, 55){\line(1,0){30}}\put(185,64){$\D_3^{a_4^+}$}
\put(190,  5){\circle{10}}\put(195,  5){\line(1,0){30}}\put(185,14){$\D_3^{a_4^-}$}
\put(230, 55){\circle{10}}\put(225,64){$\D_4^{a_4^+}$}
\put(230,  5){\circle{10}}\put(222,14){$\D_4^{a_4^-}$}
\put(233.5355,51.4645){\line(1,-1){18.1}}
\put(233.5355,8.5355){\line(1,1){18.1}}
\put(255, 30){\circle{10}}\put(251,39){$\Cti_2$}
\put(260, 30){\line(1,0){30}}
\put(295, 30){\circle{10}}\put(300, 30){\line(1,0){30}}\put(290,39){$\D_1^{a_2}$}
\put(335, 30){\circle{10}}\put(330,39){$\D_2^{a_2}$}
\end{picture}}
\end{pro}

\bigskip

We can then compute the lattices $\Pic(\Xti^{29})$ and $\Tb_{\Xti^{29}}$

\bigskip

\begin{theorem}\label{theo:picard-x29}
The Picard group $\Pic(\Xti^{29})$ admits 
\begin{multline*}
([\D_1^{a_2^+}],[\D_2^{a_2^+}],[\D_1^{a_2^-}],[\D_2^{a_2^-}],[\D_1^{a_1}],
[\D_1^{d_4}],[\D_2^{d_4}],[\D_3^{d_4}],[\D_4^{d_4}],[\Cti_3], \\ 
[\D_1^{a_4^+}],[\D_2^{a_4^+}],
[\D_3^{a_4^+}],[\D_4^{a_4^+}],[\D_1^{a_4^-}],[\D_2^{a_4^-}],[\D_3^{a_4^-}],[\D_4^{a_4^-}],
[\Cti_2],[\D_1^{a_2}])
\end{multline*}
as a $\ZM$-basis. The transcendental lattice of $\Xti^{29}$ is given by 
$$\Tb_{\Xti^{29}}=\begin{pmatrix} 6 & 0 \\ 0 & 60 \end{pmatrix}.$$
\end{theorem}

\bigskip

\begin{proof}
Let us denote by $(D_1,D_2,\dots,D_{20})$ the elements written in the statement 
of the theorem, in the same order. Let $I^\circ=(D_j\cdot D_k)_{1 \le j,k \le 20}$. Then 
$\det(I^\circ)=-360$. 
This shows that the family $(D_k)_{1 \le k \le 20}$ is $\ZM$-free 
and, as $\rhob(\Xti^{29})=20$ by Corollary~\ref{cor:q4-20}, this shows that 
$(D_k)_{1 \le k \le 20}$ is a $\QM$-basis of $ \Pic(\Xti^{29})\otimes \QM$. 
We denote by $\L$ the sublattice of $\Pic(\Xti^{29})$ generated by $(D_k)_{1 \le k \le 20}$. 
Its dual lattice $\L^\vee$ in $\Pic(\Xti^{29})\otimes \QM$ 
satisfies $|\L^\vee/\L|=\det(I^\circ)=360$ 
and $\L \subset \Pic(\Xti^{29}) \subset \L^\vee$. 

Let $m$ denote the order of $\Pic(\Xti^{29})/\L$. We must show that $m=1$. 
Assume that there exists a prime number $p$ dividing $m$. 
Then $m^2$ divides $|\L^\vee/\L|$, so $p\in \{2,3\}$. 

Assume first that $2$ divides $m$. Then $\Pic(\Xti^{29})/\L$ contains an element 
of order $2$ and a computation with {\sc Magma} shows that 
this implies that $\Pic(\Xti^{29})$ contains one of the elements 
$$\frac{1}{2}D_5=\frac{1}{2}[\D_1^{a_1}],\quad \frac{1}{2} (D_6+D_7)
=\frac{1}{2}([\D_1^{d_4}] + [\D_2^{d_4}])\quad \text{or}\quad
\frac{1}{2} (D_5+D_6+D_7).$$
But any element $D$ in this list satisfies $D \cdot D \not\in 2\ZM$: 
this contradicts the fact that $\Pic(\Xti^{29})$ is an even lattice.
So $m$ is not divisible by $2$.

Assume finally that $3$ divides $m$. Then $\Pic(\Xti^{29})/\L$ contains an element 
of order $3$ and a computation with {\sc Magma} shows that 
this implies that $\Pic(\Xti^{29})$ contains one of the elements 
$$L_{a,b}=\frac{a}{3} (D_1-D_2)+\frac{b}{3}(D_3-D_4)=
\frac{a}{3} ([\D_1^{a_2^+}]-[\D_2^{a_2^+}])+\frac{b}{3}([\D_1^{a_2^-}]-[\D_2^{a_2^-}])$$
for some $a$, $b \in \{0,1,2\}$ and $(a,b) \neq (0,0)$. 
But $L_{a,b}\cdot L_{a,b} = 2/3(a^2+b^2) \not\in \ZM$, 
so we also get a contradiction. This shows that $m$ is not divisible 
by $3$. Consequently, $m=1$, as expected.

\medskip

Let us now turn to the computation of the transcendental lattice of $\Xti^{29}$. 
First, as there is a finite rational map $X_\mukai \dashrightarrow \Xti^{29}$, 
the transcendental lattice of $\Xti^{29}$ is proportional (by some rational number) to the one 
of $X_\mukai$ by~\cite[Proposition~1.1]{inose}. But the transcendental lattice 
of $X_\mukai$ is given by
$$\Tb_{X_\mukai}=\begin{pmatrix} 4 & 0 \\ 0 & 40 \end{pmatrix}$$
(see for instance~\cite[Proposition~4.4(1)]{bonnafe-sarti-m20}). 
As the discriminant of $\Tb_{\Xti^{29}}$ is equal to the discriminant of $\Pic(\Xti^{29})$, 
this shows that ${\mathrm{disc}}(\Tb_{\Xti^{29}})=360$, and so the only possibility is 
$$\Tb_{\Xti^{29}}=\begin{pmatrix} 6 & 0 \\ 0 & 60 \end{pmatrix},$$
as expected.
\end{proof}

\bigskip

\begin{remark}\label{rem:x29-d2}
Note that one can write 
\begin{multline*}
[\D_2^{a_2}]=
[\D_1^{d_4}]+[\D_2^{d_4}]+2 [\D_3^{d_4}] + 2 [\D_4^{d_4}] + 2 [\Cti_3] + [\D_1^{a_4^+}] 
+ [\D_1^{a_4^-}] \\ - [\D_3^{a_4^+}] - [\D_3^{a_4^-}] - 2 [\D_3^{a_4^+}] - 2 [\D_3^{a_4^-}] 
- 3 [\Cti_2] - 2 [\D_1^{a_2}].~{\SS{\blacksquare}}
\end{multline*}
\end{remark}

We conclude this section by determining the Mordell-Weil group of $\pht$, 
with respect to the section $\thet^+$:

\bigskip

\begin{pro}\label{prop:mordell-g29}
$\mweil_{\thet^+}(\pht) = \ZM [\D_2^{a_4^-}] \simeq \ZM$.
\end{pro}

\bigskip

\begin{proof}
First, it follows from~\cite[Nr. 2493]{shimada} that the torsion group of 
$\mweil_{\thet^+}(\pht)$ is trivial. By the description of the singular fibers 
of the fibration $\pht$ given in Corollary~\ref{coro:singular-fibers-g29}, 
the rank of the group $\triv_{\thet^+}(\pht)$ is equal to $19$. 
Hence $\mweil_{\thet^+}(\pht) \simeq \ZM$. To determine the generators, 
one just needs to notice that $\Pic(\Xti^{29})$ is generated by all 
the classes given in Theorem~\ref{theo:picard-x29} while $\triv_{\thet^+}(\pht)$ 
is generated by all these classes except $[\D_2^{a_4^-}]$ (see Remark~\ref{rem:x29-d2}) and we fix $\D_2^{a_4^+}$ 
as the zero section of the fibration. 
\end{proof}

\bigskip

\bigskip

\subsection{Complements: conics in ${\boldsymbol{X_\mukai}}$} 
As explained in~\cite[Proposition~4.3]{bonnafe-sarti-m20}, 
the K3 surface $X_\mukai$ is the Kummer surface of the abelian surface 
$\EC_{i\sqrt{10}} \times \EC_{i\sqrt{10}}$, where $\EC_\a$ denotes the elliptic 
curve $\CM/(\ZM \oplus \ZM \a)$. Therefore, there exists a Nikulin configuration 
in $X_\mukai$ (i.e., 16 two by two disjoint smooth rational curves). Since~\cite{bonnafe-sarti-m20} 
appeared, it has been shown by Degtyarev~\cite[Theorem~1.1~and~Introduction]{degtyarev} 
that $X_\mukai$ contains $800$ irreducible conics 
(note that 320 conics were already found in~\cite[Remark~4.4]{bonnafe-sarti-m20} but this set 
of conics contains no Nikulin configuration). Later, Naskr\k{e}cki found explicit equations 
for the 800 conics, and showed that one can extract from this set a Nikulin configuration, \cite{Nas}. Let us describe them here. For this, let 
$$C_0=\{[x:y:z:t] \in \PM^3(\CM)~|~
z+i\frac{1+\sqrt{5}}{2}t=x^2+2\sqrt{2}xy+y^2+3\frac{1+\sqrt{5}}{2}t^2=0\},$$
$$C_1=\{[x:y:z:t] \in \PM^3(\CM)~|~x+y+z=y^2 + yz + z^2 + 
\frac{3+\sqrt{10}}{2} t^2=0\}$$
$$C_2=\{[x:y:z:t] \in \PM^3(\CM)~|~x+y+z=y^2 + yz + z^2 + 
\frac{3-\sqrt{10}}{2} t^2=0\}.
\leqno{\text{and}}$$
Then $C_0$, $C_1$ and $C_2$ are conics contained in $X_\mukai$ and belonging to different 
$G_{29}$-orbits. Moreover, the $G_{29}$-orbit 
of $C_0$ (resp. $C_1$, resp. $C_2$) has cardinality 480 (resp. 160, resp. 160).

\bigskip

\section{The group $G_{30}=\Wrm(H_4)$}\label{sec:h4}

\medskip

\boitegrise{{\bf Hypothesis.} {\it We assume in this section, 
and only in this section, that $W=G_{30}=\Wrm(H_4)$.}}{0.75\textwidth}

\medskip

Recall that $G_{30}$ is the Coxeter group $\Wrm(H_4)$ of type $H_4$. 
In other words, we have $G_{30}=\langle s_1,s_2,s_3,s_4\rangle$ in its natural 
representation of dimension $4$ associated with 
the Coxeter graph of type $H_4$, i.e. given by

\centerline{\begin{picture}(150,30)
\put( 15, 10){\circle{10}}\put(11,20){$s_1$}
\put( 20, 10){\line(1,0){30}}\put( 32, 13){$5$ }
\put( 55, 10){\circle{10}}\put(51,20){$s_2$}
\put( 60, 10){\line(1,0){30}}
\put( 95, 10){\circle{10}}\put(91,20){$s_3$}
\put(100, 10){\line(1,0){30}}
\put(135, 10){\circle{10}}\put(131,20){$s_4$}
\end{picture}}

\noindent (see~\cite[Chapter~IV]{bourbaki} for the definition of a 
Coxeter graph and~\cite[Chapter~V,~\S{5}]{bourbaki} for the definition of its 
associated representation). Explicit matrices may be found in~\cite{bonnafemagma}. 
We refer to~\cite[Table~\ref{table:degres}]{bonnafe-sarti-1} 
for the numerical 
informations used here. First, recall that $G_{30}'=G_{30}^\slv$ 
and that $G_{30}/G_{30}' \simeq \mub_2$. 
As the group is a Coxeter group, there exists a real vector subspace $V_\RM$ 
of $V$ such that $V=\CM \otimes_\RM V_\RM$ and 
which is stabilized by $G_{30}$. This also implies that 
$G_{30}$ admits an invariant $f_1$ of degree $2$, which is 
the scalar extension of a positive definite quadratic form on $V_\RM$. 
We fix a fundamental invariant $f_2$ of degree $12$.
If $\l \in \CM$, we set $f_{2,\l}=f_2+\l f_1^6$: this describes (up to scalar) 
all the fundamental invariants of degree $12$. We set
$$X_\l^{30}=\ZC(f_{2,\l})/G_{30}'$$
We proved in~\cite[Theorem~\ref{theo:main}]{bonnafe-sarti-1} that $X_\l^{30}$ 
is a K3 surface with ADE singularities (retrieving a result 
of Barth and the second author~\cite{barth-sarti}). Let $\pi_\l : \Xti_\l^{30} \to X_\l^{30}$ 
denote its minimal resolution: it is a smooth K3 surface. 
As this example was already studied in~\cite{barth-sarti}, we will not 
compute again the singularities of $X_\l^{30}$ 
as well as the transcendental lattices given in Table~\ref{table:k3}. 
We will just give some complementary informations coming from the 
general theory of complex reflection group (equations, base locus, ramification) 
as well as a description of an elliptic fibration together 
with its singular fibers in most cases.

\bigskip

\subsection{Singular dodecics}\label{sub:h4-singulier}
If $1 \le k \le 4$, 
we denote by $W_k$ the subgroup of $G_{30}$ generated by 
$\{s_1,s_2,s_3,s_4\} \setminus \{s_k\}$. Then
$$W_1 \simeq \SG_4,\qquad W_2 \simeq \langle s_1 \rangle \times \SG_3,$$
$$W_3 \simeq \Wrm(I_2(5)) \times \langle s_4 \rangle\qquad\text{and}\qquad 
W_4 \simeq \Wrm(H_3).$$
Here, $I_2(5)$ (resp. $H_3$) denotes the complete subgraph of $H_4$ 
whose vertices are $s_1$ and $s_2$ (resp. $s_1$, $s_2$ and $s_3$) 
and $W(I_2(5))=\langle s_1, s_2 \rangle$ (resp. $W(H_3)=\langle s_1,s_2,s_3 \rangle$) 
is its associated Coxeter group. Note that 
$\Wrm(I_2(5))$ is the dihedral group of order $10$. 
Each maximal parabolic subgroup 
is conjugate to one of the $W_k$'s, and only to one of them because 
they are two by two non-isomorphic. 
Let $v_k \in V_\RM \setminus \{0\}$ 
be such that $V^{W_k}=[v_k]$. We denote by $\O_k$ the $W$-orbit of $[v_k]$ 
in $\PM(V)$. Since $-\Id_V \in W$ 
by~\cite[Table~\ref{table:degres}]{bonnafe-sarti-1} and 
$\Nrm_{G_{30}}(W_k)/W_k$ acts faithfully on $V_\RM^{W_k}=\RM v_k$ (which is of dimension 
$1$), it follows that 
$$\Nrm_{G_{30}}(W_k) = W_k \times \langle -\Id_V \rangle.$$
Since $\Nrm_W(W_k)=W_{[v_k]}$ 
by~\cite[Remark~\ref{rem:maximaux-singuliers}]{bonnafe-sarti-1}, 
we get 
\equat\label{eq:orbites}
|\O_k|=
\begin{cases}
300 & \text{if $k=1$,}\\
600 & \text{if $k=2$,}\\
360 & \text{if $k=3$},\\
60  & \text{if $k=4$.}
\end{cases}
\endequat
Now, $f_1(v_k) \neq 0$ because $f_1$ is positive definite and $v_k \in V_{\RM}$, 
and we can define $\l_k=-f_2(v_k)/f_1(v_k)^6$. 
Therefore,~\cite[Corollary~\ref{coro:v-fixe-singulier}]{bonnafe-sarti-1} shows that 
\equat\label{eq:sarti-surfaces}
\text{\it The singular locus of the 
surface $\ZC(f_{2,\l_k})$ contains $\O_k$.}
\endequat
An explicit computation shows that 
$\l_k \neq \l_l$ if $k \neq l$. So this example explains by 
general theory and simple counting arguments the construction of 
the four singular dodecics constructed by the second author~\cite{sartipencil}. 
It also explains why the singular points are real. However, it does 
not explain why there is no more singular point, why they are 
all nodes, or why there is no more value of $\l$ such that 
$\ZC(f_{2,\l})$ is singular. All these later facts were explained 
in~\cite{sartipencil}.

As a consequence of the above discussion, we get:

\bigskip

\begin{lemma}\label{lem:g30-sing}
If $v \in V \setminus \{0\}$ is such that $[v]$ is a singular point 
of $\ZC(f_{2,\l})$ for some $\l \in \CM$, then $W_v$ 
is a maximal parabolic subgroup of $W$ (in particular, $W_v \neq 1$). 
\end{lemma}

\bigskip

\subsection{Equations}\label{sub:equations-g30}
It follows from~\cite[Proposition~\ref{prop:zf-wsl}]{bonnafe-sarti-1} that 
$$X_\l^{30}=\{[x_1:x_3:x_4:j] \in \PM(2,20,30,60)~|~
j^2=P_\fb(x_1,-\l x_1^6,x_3,x_4)\}.$$
But $\PM(2,20,30,60)=\PM(1,10,15,30)=\PM(1,2,3,6)$. 
Through this sequence of isomorphisms, there exists a polynomial 
$r_\l$ in variables $y_1$, $y_3$, $y_4$ 
which is homogeneous of degree $12$ 
if we assign to $y_1$, $y_3$, $y_4$ the weights $1$, $2$, $3$ 
respectively, and such that 
$P_{\fb}(x_1,-\l x_1^6,x_3,x_4)=
r_\l(x_1^5,x_3,x_4)$. Therefore, 
\equat\label{eq:zf2-h4}
X_\l^{30}=\{[y_1:y_3:y_4:j] \in \PM(1,2,3,6)~|~
j^2=r_\l(y_1,y_3,y_4)\}.
\endequat
We denote by $\s$ the unique non-trivial element of 
$G_{30}/G_{30}' \simeq \mub_2$: through the model of $X_\l^{30}$ 
given by~(\ref{eq:zf2-h4}), the action of $\s$ is described by
$$\s([y_1:y_3:y_4:j])=[y_1:y_3:y_4:-j].$$
Note moreover that 
\equat\label{eq:zh2-h4-total}
\ZC(f_{2,\l})/G_{30} = X_\l^{30}/\mub_2 \simeq \PM(2,20,30) \simeq 
\PM(1,2,3)
\endequat
(see~\cite[Proposition~\ref{prop:projectif-a-poids}]{bonnafe-sarti-1}). 
The branch locus $R_\l$ of the quotient morphism $\xi_\l : X_\l^{30} \longto \PM(1,2,3)$ 
is the zero set of $r_\l$.

\bigskip

\subsection{Base locus}
Let $B$ denote the {\it base locus} of the family 
of dodecic surfaces $(\ZC(f_{2,\l}))_{\l \in \CM}$, that is, the subvariety of $\PM(V)$ 
which is contained in all the members of this family. Namely,
$$B=\{p \in \PM(V)~|~f_1(p)=f_2(p)=0\}.$$
Note that $\d(10)=\d^*(10)=2$, so that $\dim V(10)=2$. 
We denote by $L_{10}$ the line $\PM(V(10))$ in $\PM(V)$. 
The next result was already obtained by Barth and the second author~\cite{barth-sarti}, 
but we give a proof that makes it an application 
of Lehrer-Springer theory.

\bigskip

\begin{pro}\label{prop:g30-base}
The stabilizer $W(10)$ of $L_{10}$ in $W$ is equal to $C_W(w_{10})$ 
and has order $600$. Moreover, 
$$B=\bigcup_{x \in W} x(L_{10})$$
consists of $24$ lines, which split into two $G_{30}'$-orbits 
of cardinality $12$. 
\end{pro}

\bigskip

\begin{proof}
This is mainly a consequence 
of~\cite[Theorem~\ref{theo:springer}]{bonnafe-sarti-1}. 
Indeed, the fact that 
$$B=\bigcup_{x \in W} x(L_{10})$$
follows from~\cite[Theorem~\ref{theo:springer}(d)]{bonnafe-sarti-1}. 
Moreover, by~\cite[Theorem~\ref{theo:springer}(f)]{bonnafe-sarti-1}, 
we have that $W(10)=C_W(w_{10})$ is a reflection group for its action 
on $V(10)$, and admits $(20,30)$ as list of degrees. 
So $|W(10)|=20 \cdot 30 = 600$ 
by~\cite[(\ref{eq:degres})]{bonnafe-sarti-1}. 

The only fact that is not covered 
by~\cite[Theorem~\ref{theo:springer}]{bonnafe-sarti-1} 
is that the $24$ lines forming $B$ split into 
two $G_{30}'$-orbits of cardinality $12$: but this follows 
from the fact that $W(10) \subset G_{30}'$ (which 
can be checked for instance with {\sc Magma}).
\end{proof}

\bigskip

Let $B'$ denote the image of $B$ in $\PM(V)/W'$. Then it follows 
from Proposition~\ref{prop:g30-base} that $B'$ is the union 
of two irreducible components $B^+$ and $B^-$. We denote by 
$\Bti^+$ and $\Bti^-$ their respective strict transforms in 
$\Xti_\l^{31}$. 

Let us examine some particular points of $B$. First, note that 
$B$ does not contain a singular point of $\ZC(f_{2,\l})$ since 
we have seen in~\S\ref{sub:h4-singulier} that $f_1(v) \neq 0$ 
for any $v \in V \setminus \{0\}$ such that $[v]$ is a singular 
point of $\ZC(f_{2,\l})$. 

Now, let $k \in \{20,30\}$. Examining Table~\ref{table:degres-29-30-31}, 
we see that $\d(k)=\d^*(k)=1$. 
By Springer Theory~\cite[Theorem~\ref{theo:springer}]{bonnafe-sarti-1}, 
this implies that $\dim V(k)=1$, that $W(V(k))=1$ and 
that $W(k)=\langle w_k \rangle$. Let $z_k$ denote the image of $V(k)$ 
in $\PM(V)$. Then the stabilizer of $z_k$ in $W$ is $W(k)$ and 
since $\det w_k=\z_k^{-60}=1$ 
(see \cite[Theorem~\ref{theo:springer}$($f$)$]{bonnafe-sarti-1}), 
this implies that the $W$-orbit $\O_k$ of $z_k$ has cardinality $14400/k$ 
and splits into two $W'$-orbits. We denote by $a_{(k/10)-1}$ 
the image of $z_k$ in $\ZC(f_{2,\l})/W \simeq \PM(1,2,3)$: 
it follows from~\cite[Theorem~\ref{theo:springer}$($d$)$]{bonnafe-sarti-1} 
that 
$$a_1=[0:1:0] \qquad \text{and}\qquad a_2=[0:0:1].$$
Note that $a_r$ is an $A_r$ singularity 
of $\PM(1,2,3)$. Now, the morphism $X_\l^{30} \to \PM(1,2,3)$ is unramified 
above $a_r$ because $W(20)$ and $W(30)$ are contained in $W'$. 
So let $a_r^\pm$ denote the two points of $X_\l^{30}$ above $a_r$: 
in the model given in~\S\ref{sub:equations-g30}, we have 
$$a_1^\pm=[0:1:0:\pm j_1]\qquad\text{and}\qquad a_2^\pm=[0:0:1:\pm j_2]$$
for some $j_r \in \CM^\times$. 
They are both $A_r$ singularities of $X_\l^{30}$ (note that this is true 
for any value of $\l$). We choose the value of $j_r$ so that 
$a_r^+ \in B^+$ (and then $a_r^- \in B^-$). Recall 
from~\cite{barth-sarti} that
\equat\label{eq:zsing-base-g30}
(X_\l^{30})_\sing \cap B' = \{a_1^+,a_1^-,a_2^+,a_2^-\}.
\endequat
Again, this fact holds for any value of $\l$.

\bigskip

\begin{lemma}\label{lem:sing-30}
Let $x \in X_\l^{30} \setminus \{a_1^\pm,a_2^\pm\}$. Then $x$ is singular if and 
only if $\xi_\l(x)$ is a singular point of the branch locus $R_\l$. In this case, 
the singularity $x$ is of the same type as the singularity $\xi_\l(x)$ of the curve $R_\l$. 
\end{lemma}

\bigskip

\begin{proof}
Since the only singular points of $\PM(1,2,3)$ are $a_1$ and $a_2$, 
the result follows from~\cite[Proposition~\ref{prop:double cover}]{bonnafe-sarti-1}. 
\end{proof}

\bigskip

\subsection{Elliptic fibration}
With the model of $X_\l^{30}$ given in~\S\ref{sub:equations-g30}, we can 
define a map 
$$\fonction{\ph_\l}{X_\l^{30} \setminus 
\{a_2^+,a_2^-\}}{\PM^1(\CM)}{[y_1:y_3:y_4:j]}{[y_1^2:y_3].}$$
Since this map factorizes through the quotient $\PM(1,2,3)$ 
of $X_\l^{30}$, the same argument as in the proof of 
Proposition~\ref{pro:elliptic-g29} shows that:

\bigskip

\begin{pro}\label{pro:elliptic-g30}
The map $\ph_\l \circ \pi_\l : 
\Xti_\l^{30} \setminus \bigl(\pi_\l^{-1}(a_2^+) \cup \pi_\l^{-1}(a_2^-)\bigr) 
\longto \PM^1(\CM)$ extends to a morphism of algebraic varieties 
$$\pht_\l : \Xti_\l^{30} \longto \PM^1(\CM).$$
\end{pro}

\bigskip

\begin{remark}\label{rem:g30-section}
By the same argument as in Remark~\ref{rem:section-g29}, the elliptic 
fibration $\pht_\l  : \Xti_\l^{30} \longto \PM^1(\CM)$ admits two sections 
$\thet_\l^\pm : \PM^1(\CM) \longto \Xti_\l^{30}$ which satisfy $\thet^-=\s \circ \thet^+$.\finl
\end{remark}

\bigskip

Note that the above result is independent of $\l$. However, 
we will see in the next corollary that the singular fibers 
of the elliptic fibration $\pht_\l$ depend on $\l$. We will not determine 
the fiber in all cases, but only whenever the following hypothesis is satisfied:

\medskip

\begin{quotation}
{\bf Hypothesis ($\Hb_\l$).} {\it If $x$ and $y$ are two different singular points 
of $X_\l^{30} \setminus \{a_1^\pm,a_2^\pm\}$, then $\ph_\l(x) \neq \ph_\l(y)$.}
\end{quotation}

\medskip

\noindent Note that Hypothesis ($\Hrm_\l$) holds for all but a finite number of values of $\l$. 
Moreover, an explicit computation with {\sc Magma} shows that it holds 
for $\l \in \{\l_1,\l_2,\l_3,\l_4\}$. 

\bigskip

\begin{cor}\label{coro:singular-g30}
Let $\l \in \CM$ be such that $(\Hrm_\l)$ holds. Then 
the singular fibers of the elliptic fibration $\pht_\l : \Xti_\l^{30} \to \PM^1(\CM)$ 
are given by Table~\ref{table:k3-g30}.
\end{cor}

\bigskip

\begin{proof}
Let us first examine the fiber at $[0:1]$. For this particular fiber, 
the description will not depend on $\l$. Note that 
$$\overline{\ph_\l^{-1}([0:1])}=B'=B^+ \cup B^-.$$
We now apply results from Appendix~\ref{appendix-2} in the case where $(k,l)=(1,2)$.  
Let $\D_1$ and $\D_2$ denote the lines in $\hat{\PM}(1,2,3)$ described in the 
Appendix~\ref{appendix-2} and let $\phh_{1,2} : \hat{\PM}(1,2,3) \to \PM^1(\CM)$ 
denote the map constructed in~\eqref{eq:map-p1}. It follows from 
Proposition~\ref{prop:graphe-eclatement} that 
$$\phh_{1,2}^{-1}([0:1])=\D_2 \cup \Delt^{(1)},\leqno{(\clubsuit)}$$
where $\D^{(1)}=\{[y_1:y_3:y_4] \in \PM(1,2,3)~|~y_1=0\}$ 
and $\Delt^{(1)}$ is the strict transform of $\D^{(1)}$ in $\hat{\PM}(1,2,3)$. 
By the argument in Remark~\ref{rem:section-g29}, the two smooth rational 
curves $\D_1^{a_2^\pm}$ and $\D_2^{a_2^\pm}$ above the point $a_2^\pm$ 
can be numbered so that $\D_k^{a_2^\pm}$ is mapped isomorphically to $\D_k$ 
through the quotient morphism $\Xti_\l^{30} \to \PM(1,2,3)$, 
moreover, the inverse image of $\D^{(1)}$ in $X_\l$ is $B^+ \cup B^-$. 
So, if we denote 
by $\D^{a_1^\pm}$ the smooth rational curve above the points $a_1^\pm$ 
and by $\Bti$ (resp. $\Bti^\pm$) the strict transform of $B'$ (resp. $B^\pm$) 
in $\Xti_\l^{30}$, then it follows from~$(\clubsuit)$ and the construction of 
$\pht_\l$ that
$$\pht_\l^{-1}([0:1])=\D^{a_1^+} \cup \D^{a_1^-} \cup \D_2^{a_2^+} \cup \D_2^{a_2^-} 
\cup \Bti^+ \cup \Bti^-.$$
Since $\Bti^+ \cap \Bti^- \neq \vide$, since $\D^{a_1^\e} \cap B^\eta \neq \vide$ 
if and only if $\e=\eta$, since $\D_2^{a_2^\e} \cap B^\eta \neq \vide$ 
if and only if $\e=\eta$ and since $\D^{a_1^\e} \cap \D_2^{a_2^\eta} = \vide$, 
the Kodaira-N\'eron classification of singular fibers implies that 
$$\text{$\pht_\l^{-1}([0:1])$ is of type $\Dti_5$.}\leqno{(\diamondsuit)}$$
Note that $(\diamondsuit)$ holds for any value of $\l$.
We will now start the discussion according to the value of $\l$.

\medskip

Assume that $\l \not\in \{\l_1,\l_2,\l_3,\l_4\}$. 
Then $X_\l^{30} \setminus \{a_1^\pm,a_2^\pm\}$ has 6 singular points $x_1$,\dots, $x_6$, 
of respective type $A_1$, $A_1$, $A_1$, $A_2$, $A_2$ and $A_4$. Let $x$ be one of these 
6 points and let $m$ denote its Milnor number. Then $\pi_\l^{-1}(x)$ is the union of $m$ 
smooth rational curves $\D_1^x$,\dots, $\D_m^x$. 
Let $E^x$ denote the closure of $\ph_\l^{-1}(\ph_\l(x))$ and let $\Eti^x$ denote its 
strict transform in $\Xti_\l^{30}$. Then 
$$\pht_\l^{-1}(\ph_\l(x)) = \Eti^x \cup \D_1^x \cup \cdots \cup \D_m^x.\leqno{(\heartsuit)}$$
So $\pht_\l^{-1}(\ph_\l(x))$ is a singular fiber. Let us determine its type. 

Note that $\phh_{1,2}^{-1}(\ph_\l(x))$ is a projective line by Proposition~\ref{prop:graphe-eclatement} 
and Remark~\ref{rem:p1-fibre-generale}. Therefore, its double cover $\Eti^x$ has at 
most two irreducible components. Note also that the multiplicity of $\D_k^m$ 
in the singular fiber $\pht_\l^{-1}(\ph_\l(x))$ is equal to one. Therefore, 
according to the Kodaira-N\'eron classification of singular fibers, 
$(\heartsuit)$ gives the following possibilities:
\begin{itemize}
\item If $x=x_1$, $x_2$ or $x_3$ is an $A_1$ singularity, then $\pht_\l^{-1}(\ph_\l(x))$ is of type 
$\Ati_1$ or $\Ati_1^\#$ if $\Eti^x$ is irreducible or of type 
$\Ati_2$ or $\Ati_2^\#$ if $\Eti^x$ has two irreducible components.

\item If $x=x_4$ or $x_5$ is of type $A_2$, then $\pht_\l^{-1}(\ph_\l(x))$ is of type $\Ati_2$ 
or $\Ati_2^\#$ if $\Eti^x$ is irreducible or of type $\Ati_3$ if $\Eti^x$ 
has two irreducible components.

\item If $x=x_6$ is of type $A_4$, then $\pht_\l^{-1}(\ph_\l(x))$ is of type $\Ati_4$ 
if $\Eti^x$ is irreducible or of type $\Ati_5$ if $\Eti^x$ 
has two irreducible components.
\end{itemize}
Let $\chi_k$ denote the Euler characteristic of the singular fiber above $x_k$. 
Since the Euler characteristic of $\pht_\l^{-1}([0:1])$ is equal to $7$ by~$(\diamondsuit)$, 
we have
$$\chi_1+\chi_2+\chi_3+\chi_4+\chi_5+\chi_6 \le 24-7=17.\leqno{(\spadesuit)}$$
But it follows from the above discussion that 
$$\chi_1 \ge 2,\quad \chi_2 \ge 2,\quad \chi_3 \ge 2,\quad \chi_4 \ge 3,\quad \chi_5 \ge 3
\quad \text{and}\quad \chi_6 \ge 5.$$
Therefore, $(\spadesuit)$ forces $\chi_1=\chi_2=\chi_3=2$, $\chi_4=\chi_5=3$ and $\chi_6=5$. And so 
the singular fibers are of the types described in the first line of Table~\ref{table:k3-g30}.

\medskip

The cases mentioned in the last four lines of Table~\ref{table:k3-g30} follow from a similar discussion, 
the conclusion using the same argument based on the Euler characteristic.
\end{proof}

\bigskip

\begin{table}\label{table3}
$$\begin{array}{!{\vline width 2pt} c !{\vline width 1.3pt} c|c|c|c!{\vline width 2pt}}
\hlinewd{2pt}
\ZC_\sing(f_{2,\l}) & \text{singularities of $X_\l^{30}$} & 
\Tb_{\Xti_\l^{30}} & \text{singular fibers of $\pht_\l$} & \mweil_{\thet_\l^+}(\pht_\l) \petitespace
\\
\hlinewd{2pt}
\vide & A_4+4\,A_2+5\,A_1 & {\rm Theorem}\,\,\ref{teo30}
& \Dti_5 + \Ati_4 + 2\,\Ati_2 + 3\,\Ati_1^{(\dagger)} & \ZM \petitespace \\ \hline
60\,A_1 & E_8+3\,A_2+4\,A_1 &  \begin{pmatrix} 4 & 2 & 34 \end{pmatrix} & \Eti_8+\Dti_5+\Ati_2+2\,\Ati_1 & \ZM \petitespace\\ \hline
300\,A_1 & E_6+A_4+2\,A_2+4\,A_1 & \begin{pmatrix} 12 & 6 & 58 \end{pmatrix} & \Eti_6+\Dti_5+\Ati_4+2\,\Ati_1 & \ZM \petitespace\\ \hline
360\,A_1 & D_7+4\,A_2+3\,A_1 & \begin{pmatrix} 6 & 0 & 132 \end{pmatrix} & \Dti_7+\Dti_5 +2\,\Ati_2+\Ati_1 & \ZM \petitespace\\ \hline
600\,A_1 & D_5+A_4+3\,A_2+3\,A_1 & \begin{pmatrix} 6 & 0 & 220 \end{pmatrix} & 2\,\Dti_5+\Ati_4+\Ati_2+\Ati_1 & \ZM \petitespace\\
\hlinewd{2pt}
\end{array}$$
\caption{Some numerical data for the family of K3 surfaces 
$(\Xti_\l^{30})_{\l \in \CM}$}\label{table:k3-g30}
\vskip-0.4cm
$^{(\dagger)}$ Only for $\l$ generic
\end{table}

\section{The group $G_{31}$}\label{sec:g31}

\medskip

\boitegrise{{\bf Hypothesis.} 
{\it We assume in this section, and only in this 
section, that $W=G_{31}$.}}{0.75\textwidth}

\medskip

Let
$$s_1=
\begin{pmatrix}
0 & 1 & 0 & 0 \\
1 & 0 & 0 & 0 \\
0 & 0 & 1 & 0 \\
0 & 0 & 0 & 1\\
\end{pmatrix},\qquad
s_2=\begin{pmatrix}
1 & 0 & 0 & 0 \\
0 & 0 & 1 & 0 \\
0 & 1 & 0 & 0 \\
0 & 0 & 0 & 1\\
\end{pmatrix},$$
$$s_3=
\begin{pmatrix}
1 & 0 & 0 & 0 \\
0 & 1 & 0 & 0 \\
0 & 0 & -1 & 0 \\
0 & 0 & 0 & 1 \\
\end{pmatrix}\quad\text{and}\quad
s_4=\frac{1}{2}\begin{pmatrix}
1 & -1 & -1 & -1 \\
-1 & 1 & -1 & -1 \\
-1 & -1&  1 & -1 \\
-1 & -1& -1 &  1 \\
\end{pmatrix}.$$
Then $\Wrm(F_4)=G_{28}=\langle s_1,s_2,s_3,s_4 \rangle$. We set
$$s_5=\begin{pmatrix}
0 & i & 0 & 0 \\
-i& 0 & 0 & 0 \\
0 & 0 & 1 & 0 \\
0 & 0 & 0 & 1 \\
\end{pmatrix}.$$
Then $G_{31}=\langle s_1,s_2,s_3,s_4,s_5\rangle$. Note that, even though 
$G_{31}$ is of rank $4$, it cannot be generated by only $4$ reflections.
Note also that $|G_{31}/G_{31}'|=2$ and that $G_{31}'=G_{31}^\slv$. 

Recall from~\cite[Table~\ref{table:degres}]{bonnafe-sarti-1} that 
$\Deg(W)=(8,12,20,24)=(d_1,d_2,d_3,d_4)$ 
and we denote by $\fb=(f_1,f_2,f_3,f_4)$ a family of fundamental 
invariants such that $\deg(f_i)=d_i$. Then $f_1$ and $f_2$ are 
uniquely determined (up to a scalar). We have
$$f_1=\Sigma(x^8)+14\Sigma(x^4y^4)+168 x^2y^2z^2t^2$$
$$f_2=\Sigma(x^{12})-33\Sigma(x^8y^4)+792\Sigma(x^6y^2z^2t^2)
+330 \Sigma(x^4y^4z^4).$$
We will make a special choice 
for $f_3$ as follows. First, let $N$ denote the normalizer of $G_{28}$ 
in $G_{31}$. Then $N$ has index $10$ in $G_{31}$ and we denote by 
$[G_{31}/N]$ a set of representatives of the cosets in $G_{31}/N$. 
Then $x^2+y^2+z^2+t^2$ is $G_{28}$-invariant (but not $N$-invariant) 
and it turns out that 
$$f_3=\prod_{g \in [G_{31}/N]} \lexp{g}{(x^2+y^2+z^2+t^2)}$$
is a fundamental invariant of degree $20$ of $G_{31}$. 
Of course, $\ZC(f_3)$ is not irreducible (it is the union of $10$ 
quadrics). We choose the set of representatives $[G_{31}/N]$ such that 
the coefficient of $x^{14}y^2z^2t^2$ in $f_3$ is equal to $648$. Then
\eqna
f_3&=& 648(\Sigma(x^{14} y^2 z^2 t^2) - \Sigma(x^{12} y^4 z^4) - \Sigma(x^{10} y^6 z^2 t^2) \\
&& + 2 \Sigma(x^8 y^8 z^4) + 13 \Sigma(x^8 y^4 z^4 t^4) - 14\Sigma(x^6 y^6 z^6 t^2)).
\endeqna
Finally, we set
\eqna
f_4&=&3888(\Sigma(x^{18}y^2z^2t^2) + 2\Sigma(x^{16}y^4z^4) 
- 12\Sigma(x^{14}y^6z^2t^2) - 2\Sigma(x^{12}y^8z^4) \\
&& + 76\Sigma(x^{12}y^4z^4t^4) + 22\Sigma(x^{10}y^{10}z^2t^2) 
- 52\Sigma(x^{10}y^6z^6t^2) \\ 
&& + 36\Sigma(x^8y^8z^8) + 36\Sigma(x^8y^8z^4t^4) 
- 8 x^6y^6z^6t^6).
\endeqna
Then $\fb=(f_1,f_2,f_3,f_4)$ is a family of fundamental 
invariants of $G_{31}$. Note that the coefficients $648$ (for $f_3$) and 
$3888$ (for $f_4$) are just for simplifying the general 
equation of the surfaces studied in this section. 

\medskip

If $\l \in \CM$, we set $f_{3,\l}=f_3+\l f_1f_2$. 
Recall from~\cite{bonnafesingular}  
that there are only $6$ values of $\l$ such that $\ZC(f_{3,\l})$ 
is singular: one of them is $0$, which is the only value of $\l$ for which 
$\ZC(f_{3,\l})$ is not irreducible. We set 
$$X_\l^{31}=\ZC(f_{3,\l}) / G_{31}'$$
(it is a K3 surface with ADE singularities 
by~\cite[Theorem~\ref{theo:main}]{bonnafe-sarti-1}) 
and we denote by $\Xti_\l^{31}$ its minimal resolution (it is a smooth 
K3 surface). 
We aim in this section to prove the results stated in 
Table~\ref{table:k3}, namely compute the singularities 
of $X_\l^{31}$, the Picard number 
and the transcendental lattice of $\Xti_\l^{31}$. 
We will also provide some more informations about 
the geometry of $\ZC(f_{3,\l})$ and $X_\l^{31}$ (lines, 
branch locus of the double cover 
$X_\l^{31} \to \ZC(f_{3,\l})/G_{31}=\PM^2(\CM)$,\dots).

\bigskip

\subsection{Equations, branch locus}\label{sub:g31-ram}
Let $\xi_\l : X_\l^{31} \to \ZC(f_{3,\l})/G_{31} = \PM^2(\CM)$ be the 
natural map. This is a double cover, whose branch  
locus $R_\l \subset \PM^2(\CM)$ is a sextic that 
will be described below. 
First~\cite[Proposition~\ref{prop:zf-wsl}]{bonnafe-sarti-1}.
$$X_\l^{31} = \{[x_1:x_2:x_4:j] \in \PM(8,12,24,60)~|~
j^2=P_\fb(x_1,x_2,-\l x_1x_2,x_4)\}.$$
But $\PM(8,12,24,60) \simeq \PM(2,3,6,15) \simeq \PM(2,1,2,5)$. 
So there exists a polynomial $q_\l \in \CM[y_1,y_2,y_4]$ 
which is homogeneous of degree $10$ if we assign 
to $y_1$, $y_2$, $y_4$ the degrees $2$, $1$, $2$ 
respectively, and such that 
$$X_\l^{31} = \{[y_1:y_2:y_4:j] \in \PM(2,1,2,5)~|~
j^2=q_\l(y_1,y_2,y_4)\}.$$
But $\PM(2,1,2) \simeq \PM(1,1,1) =\PM^2(\CM)$, so there 
exists a polynomial $r_\l(z_1,z_2,z_4) \in \CM[z_1,z_2,z_4]$, 
which is homogeneous of degree $5$ if we assign 
to $z_1$, $z_2$, $z_4$ the degrees $1$, $1$, $1$ 
respectively, and such that 
\equat\label{eq:xl-g31}
X_\l^{31} = \{[y_1:y_2:y_4:j] \in \PM(2,1,2,5)~|~
j^2=r_\l(y_1,y_2^2,y_4)\}.
\endequat
Through this description, the action of the 
unique non-trivial element $\s$ of $G_{31}/G_{31}'$ 
is given by
$$\s([y_1:y_2:y_4:j])=[y_1:y_2:y_4:-j] =[y_1:-y_2:y_4:j],$$
and the morphism $X_\l^{31} \longto \PM^2(\CM)$ is given explicitly by
$$[y_1:y_2:y_4:j] \longmapsto [y_1:y_2^2:y_4].$$
So the branch locus of $\xi_\l$ is 
\equat\label{eq:rl-g31}
R_\l=\{[z_1:z_2:z_4] \in \PM^3(\CM)~|~z_2 r_\l(z_1,z_2,z_4)=0\}.
\endequat
In other words, the sextic $R_\l$ is the union of the projective 
line $\Bba_2$ defined by $z_2=0$ and of the quintic $R_\l'=\ZC(r_\l)$. 

\medskip

\subsubsection{Other model} Let 
$$\XC_\l = \{[z_1 : z_2 : z_4:t] \in \PM(1,1,1,3)~|~t^2=z_2r_\l(z_1,z_2,z_4)\}.$$
Then $\XC_\l \to \PM^2(\CM)$, $[z_1 : z_2 : z_4:t] \mapsto [z_1 : z_2 : z_4]$ 
is a double cover of $\PM^2(\CM)$ ramified on the sextic $R_\l=\Bba_2 \cup R_\l'$. 
The rational map
$$\fonction{\iota}{\PM(2,1,2,5)}{\PM(1,1,1,3)}{
[y_1:y_2:y_4:j]_{\SSS{2,1,2,5}}}{[y_1:y_2^2:y_4:y_2j]_{\SSS{1,1,1,3}}}$$
is well-defined outside of $[0:0:0:1]_{2,1,2,5}$ and is birational 
(it is for instance an isomorphism between the open 
subsets defined respectively by $y_2 \neq 0$ and $z_2 \neq 0$). But note that 
$[0:0:0:1]_{2,1,2,5} \not\in X_\l^{31}$ and that $\iota(X_\l^{31})=\XC_\l^{31}$. 
Also $X_\l^{31}$ (resp. $\XC_\l$) is contained in the open subsets defined 
by $(y_i \neq 0)_{i \in \{1,2,4\}}$ (resp. $(z_i \neq 0)_{i \in \{1,2,4\}}$). 
An immediate computation in all these open subsets show that $\iota$ induces 
an isomorphism $X_\l^{31} \longiso \XC_\l$. As this second model is somewhat 
simpler to work with, we will now identify $X_\l^{31}$ with $\XC_\l$ 
and so view $X_\l^{31}$ in the more classical model for double covers of $\PM^2(\CM)$ 
ramified above a sextic:
\equat\label{eq:xl-g31-bis}
X_\l^{31} = \{[z_1 : z_2 : z_4:t] \in \PM(1,1,1,3)~|~t^2=z_2r_\l(z_1,z_2,z_4)\}.
\endequat
Through this model, the double cover morphism $\xi_\l : X_\l^{31} \longto \PM^2(\CM)$ is just given by 
$\xi_\l([z_1 : z_2 : z_4:t])=[z_1 : z_2 : z_4]$.

\bigskip

\subsubsection{Value of $r_\l$}\label{subsub:rl}
The explicit value of the polynomial $r_\l$ is given below 
(recall that it depends on our special choice of the family 
$\fb$ of fundamental invariants and a suitable normalization for $J$):
\eqna
r_\l&\!\!=\!\!&-432\,\l^3(\l+1)\,z_1^3 z_2^2 - 108\,\l^2\,z_1^3 z_2 z_4 \\
&&    + (12500\,\l^6 + 22500\,\l^5 + 10800\,\l^4 + 864\,\l^3)\,z_1^2 z_2^3 \\
&&    + (4125\,\l^4 + 3420\,\l^3 + 216\,\l^2)\,z_1^2 z_2^2 z_4 
      + 222\,\l^2\,z_1^2 z_2 z_4^2 + z_1^2 z_4^3 \\ 
&&    - 432\,\l^3\,z_1 z_2^4 + (900\,\l^3 - 108\,\l^2)\,z_1 z_2^3 z_4 
      + (-500\,\l^3 + 210\,\l^2)\,z_1 z_2^2 z_4^2 \\
&&    + (-150\,\l^2 - 24\,\l - 2)\,z_1 z_2 z_4^3 - 2\,z_1 z_4^4 + z_2^2 z_4^3 - 2\,z_2 z_4^4 + z_4^5.\\
\endeqna

\bigskip

\begin{remark}\label{rem:lambda=0}
Assume in this remark, and only in this remark, that $\l=0$. Then
$$X_0^{31}\!=\!\{[z_1:z_2:z_4:t] \in \PM(1,1,1,3)\,|\, 
t^2=z_4^3(z_1^2+z_2^2+z_4^2-2z_1z_2-2z_1z_4-2z_2z_4)\}.$$
The singular locus is a union of the point $[1:0:1:0]$ and the smooth rational 
curve defined by $z_4=t=0$. So the singular locus has dimension $1$ and the surface 
$X_0^{31}$ will not be considered in this section.\finl
\end{remark}

\bigskip

\boitegrise{{\bf Hypothesis.} {\it From now on, and until the end of this section, 
we assume that $\l \neq 0$.}}{0.75\textwidth}

\bigskip

\subsection{Singular icosics}\label{sub:g31-sing}
As explained in the introduction of this section, 
it follows from~\cite{bonnafesingular} that there are 
$5$ values of $\l \in \CM^\times$ such that $\ZC(f_{3,\l})$ is singular. 
We explain here what are these special values, and how we can recover 
the singularities of $\ZC(f_{3,\l})/G_{31}'$ thanks 
to~\cite[Proposition~\ref{prop:double cover}]{bonnafe-sarti-1} and {\sc Magma} computations. 

First, we set 
$$W_{145}=\langle s_1,s_4,s_5 \rangle,\quad W_{245}=\langle s_2,s_4,s_5 \rangle
\quad\text{and}\quad W_{1234}=\langle s_1,s_2,s_3,s_4 \rangle.$$
Note that these are representatives of conjugacy classes 
of maximal parabolic subgroups of $W$. If $k \in \{145,245,1234\}$, 
we denote by $v_k$ a generator of the line $V^{W_k}$, and we 
set $z_k=[v_k]=V^{W_k} \in \PM(V)=\PM^3(\CM)$. We also set 
$N_k=\Nrm_W(W_k)$ and we denote by $\O_k$ the $W$-orbit of $z_k$: 
it follows from~\cite[Remark~\ref{rem:maximaux-singuliers}]{bonnafe-sarti-1} 
that $|\O_k|=|W|/|N_k|$. Concretely, we have:
\equat\label{eq:omega-g31}
|\O_k|=
\begin{cases}
960 & \text{if $k=145$,}\\
480 & \text{if $k=245$.}\\
60 & \text{if $k=1234$.}\\
\end{cases}
\endequat
A {\sc Magma} computation 
shows that 
\equat\label{eq:f1-f2-lisses}
\text{\text{\it $\ZC(f_1)$ and $\ZC(f_2)$ are smooth.}}
\endequat
In particular, $v_k \not\in \ZC(f_1) \cup \ZC(f_2)$ 
by~\cite[Corollary~\ref{coro:v-fixe-singulier}]{bonnafe-sarti-1}. 
So we can define 
$\l_k=-f_3(v_k)/(f_1f_2)(v_k)$. It turns out that $\l_{1234}=0$, 
so that $f_{3,\l_{1234}}=f_3$ is not irreducible: this case does not 
lead to a K3 surface and will not be studied here. 

Therefore, we have found in this way two values of $\l$, namely 
$\l_{145}$ and $\l_{245}$, such that $\ZC(f_{3,\l})$ is irreducible 
and singular. But there are three more values of $\l$ such that $\ZC(f_{3,\l})$ 
is irreducible and singular~\cite{bonnafesingular}: this shows that, 
by opposition with the cases of $G_{29}$ (in degree $8$) and $G_{30}$ 
(in degree $12$),~\cite[Corollary~\ref{coro:v-fixe-singulier}]{bonnafe-sarti-1} 
is not sufficient to explain all the singular icosics that can 
be constructed from fundamental invariants of $G_{31}$ 
of degree $20$. With our choice of the family $\fb$ of fundamental 
invariants of $G_{31}$, we have 
$$\l_{145}=-\frac{8}{25}\qquad\text{and}\qquad \l_{245}=-\frac{81}{175}.$$
We set
$$\l_1=1,\qquad \l_2=-\frac{1}{3}\qquad\text{and}\qquad 
\l_3=-\frac{1}{2}.$$
Then $\l_{145}$, $\l_{245}$, $\l_1$, $\l_2$, $\l_3$ are the five 
values of $\l$ such that $\ZC(f_{3,\l})$ is irreducible and singular. 
By~\cite[Proposition~3.6~and~Table~4]{bonnafesingular} and the correction statement 
at {\tt https://doi.org/10.1080/10586458.2018.1555778}, 
the singularities of $\ZC(f_{3,\l_k})$ are given by 
\equat\label{eq:g31-icosics}
\ZC_\sing(f_{3,\l_k}) =
\begin{cases}
960\,A_1 & \text{if $k=145$,}\\
480\,A_1 & \text{if $k=245$,}\\
1920\,A_1 & \text{if $k=1$,}\\
1440\,A_2 & \text{if $k=2$,}\\
640\,A_3 & \text{if $k=3$.}\\
\end{cases}
\endequat

\bigskip

\subsection{Springer theorem, base locus} 
Recall from~\cite[Table~\ref{table:degres}]{bonnafe-sarti-1} that
$$
\begin{cases}
\Deg(W)=(8,12,20,24),\\
\Codeg(W)=(0,12,16,28).
\end{cases}$$
The following facts can be deduced immediately from 
this and from~\cite[Theorem~\ref{theo:springer}]{bonnafe-sarti-1}:
\begin{itemize}
\itemth{a} $\d(8)=\d^*(8)=2$, so $\dim V(8)=2$. 
We denote by $L_8$ the line $\PM(V(8)) \subset \PM(V)=\PM^3(\CM)$. 
Then $W(8)=C_W(w_8)$ is a reflection group for its action 
on $V(8)$, and its degrees are $(8,24)$. So $|W(8)|=8\cdot 24=192$ 
by~\cite[(\ref{eq:degres})]{bonnafe-sarti-1}, 
and so the $W$-orbit of $L_8$ contains $240$ lines.

\medskip

\itemth{b} $\d(12)=\d^*(12)=2$, so $\dim V(12)=2$. 
We denote by $L_{12}$ the line $\PM(V(12)) \subset \PM(V)=\PM^3(\CM)$. 
Then $W(12)=C_W(w_{12})$ is a reflection group for its action 
on $V(12)$, and its degrees are $(12,24)$. So $|W(12)|=12\cdot 24=288$ 
by~\cite[(\ref{eq:degres})]{bonnafe-sarti-1}, and so the $W$-orbit of 
$L_{12}$ contains $160$ lines.

\medskip

\itemth{c} $\d(20)=\d^*(20)=1$, so $\dim V(20)=1$. 
We denote by $z_{20} \in \PM^3(\CM)$ the point defined by the line $V(20)$. Then 
$W(20)=C_W(w_{20})=\langle w_{20} \rangle$ is cyclic of order $20$.

\medskip

\itemth{d} $\d(24)=\d^*(24)=1$, so $\dim V(24)=1$. 
We denote by $z_{24} \in \PM^3(\CM)$ the point defined by the line $V(24)$. Then 
$W(24)=C_W(w_{24})=\langle w_{24} \rangle$ is cyclic of order $24$.
\end{itemize}
It follows from the above discussion 
and~\cite[Theorem~\ref{theo:springer}(f)]{bonnafe-sarti-1} that, 
if $e \in \{8,12,20,24\}$, 
then the eigenvalues of $w_e$ are $\z_e^{-7}$, $\z_e^{-11}$, 
$\z_e^{-19}$ and $\z_e^{-23}$ and so
\equat\label{eq:det-we-g31}
\det(w_e)=\z_e^{-60}=
\begin{cases}
-1 & \text{if $e \in \{8,24\}$,}\\
1 & \text{if $e \in \{12,20\}$.} \\
\end{cases}
\endequat

\bigskip

\begin{pro}\label{prop:lines}
With the above notation, we have that $z_{24} \in \ZC(f_{3,\l})$, 
that $L_8$ and $L_{12}$ are contained in $\ZC(f_{3,\l})$ 
and that $z_{20} \not\in \ZC(f_{3,\l})$.
\end{pro}

\bigskip

\begin{proof}
The facts that $z_{24} \in \ZC(f_{3,\l})$ and 
that $L_8$ and $L_{12}$ are contained in $\ZC(f_{3,\l})$ 
follow from~\cite[Lemma~\ref{lem:trivial}]{bonnafe-sarti-1}. 

Again by~\cite[Lemma~\ref{lem:trivial}]{bonnafe-sarti-1}, we have 
$f_1(z_{20})=f_2(z_{20})=f_4(z_{20})=0$, so that we cannot have 
$f_3(z_{20})=0$ since $0$ is the only common zeros of the fundamental 
invariants $(f_k)_{1 \le k \le 4}$. Hence $z_{20} \not\in \ZC(f_{3,\l})$.
\end{proof}

\bigskip

This shows in particular that $\ZC(f_{3,\l})$ contains 
at least $400$ lines (the $W$-orbit of $L_8$ of length $240$ 
and the $W$-orbit of $L_{12}$ of length $160$): it can be shown 
that, for $\l$ generic, these are the only lines 
contained in $\ZC(f_{3,\l})$.

\bigskip

Now, let 
$$B = \ZC(f_3) \cap \ZC(f_1f_2)$$
denote the {\it base locus} of the family 
$(\ZC(f_{3,\l}))_{\l \in \CM^\times}$. We write 
$$B_1=\ZC(f_3) \cap \ZC(f_1)\qquad\text{and}\qquad
B_2=\ZC(f_3) \cap \ZC(f_2),$$
so that 
$$B=B_1 \cup B_2.$$
By~\cite[Theorem~\ref{theo:springer}$($d$)$]{bonnafe-sarti-1}, we have
\equat\label{eq:z1-z2}
B_1=\bigcup_{x \in W} x(L_{12})\qquad\text{and}\qquad
B_2=\bigcup_{x \in W} x(L_8).
\endequat
We denote by $B'$ the image of $B$ in $\PM(V)/W'$ 
(it is the base locus of the family $(X_\l^{31})_{\l \in \CM}$). 
Then 
$$B'=B_1' \cup B_2',$$
where $B_j'$ denotes the image of $B_j$. 
It can be checked that the stabilizers $C_W(w_8)$ and $C_W(w_{12})$ 
of $L_8$ and $L_{12}$ in $W$ respectively are not contained in $W'$. 
So $B_1'$ (resp. $B_2'$) is also the image of $L_{12}$ 
(resp. $L_8$), hence it is a (possibly singular) rational curve.

\bigskip

\begin{pro}\label{prop:x1-x2}
The rational curve $B_2'$ is smooth, while the rational 
curve $B_1'$ has singularities $A_1 + A_2$.
\end{pro}

\bigskip

\begin{proof}
From the explicit formula for $r_\l$ given in~\S\ref{subsub:rl}, 
we have
$$B_1'=\{[z_1:z_2:z_4:t] \in \PM(1,1,1,3)~|~z_1=0~\text{and}~t^2=z_2z_4^3(z_2-z_4)^2\}$$
$$B_2'=\{[z_1:z_2:z_4:t] \in \PM(1,1,1,3)~|~z_2=t=0\}.
\leqno{\text{and}}$$
So $B_2' = \PM(1,1) = \PM^1(\CM)$ as expected. 

Let us now consider the case of $B_1'$. An easy computation in the affine 
charts defined by $z_2 \neq 0$ and $z_4 \neq 0$ gives 
two singular points $[0:1:0:0]$ and $[0:1:1:0]$ which are singularities 
of type $A_2$ and $A_1$ respectively.
\end{proof}

\bigskip

Note that the set theoretic intersection of $B_1'$ and $B_2'$ consists of only one point 
(let us call it $z_{24}'$ as it is the image of $z_{24} \in \ZC(f_{3,\l}) \subset \PM(V)$). 
Its coordinates are given by
$$z_{24}'=[0:0:1:0] \in X_\l^{31} \subset \PM(1,1,1,3).$$
Its image $\zba_{24}=[0:0:1] \in \PM^2(\CM)$ is a smooth point of the 
branch locus $R_\l$ (for all values of $\l$, because $r_\l(0,0,1)=1 \neq 0$). 

\bigskip

\begin{remark}\label{rem:intersection-g31}
Let $\Bba_1$ and $\Bba_2$ denote the respective images of $B_1'$ and $B_2'$ 
in $X_\l^{31}/\langle \s \rangle = \PM^2(\CM)$. Then $\Bba_1$ (resp. $\Bba_2$) 
is the line defined by the equation $z_1=0$ (resp. $z_2=0$). Note that 
the morphism $B_2' \longto \Bba_2$ is an isomorphism (as $B_2'$ is contained 
in the ramification locus) while the morphism $B_1' \longto \Bba_1$ is a morphism 
of degree $2$.

Recall that the branch locus of $X_\l^{31} \to \PM^2(\CM)$ is the union of 
$\Bba_2$ and $R_\l'=\ZC(r_\l)$. So 
$$\Bba_2 \cap R_\l' = \{[z_1:0:z_4] \in \PM^2(\CM)~|~
z_4^3(z_1-z_4)^2=0\}.$$
The {\it set} $\Bba_2 \cap R_\l'$ contains two points 
$d_6$ and $a_3$ of respective multiplicity $3$ and $2$
and whose coordinates are given by
$$d_6=[1:0:0]\qquad\text{and}\qquad a_3=[1:0:1].$$
They do not depend on $\l$. We will see in Corollary~\ref{coro:g31-singularites} and 
Proposition~\ref{prop:g31-icosic} that, if $\l \neq 0$,  
then $d_6$ is always a $D_6$ singularity of $R_\l$ while $a_3$ is an $A_3$ singularity 
except whenever $\l = \l_2$ (in which case it is a $D_5$ singularity).\finl
\end{remark}

\bigskip

\subsection{Singularities}\label{sub:sing-g31}
We wish to determine the list of singularities of $X_\l^{31}$. 
We gather in the next proposition some helpful general facts, 
from which we can deduce the list of singularities of $X_\l^{31}$ 
thanks to a few computations with {\sc Magma}.

\bigskip

\begin{pro}\label{prop:g31-prelim}
Let $v \in V \setminus \{0\}$ and let $z=[v]$. We assume that 
$z$ is a smooth point of $\ZC(f_{3,\l})$ and we denote 
by $z'$ its image in $X_\l^{31}$. 
\begin{itemize}
\itemth{a} If $|W_v|=1$ or $2$, then $z'$ is smooth.

\itemth{b} If $z \in B$ and 
$W_v$ has rank $2$, then $\Trm_z(\ZC(f_{3,\l}))$ 
together with its action of $W_z$ does not depend on $\l$.

\itemth{c} If $P$ is a parabolic subgroup of rank $2$ and 
if $z \in (\ZC(f_{3,\l}) \setminus B) \cap \PM(V^P)$, 
then $W_v=P$ and $W_z=P\langle w_4 \rangle$.
\end{itemize}
\end{pro}

\bigskip

\begin{proof}
(a) Assume first that $W_v=\{1\}$. Then $W_z=<w_{e_z}>$ and 
$4$ divides $e_z$ 
(see~\cite[\S\ref{sub:stab}, Fact~(a)]{bonnafe-sarti-1}). 
Since $e_z$ divides one of the degrees, we have 
$e_z \in \{4,8,12,20,24\}$. Note that $e_z \neq 20$ by 
Proposition~\ref{prop:lines}.

If $e_z=4$, then $(PW)_z=\{1\}$ and so $z'$ is smooth. 
If $e_z \in \{8,12,24\}$, then $\d(e_z)=\d^*(e_z)$ 
and the eigenvalues of $w_{e_z}$ on the tangent space 
$\Trm_z(\ZC(f_{3,\l}))$ are given 
by~\cite[Corollary~\ref{coro:vp-tangent}(b)]{bonnafe-sarti-1} 
and the determinant of $w_{e_z}$ is given by~(\ref{eq:det-we-g31}). 
So we get:
\begin{itemize}
\item If $e_z=8$, then $\det(w_{e_z})=-1$ and so $W_z'=\langle w_8^2 \rangle
= \langle w_4 \rangle$. So $(PW')_z=\{1\}$, which implies that $z'$ is smooth.

\item If $e_z=12$, then $\det(w_{e_z})=1$ and the eigenvalues 
of $w_{e_z}$ on $\Trm_z(\ZC(f_{3,\l}))$ are $\z_{12}^{-8}$ 
and $\z_{12}^{-24}=1$, so $w_{e_z}$ acts as a reflection 
on $\Trm_z(\ZC(f_{3,\l}))$. This implies that $z'$ is smooth.

\item If $e_z=24$, then $\det(w_{e_z})=-1$ so 
$W_z'=\langle w_{e_z}^2 \rangle$. Moreover, the eigenvalues 
of $w_{e_z}^2$ on $\Trm_z(\ZC(f_{3,\l}))$ are $\z_{24}^{-16}$ 
and $\z_{24}^{-24}=1$, so $w_{e_z}^2$ acts as a reflection 
on $\Trm_z(\ZC(f_{3,\l}))$. This implies that $z'$ is smooth.
\end{itemize}
This shows~$($a$)$ whenever $W_v=\{1\}$.

\smallskip

Let us now assume that $|W_v|=2$. Since $w_{e_z}$ normalizes 
$W_v$, this means that $w_{e_z}$ commutes with the non-trivial 
element of $W_v$, which is a reflection. But a {\sc Magma} 
computation shows that $w_e$ does not commute 
with any reflection if $e \in \{8,12,24\}$. So $e_z=4$, 
which means that $(PW')_z=\{1\}$. So $z'$ is smooth.

\medskip

(b) Assume that $z \in B$ and that 
$W_v$ has rank $2$. Then $\Trm_z(\ZC(f_{3,\l}))$ 
is a dimension $2$ subspace of $\Trm_z(\PM(V))$ 
which is stable under the action of $W_v$: but 
$\Trm_z(\PM(V))=V/z$ endowed with the natural 
action of $W_v$ which is of rank $2$, so 
there is a unique $W_v$-stable dimension $2$ subspace 
of $\Trm_z(\PM(V))$. This shows~(b).

\medskip

(c) Assume that $P$ is a parabolic subgroup of rank $2$ and 
that $z \in (\ZC(f_{3,\l}) \setminus B) \cap \PM(V^P)$. 
The fact that $z \not\in B$ implies that 
$e_z \not\in \{8,12,24\}$ by~(\ref{eq:z1-z2}). 
This shows that $W_z=W_v\langle w_4 \rangle$. 
On the other hand, $P = W_v$ by~(\ref{eq:stab}).
\end{proof}

\bigskip

\begin{cor}\label{coro:g31-singularites}
If $\l \in \CM$ is such that $\ZC(f_{3,\l})$ is smooth, then 
$X_\l^{31}$ has singularities $D_6+A_3+3\,A_2+2\,A_1$.
\end{cor}

\bigskip

\begin{proof}
The previous proposition shows that it is sufficient to 
determine a set of representatives of conjugacy classes 
of parabolic subgroups $P$ of rank $2$ and to determine 
the action of $W_z$ on $\Trm_z(\ZC(f_{3,\l}))$ for all 
$z \in \ZC(f_{3,\l}) \cap \PM(V^P)$. Let 
$$W_{14}=\langle s_1,s_4 \rangle,\quad
W_{15}=\langle s_1,s_5 \rangle\quad\text{and}\quad
W_{123}=\langle s_1,s_2,s_3 \rangle.$$
We set $N_k=\Nrm_W(W_k)$ 
and $L_k=\PM(V^{W_k})$ for $k \in \{14,15,123\}$. 
Computations with {\sc Magma} show that:
\begin{itemize}
\item $W_{14}$, $W_{15}$, $W_{123}$ are representatives 
of conjugacy classes of parabolic subgroups of rank $2$.

\item $W_{14}$ is a Coxeter group of type $A_2$ and 
$|N_{14}/W_{14}\langle w_4 \rangle|=6$. Moreover:
\begin{itemize}
\item $L_{14} \cap B_1$ contains $2$ elements 
which form a single $N_{14}$-orbit. If $z \in L_{14} \cap B_1$, 
then the action of $W_z'$ on $\Trm_z(\ZC(f_{3,\l}))$ 
can be computed for a single value of $\l$ thanks 
to Proposition~\ref{prop:g31-prelim}(b), and it can be checked 
that it acts as a reflection group, so the image of $z$ is smooth.

\item $L_{14} \cap B_2 = \vide$.

\item So it remains $18$ points in 
$(\ZC(f_{3,\l}) \setminus B) \cap L_{14}$: since the 
stabilizers of these points are equal to $W_{14}\langle w_4 \rangle$ 
by Proposition~\ref{prop:g31-prelim}(c), their 
$N_{14}$-orbits have cardinality $6$, 
so there are $3$ such orbits, each leading to an $A_2$-singularity 
because $W_{14}$ is of type $A_2$.
\end{itemize}

\item $W_{15}$ is a Coxeter group of type $A_1 \times A_1$ and 
$|N_{15}/W_{15}\langle w_4 \rangle|=8$. Moreover:
\begin{itemize}
\item $L_{15} \cap B_1 = \vide$.

\item $L_{15} \cap B_2$ contains $4$ elements 
which form a single $N_{15}$-orbit. If $z \in L_{15} \cap B_2$, 
then the action of $W_z'$ on $\Trm_z(\ZC(f_{3,\l}))$ 
can be computed for a single value of $\l$ thanks 
to Proposition~\ref{prop:g31-prelim}(b), and it can 
then be checked that the image of $z$ is an $A_3$-singularity.

\item So it remains $16$ points in 
$(\ZC(f_{3,\l}) \setminus B) \cap L_{15}$: since the 
stabilizers of these points are equal to $W_{15}\langle w_4 \rangle$ 
by Proposition~\ref{prop:g31-prelim}(c), their 
$N_{15}$-orbits have cardinality $8$, 
so there are $2$ such orbits, each leading to an $A_1$-singularity 
because $W_{15}$ is of type $A_1 \times A_1$.
\end{itemize}

\item $W_{123}$ is a complex reflection group of type $G(4,2,2)$ and 
$|N_{123}/W_{123}\langle w_4 \rangle|=24$. Moreover:
\begin{itemize}
\item $L_{123} \cap B_1$ contains $8$ elements 
which form a single $N_{123}$-orbit. Again, 
Proposition~\ref{prop:g31-prelim}(b) allows an easy computation 
which implies that the image of $z$ is smooth.

\item $L_{123} \cap B_2$ contains $12$ elements 
which form a single $N_{123}$-orbit. Again, 
Proposition~\ref{prop:g31-prelim}(b) allows an easy computation 
which implies that the image of $z$ is a $D_6$-singularity.

\item It remains no point in 
$(\ZC(f_{3,\l}) \setminus B) \cap L_{123}$.
\end{itemize}
\end{itemize}
The proof of the corollary is complete.
\end{proof}

\bigskip

%

\begin{pro}\label{prop:g31-icosic}
If $k \in \{145,245,1,2,3\}$, then the singularities of 
$X_{\l_k}^{31}$ are given by Table~\ref{table:k3}, i.e. 
\begin{itemize}
\itemth{a} The surface $X_{\l_{145}}^{31}=X_{-8/25}^{31}$ has singularities 
$D_6+D_5+A_3+2\,A_2$. 

\smallskip

\itemth{b} The surface $X_{\l_{245}}^{31}=X_{-81/175}^{31}$ has singularities 
$E_6+D_6+A_3+A_2+A_1$. 

\smallskip

\itemth{c} The surface $X_{\l_1}^{31}=X_{1}^{31}$ has singularities 
$D_6+A_5+A_3+A_2+2\,A_1$.

\smallskip

\itemth{d} The surface $X_{\l_2}^{31}=X_{-1/3}^{31}$ has singularities 
$D_6+D_5+3\,A_2+A_1$.

\smallskip

\itemth{e} The surface $X_{\l_3}^{31}=X_{-1/2}^{31}$ has singularities 
$D_6+2\,A_3+2\,A_2+2\,A_1$.
\end{itemize}
\end{pro}

\bigskip

\begin{proof}
Using the formula for $r_\l$ given in the previous subsection, 
one can can easily obtain the equation of the branch locus 
$R_{\l_k}$ for the five values of $k$. The singularities 
of the curve $R_{\l_k}$ are then easily determined thanks to {\sc Magma} 
and we conclude thanks 
to~\cite[Proposition~\ref{prop:double cover}]{bonnafe-sarti-1}.
\end{proof}

\bigskip

This has the following consequence, which confirms some of the results of 
Table~\ref{table:k3}:

\bigskip

\begin{cor}\label{coro:picard-g31}
Let $\l \in \CM^\times$. Then:
\begin{itemize}
\itemth{a} If $\l \not\in \{\l_{145},\l_{245},\l_1,\l_2,\l_3\}$, then $\rhob(\Xti_\l^{31}) \ge 18$. 

\itemth{b} For generic $\l$, we have $\rhob(\Xti_\l^{31}) = 18$.

\itemth{c} If $\l \in \{\l_{145},\l_{245},\l_1,\l_2,\l_3\}$, then $\rhob(\Xti_\l^{31}) = 19$.
\end{itemize}
\end{cor}

\bigskip

\begin{proof}
Let $\l \in \CM^\times$. We denote by $m$ the sum of the Milnor numbers of the 
singularities of $X_\l^{31}$ (i.e., $m$ is the number of smooth rational 
curves in the exceptional divisors of the resolution 
$\pi_\l : \Xti_\l^{31} \longto X_\l^{31}$). Then $\rhob(\Xti_\l^{31}) \ge 1+m$ 
since $X_\l^{31}$ is projective, 
so one can check from Corollary~\ref{coro:g31-singularites} and 
Proposition~\ref{prop:g31-icosic} the following two facts:
$$\begin{cases}
\text{If $\l \not\in \{\l_{145},\l_{245},\l_1,\l_2,\l_3\}$, 
then $\rhob(\Xti_\l^{31}) \ge 18$;}\\
\text{If $\l \in \{\l_{145},\l_{245},\l_1,\l_2,\l_3\}$, then $\rhob(\Xti_\l^{31}) \ge 19$.}
\end{cases}
\leqno{(\clubsuit)}$$
Note that $(\clubsuit)$ proves the inequality stated in~(a).

\medskip

Let us now prove the equalities stated in~(b) and~(c). 
We shall use the methods developed 
by van Luijk~\cite{luijk} and Elsenhans and Jahnel~\cite[\S{3.3.1}]{elsenhans-jahnel}, 
based on the Artin-Tate Conjecture (proved by Nygaard and Ogus for K3 surfaces~\cite{nygaard-ogus} 
in characteristic~$\ge 5$), but we adapt them to the singular case. 
For this, assume that $\l \in \QM$ and let $\PC_\l$ denote the set of 
prime numbers $p$ such that:
\begin{itemize}
\itemth{1} $p \ge 5$ and $p$ does not divide any denominator of any coefficient of $r_\lambda$. 
(so that we can define a reduction of $X_\l^{31}$ modulo $p$, which will be defined 
over $\FM_{\! p}$ and will be denoted by $(X_\l^{31})_p$: we also denote by $(R_\l)_p$ 
the reduction modulo $p$ of the ramification locus of $\pi_\l$).

\itemth{2} If $\OC_\l$ is the ring of integers of the minimal number field $K_\l$ containing 
the coordinates of all the singular points of $X_\l^{31}$ and if $\pG_\l$ is a prime ideal 
of $\OC_\l$ lying over $p$, then $\OC_\l/\pG_\l=\FM_{\! p}$ and all the singular points of 
$(X_\l^{31})_p$ have coordinates in $\FM_{\! p}$ 
and are the reduction modulo $\pG_\l$ of the singular points of $X_\l^{31}$. 

\itemth{3} If $x \in X_\l^{31}$ is a singular point, then its reduction modulo $p$ 
is an ADE singularity of $(X_\l^{31})_p$ of the same type as $x$.
\end{itemize}
So let $p \in \PC_\l$. We denote by $(\Xti_\l^{31})_p$ the minimal resolution 
of the K3 surface $(X_\l^{31})_p$. Then $(\Xti_\l^{31})_p$ is the reduction modulo $p$ 
of $X_\l^{31}$ by~(1),~(2) and~(3), because $(\Xti_\l^{31})_p$ is obtained 
from $X_\l^{31}$ by the same sequence of blow-ups. This shows in particular 
that $\Xti_\l^{31}$ has good reduction modulo $p$ (i.e. remains smooth) 
and that its reduction modulo $p$ 
is exactly $(\Xti_\l^{31})_p$. 

We denote by $P_{\l,p} \in \ZM[T]$ (resp. $\Pti_{\l,p} \in \ZM[T]$) the Weil polynomial of $(X_\l^{31})_p$ 
(resp. $(\Xti_\l^{31})_p$), namely the characteristic polynomial of the Frobenius 
map on the second $\ell$-adic cohomology group of $(X_\l^{31})_p$ (resp. $(\Xti_\l^{31})_p$). 
Note that the polynomial $P_{\l,p}$ can be computed explicitly (and efficiently!) 
thanks to the command {\tt WeilPolynomialOfDegree2K3Surface} in 
{\sc Magma} and that
$$\Pti_{\l,p}=(T-p)^m P_{\l,p},\leqno{(\diamondsuit)}$$
where we recall that $m$ is the number of irreducible components of the 
exceptional divisors of the minimal resolution of $X_\l^{31}$ (or of $(X_\l^{31})_p$, 
as they are all defined over $\FM_{\! p}$ by~(2) and~(3)). Let $\r_{\l,p}$ denote the $(T-p)$-valuation 
of $P_{\l,p}$ and let $Q_{\l,p}=P_{\l,p}/(T-p)^{\r_{\l,p}}$. Let $\r_{\l,p}^g$ 
denote the number of root of $Q_{\l,p}$ of the form $\z p$, where $\z$ is a root of unity 
(note that $\r_{\l,p}^g \ge \r_{\l,p}$). 
Also, we denote by $D_\l \in \QM^\times$ the discriminant of 
the Picard group of $X_\l^{31}$. 

We denote by $\Pic_g((\Xti_\l^{31})_p)$ the geometric Picard group of $(\Xti_\l^{31})_p$, 
namely the Picard group of $\overline{\FM}_{\! p} \times_{\FM_{\! p}} (\Xti_\l^{31})_p$. 
Then Artin-Tate Conjecture and $(\diamondsuit)$ say that
$$m+\r_{\l,p} = {\mathrm{rk}} \Pic((\Xti_\l^{31})_p)
\qquad\text{and}\qquad
m+\r_{\l,p}^g = {\mathrm{rk}} \Pic_g((\Xti_\l^{31})_p).\leqno{(\heartsuit)}$$
Reduction modulo $p$ induces an injective map $\Pic \Xti_\l^{31} \injto \Pic_g((\Xti_\l^{31})_p)$ 
(see~\cite[Proposition~6.2]{luijk-quartic}). 
Hence
$$\rhob(\Xti_\l^{31}) \le m+\r_{\l,p}^g.\leqno{(\spadesuit)}$$
Moreover, if these two groups have the same rank, then their discriminant are equal modulo 
$\QM^{\times 2}$. By Artin-Tate Conjecture and~$(\heartsuit)$, this forces
$$\text{If $\rhob(\Xti_\l^{31}) \!=\! m+\r_{\l,p}^g\!=\! m+\r_{\l,p}$, 
then $D_\l \equiv p^{m+\r_{\l,p}^g-21} Q_{\l,p}(p) \!\!\mod \QM^{\times 2}$.}\leqno{(\spadesuit^+)}$$
With all these tools in hand, we proceed as follows (numerical results 
stated below are obtained with {\sc Magma}).

\medskip

(b) By~$(\clubsuit)$, $\rhob(\Xti_{-1/4}^{31}) \ge 18$. Note that $m=17$ in this case. 
On the other hand, $193 \in \PC_{-1/4}$ and 
$$P_{-1/4,193}=(T - 193)(T^4 + 212\,T^3 + 10422\,T^2 + 7896788\,T + 1387488001).$$
This shows that $Q_{-1/4,193}=T^4 + 212\,T^3 + 10422\,T^2 + 7896788\,T + 1387488001$. 
Since this polynomial has no root of the form $193\z$ with $\z$ a root of unity, 
we get that $\r_{-1/4,193}=\r_{-1/4,193}^g=1$ and so $\rhob(\Xti_{-1/4}^{31}) \le 18$ by $(\spadesuit)$. 
This proves~(b) for $\l =-1/4$ and so this proves~(b) for $\l$ generic.

\medskip

(c) We explain how to prove~(c) whenever $\l=\l_{145}=-8/25$, the other cases 
being treated similarly. Note first that $m=18$ in this case. 
By~$(\clubsuit)$, $\rhob(\Xti_{-8/25}) \in \{19,20\}$. 
Note that $23$ and $47$ belong to $\PC_{-8/25}$. 
We have
$$P_{-8/25,23}=(T - 23)^2(T^2 + 38\,T + 529)\quad\text{and}\quad 
Q_{-8/25,23}(23)/23 \equiv 21 \!\!\mod \QM^{\times2},$$
$$P_{-8/25,47}=(T - 47)^2(T^2 + 22\,T + 2209)\quad\text{and}\quad 
Q_{-8/25,47}(47)/47 \equiv 29 \!\!\mod \QM^{\times 2}.$$
Assume that $\rhob(\Xti_{-8/25})=20$. Then 
$$20=\rhob(\Xti_{-8/25})=m+\r_{-8/25,23}=m+\r_{-8/25,23}^g=m+\r_{-8/25,47}=m+\r_{-8/25,47}^g,$$
so it follows from~$(\spadesuit^+)$ that $21 \equiv 29 \!\!\mod \QM^{\times 2}$, 
which is impossible. So $\rhob(\Xti_{-8/25}^{31})=19$, as expected.
\end{proof}

\bigskip

\begin{remark}[Supersingular surfaces]
Keep the notation of the proof of Corollary~\ref{coro:picard-g31}. 
For each exceptional value of $\l$ 
(i.e. $\l \in \{\l_{145},\l_{245},\l_1,\l_2,\l_3\}$) there exist prime 
numbers $p$ such that $\Xti_\l^{31}$ has good reduction modulo $p$ 
and $(\Xti_\l^{31})_p$ is a supersingular variety (i.e. has geometric 
Picard number $22$). We give here a (non-exhaustive) list of examples. 
So assume that $(\l,p)$ is a pair where $\l \in \{\l_{145},\l_{245},\l_1,\l_2,\l_3\}$ 
and $p$ is a prime number such that:
\begin{itemize}
\item[$\bullet$] If $\l = \l_{145}=-8/25$, then $p \in \{59,73,89\}$.

\item[$\bullet$] If $\l = \l_{245}=-81/175$, then $p \in \{31,47,73\}$.

\item[$\bullet$] If $\l = \l_1=1$, then $p = 43 $.

\item[$\bullet$] If $\l = \l_2=-1/3$, then $p = 337$.

\item[$\bullet$] If $\l = \l_3=-1/2$, then $p \in \{73,79 \}$.
\end{itemize}
Then $(\Xti_\l^{31})_p$ is supersingular.\finl
\end{remark}

\bigskip

\begin{remark}\label{rem:g31-particular}
Note that, generically, $r_\l$ is irreducible. However, 
$r_{\l_1}$ and $r_{\l_3}$ are not irreducible\footnote{We do not know if there are other values 
of $\l$ such that $r_\l$ is not irreducible.}:
\begin{itemize}
\item[$\bullet$] The quintic $R_{\l_1}'=R_1'$ is the union of a smooth irreducible conic 
and an irreducible cubic. More detail about this case will be given in~\S\ref{sub:g31-cubic}.

\item[$\bullet$] The quintic $R_{\l_3}'=R_{-1/2}'$ is the union of a line and 
an irreducible quartic. More detail about this case will be given in~\S\ref{sub:g31-quartic}.\finl
\end{itemize}
\end{remark}

\bigskip

\bigskip

\subsection{Complements}\label{sub:boissiere-sarti}
The experienced reader might have noticed that 
$$f_1=\Sigma(x^8)+14\Sigma(x^4y^4)+168 x^2y^2z^2t^2$$
is the polynomial which defines the smooth octic containing 352 lines 
constructed by Boissi\`ere and the second author~\cite{boissieresarti}. 
We will revisit here this example. 

Let 
$$\s=\frac{\sqrt{2}}{2}
\begin{pmatrix}
-1 & 0 & 0 & -1 \\
0 & 1 & 1 & 0 \\
0 & 1 & -1 & 0 \\
-1 & 0 & 0 & 1 \\
\end{pmatrix}.$$
Then $\s(s_i)=s_{5-i}$ if $i \in \{1,2,3,4\}$ and $\z_8 \s \in G_{31}$. 
Moreover, $\z_8 \s$ normalizes 
the subgroup $G_{28}$. In~\cite{boissieresarti}, 
the polynomial $f_1$ was constructed as a particular invariant of 
the one-parameter family of fundamental invariants of degree $8$ of the group 
$\langle \s \rangle \ltimes G_{28}^\slv$ (which is contained in $\langle \z_8 \rangle G_{31}$), 
but it turns out that this 
is exactly the one which is invariant by $G_{31}$. 

The 352 lines on $\ZC(f_1)$ are divided into two $G_{31}$-orbits: 
one of size 160 and one of size 192. We explain here how to construct 
these two orbits. 

First, as $12$ does not divide $8$, the $G_{31}$-orbit of $L_{12}$ 
is contained in $\ZC(f_1)$ 
by~\cite[Lemma~\ref{lem:trivial}]{bonnafe-sarti-1}, 
so this explains the first orbit 
with $160$ lines. 
For constructing the second orbit, one requires some 
more material. Let $\WC=G_{37}=\Wrm(E_8)$ acting 
on a vector space $V_8$ of dimension $8$. 
The list of degrees (resp. codegrees) of $\WC$ is 
$(2,8,12,14,18,20,24,30)$ (resp. 
$(0,6,10,12,16,18,22,28)$). Applying Theorem~\ref{theo:springer} 
with $\WC$ and $e=4$ shows that there exists an element $w_4 \in \WC$ 
such that $\dim V_8(w_4,i) = 4$, $\WC(V_8(w_4,i))=\{1\}$ and 
$\WC_{V_8(w_4,i)}$ acts on $V_8(w_4,i)$ as a reflection group 
whose list of degrees is $(8,12,20,24)$: in fact, 
$\WC_{V_8(w_4,i)} \simeq G_{31}$ 
(as a reflection group). 
Therefore, we may identify $V$ with $V_8(w_4,i)$ and 
$G_{31}$ with $\WC_{V_8(w_4,i)}$. 

Now, let $\phi=(1+\sqrt{5})/2$ be the golden ratio. 
By~\cite[\S{3}]{lusztig}, there exists an automorphism $\ph$ of $V_8$ 
satisfying $\ph^2=\ph+\Id_{V_8}$ and such that $\dim V_8(\ph,\phi) = 4$ 
and $\WC_{V_8(\ph,\phi)}$ acts faithfully on $V_8(\ph,\phi)$ 
as the complex reflection group $\Wrm(H_4)=G_{30}$. 
Using again Theorem~\ref{theo:springer} with $\WC_{V_8(\ph,\phi)} \simeq G_{30}$, 
we see that we may choose the above element $v$ as belonging 
to $\WC_{V_8(\ph,\phi)}$. Moreover, $E=V_8(v,i) \cap V_8(\ph,\phi)$ 
has dimension $2$, and its stabilizer $\WC_E$ acts faithfully on $E$ 
as the complex reflection group $G_{22}$, whose list of 
degrees is $(12,20)$ (as they are the only degrees of $W$ which 
are divisible by $4$). Hence, the restriction of $f_1$ to $E$ 
having degree $8$ and being invariant under $\WC_E=C_W(w_4)$, 
this implies that $f_1$ vanishes on $E$. So, if we let $L'=\PM(E)$, then 
$L'$ is a line contained in $\ZC(f_1)$, whose stabilizer 
in $G_{31}$ has order $12 \times 20=240$. This shows that 
the $G_{31}$-orbit of $L'$ contains $192$ lines which are all contained 
in $\ZC(f_1)$.

\bigskip

\section{The group $G_{31}$ (continued): elliptic fibrations}

\medskip

We will use here the constructions of Appendix~\ref{appendix-3}. 
Let $x$ be a singular point of the branch locus $R_\l$. Since 
$x$ belongs to the branch locus, there is a unique point $\xdo \in X_\l^{31}$ 
above $x$. 
Let 
$p_x : \PM^2(\CM)\setminus\{x\} \to \PM^1(\CM)$ be the 
projection from the point $x$. We denote by 
$\hat{\PM}_x^2(\CM)$ the blow-up of $\PM^2(\CM)$ at $x$ 
and by $\Xhat_\l^x$ the blow-up of $X_\l^{31}$ at $\xdo$. Then:
\begin{itemize}
\item[$\bullet$] The projection $p_x$ lifts and extends to a morphism 
$\phat_x : \hat{\PM}_x^2(\CM) \to \PM^1(\CM)$.

\item[$\bullet$] Since $X_\l^{31}$ has only ADE singularities, the map $\xi_\l : X_\l^{31} \to \PM^2(\CM)$ 
lifts to a map $\xih_\l^x : \Xhat_\l^x \longto \hat{\PM}_x^2(\CM)$ 
(see Proposition~\ref{prop:blowing-double-cover}). 

\item[$\bullet$] Since $X_\l^{31}$ has only ADE singularities, its minimal 
resolution is obtained by successive blow-ups of singular points. 
In particular, the morphism $\pi_\l : \Xti_\l^{31} \to X_\l^{31}$ 
factorizes through $\pih_\l^x : \Xti_\l^{31} \to \Xhat_\l^x$.
\end{itemize}
Altogether, this gives a well-defined morphism of varieties
$$\pht_\l^x = \phat_x \circ \xih_\l^x \circ \pih_\l^x : 
\Xti_\l^{31} \longto \PM^1(\CM),$$
i.e. an elliptic fibration.

This gives lots of elliptic fibrations, and the particular values $\l_k$ of $\l$ 
must also be treated separately. For this reason, we will not compute the 
singular fibers in all cases. We will just provide general facts about sections, 
use them to determine the intersection graph of the curves 
contained in $\pi_\l^{-1}(B_2')$ and just focus on singular fibers of the fibration $\pht_\l^{d_6}$. 

\bigskip

\begin{quotation}
\noindent{\bf Question.} {\it Are there other elliptic fibrations on the surface $\Xti_\l^{31}$?}
\end{quotation}


\medskip

\subsection{Sections} 
Let us first discuss the question of sections of the elliptic fibration 
associated with $\pht_\l^x$, using Proposition~\ref{prop:section-double-cover}. 
For this, let $\Ehat_x$ denote the exceptional divisor of the blow-up 
$\hat{\PM}_x^2(\CM)$ (it is isomorphic to $\PM^1(\CM)$ and maps isomorphically to 
$\PM^1(\CM)$ through $\phat_x$). 
If we denote by $m$ the Milnor number 
of $\xdo$ and by $\D_1^x$,\dots, $\D_m^x$ the smooth rational 
curves of the exceptional divisor of $\Xhat_\l^\xdo$, then 
Proposition~\ref{prop:section-double-cover} implies that:
\begin{itemize}
\item If $m \ge 2$ and $x$ is an $A_m$-singularity (and if we assume that 
the smooth rational curves $\D_j^x$ are numbered so that the extremal 
vertices of their intersection graph are $\D_1^x$ and $\D_m^x$, then 
$\D_1^x$ and $\D_m^x$ are exchanged by $\s$ and are mapped isomorphically 
to $\Ehat_x$. This gives two sections $\th_x^\pm : \PM^1(\CM) \longto \Xti_\l^{31}$ 
of the elliptic fibration $\pht_\l^x$ satisfying $\th_x^-=\s \circ \th_x^+$.

\item If $x$ is not of type $A$, then only one of the smooth rational curves 
$\D_j^x$ maps isomorphically to $\Ehat_x$. This leads to a section 
$\th_x :\PM^1(\CM) \longto \Xti_\l^{31}$ of $\pht_\l^x$. 
\end{itemize}
Note also that any line $L$ of $\PM^2(\CM)$ not containing $x$ maps isomorphically through 
the projection $p_x$, so its inverse image $\Lhat\simeq L$ in $\hat{\PM}_x^2(\CM)$ 
maps isomorphically to $\PM^1(\CM)$ through $\phat_x$. Applied to the line $\Bba_2$, 
and using the fact that $\Bba_2$ lies in the branch locus (and so 
the map $B_2' \longto \Bba_2$ is an isomorphism), we see that, if $x \not\in \Bba_2$, 
then the elliptic fibration $\pht_\l^x$ admits a section $\th_x^B : \PM^1(\CM) \longto \Xti_\l^{31}$ 
whose image is the strict transform $\Bti_2'$ of $B_2'$ in $\Xti_\l^{31}$. 
We summarize the above discussion in the next proposition:

\bigskip

\begin{pro}\label{prop:section-g31}
Let $x$ be a singular point of $R_\l$. Then:
\begin{itemize}
\itemth{a} If $m \ge 2$ and if $x$ is an $A_m$ singularity, then the elliptic fibration 
$\pht_\l^x$ admits two sections $\th_x^\pm$ whose images are the two extremal smooth rational curves 
of the exceptional divisor $\pit_\l^{-1}(\xdo)$.

\itemth{b} If $x$ is not a type $A$ singularity, then the elliptic fibration 
$\pht_\l^x$ admits a section whose image is one of the smooth rational curves 
of the exceptional divisor $\pit_\l^{-1}(\xdo)$.

\itemth{c} If $x \not\in \Bba_2$ (i.e. if $x \not\in \{a_3,d_6\}$), then the elliptic fibration 
$\pht_\l^x$ admits a section whose image is $\Bti_2'$. 
\end{itemize}
\end{pro}

\bigskip

\subsection{Intersection graph in ${\boldsymbol{\pi_\l^{-1}(B_2')}}$ and the elliptic fibration ${\tilde{\varphi}}_\lambda^{a_3}$}
Recall that $a_3$ and $d_6$ are the only singular points of $R_\l$ belonging 
to $B_2'$. It will be interesting for computing Picard numbers and transcendental 
lattices to determine the intersection graph between the smooth rational 
curves of $\pi_\l^{-1}(\ddo_6)$, the ones of $\pi_\l^{-1}(\ado_3)$ and the strict transform $\Bti_2'$ 
of $B_2'$ in $\Xti_\l^{31}$. This will be done thanks to the elliptic fibrations 
constructed in this section. We need some notation. 
The point $\ddo_6 \in X_\l^{31}$ is always a $D_6$ singularity. We assume that 
the $6$ smooth rational curves $(\D_k^{d_6})_{1 \le k \le 6}$ of the exceptional 
divisor $\pi_\l^{-1}(\ddo_6)$ are numbered in such a way that the intersection graph 
is given by

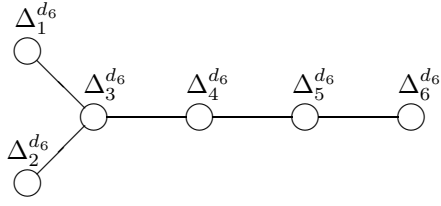
\centerline{\begin{picture}(300,78)
\put(75, 55){\circle{10}}\put(70,64){$\D_1^{d_6}$}
\put(75,  5){\circle{10}}\put(67,14){$\D_2^{d_6}$}
\put(78.5355,51.4645){\line(1,-1){18.1}}
\put(78.5355,8.5355){\line(1,1){18.1}}
\put(100, 30){\circle{10}}\put(96,39){$\D_3^{d_6}$}
\put(105, 30){\line(1,0){30}}
\put(140, 30){\circle{10}}\put(145, 30){\line(1,0){30}}\put(135,39){$\D_4^{d_6}$}
\put(180, 30){\circle{10}}\put(185, 30){\line(1,0){30}}\put(175,39){$\D_5^{d_6}$}
\put(220, 30){\circle{10}}\put(215 ,39){$\D_6^{d_6}$}
\end{picture}}
\noindent We denote by $m_3(\l)$ the Milnor number of the singularity $\ado_3$. 
If $\l \not=\l_2=-1/3$ (resp. $\l = \l_2$), then $a_3$ is an $A_3$ (resp. a $D_5$) 
singularity, so $m_3(\l)=3$ (resp. $m_3(\l)=5$) and we assume that the $m_3(\l)$ smooth rational curves 
$(\D_k^{a_3})_{1 \le k \le m_3(\l)}$ 
of the exceptional divisor $\pi_\l^{-1}(\ado_3)$ are numbered in such a way that the intersection graph 
is given by

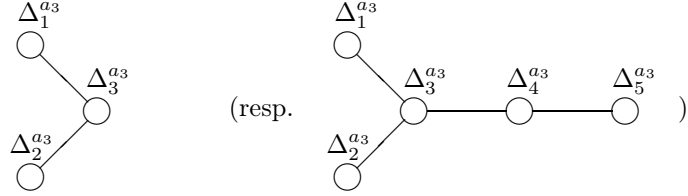
\centerline{\begin{picture}(300,75)
\put(25, 55){\circle{10}}\put(20,64){$\D_1^{a_3}$}
\put(25,  5){\circle{10}}\put(17,14){$\D_2^{a_3}$}
\put(28.5355,51.4645){\line(1,-1){18.1}}
\put(28.5355,8.5355){\line(1,1){18.1}}
\put(50, 30){\circle{10}}\put(46,39){$\D_3^{a_3}$}
\put(100,28){(resp.}
\put(145, 55){\circle{10}}\put(140,64){$\D_1^{a_3}$}
\put(145,  5){\circle{10}}\put(137,14){$\D_2^{a_3}$}
\put(148.5355,51.4645){\line(1,-1){18.1}}
\put(148.5355,8.5355){\line(1,1){18.1}}
\put(170, 30){\circle{10}}\put(166,39){$\D_3^{a_3}$}
\put(175, 30){\line(1,0){30}}
\put(210, 30){\circle{10}}\put(215, 30){\line(1,0){30}}\put(205,39){$\D_4^{a_3}$}
\put(250, 30){\circle{10}}\put(245,39){$\D_5^{a_3}$}
\put(270,28){)}
\end{picture}}
\noindent Now, if $x$ is a singular point of $X_\l^{31}$ different from $a_3$ and $d_6$ 
(there always exists such a point), then $(\pht_\l^x)^{-1}(p_x(d_6))$ and $(\pht_\l^x)^{-1}(p_x(a_3))$ 
are two singular fibers (because they contain $\pi_\l^{-1}(\ddo_6)$ and $\pi_\l^{-1}(\ado_3)$). 
Since $\Bti_2'$ is a section of the elliptic 
fibration $\pht_\l^x$ by Proposition~\ref{prop:section-g31}, $\Bti_2'$ meets 
$\pi_\l^{-1}(\ddo_6)$ and $\pi_\l^{-1}(\ado_3)$ transversally at only one curve with multiplicity $1$. 
Recall that the multiplicity $1$ curves of $\pi_\l^{-1}(\ddo_6)$ (resp. $\pi_\l^{-1}(\ado_3)$) 
are $\D_1^{d_6}$, $\D_2^{d_6}$ and $\D_6^{d_6}$
(resp. $\D_1^{a_3}$, $\D_2^{a_3}$ and $\D_{m_3(\l)}^{a_3}$, where $m_3(\l)$ denote the Milnor number of the singularity $\ado_3$). 

Since $\s(\Bti_2')=\Bti_2'$ and $\s(\D_1^{a_3})=\D_2^{a_3}$, this forces 
that $\Bti_2'$ meets $\pi_\l^{-1}(\ado_3)$ transversally at $\D_{m_3(\l)}^{a_3}$. 

To determine which curve of $\pi_\l^{-1}(\ddo_6)$ meets $\Bti_2'$, we use the 
elliptic fibration $\pht_\l^{a_3}$. First, $(\pht_\l^{a_3})^{-1}(p_{a_3}(d_6))$ contains 
$\pi_\l^{-1}(\ddo_6)$ and $\Bti_2'$. Moreover, $\pht_\l^{a_3}(\D_{m_3(\l)}^{a_3})$ is a 
point by Proposition~\ref{prop:section-double-cover}, so it must be the same point 
as $\pht_\l^{a_3}(\Bti_2')$, which is $\pht_\l^{a_3}(\ddo_6)=p_{a_3}(d_6)$. Therefore, 
$$(\pht_\l^{a_3})^{-1}(p_{a_3}(d_6))=\D_{m_3(\l)}^{a_3} \cup \Bti_2' \cup \Bigl(\bigcup_{k=1}^6 \D_k^{d_6}\Bigr).$$
The Kodaira-N\'eron classification of singular fibers then shows that the 
only possibility is that $(\pht_\l^{a_3})^{-1}(p_{a_3}(d_6))$ is of type $\Eti_7$ 
and $\Bti_2'$ meets $\D_1^{d_6}$ (by exchanging $\D_1^{d_6}$ and $\D_2^{d_6}$ if necessary). 
So we have shown most of the following lemma:

\bigskip

\begin{lemma}\label{lem:g31-graph}
The intersection graph of curves contained in $\pi_\l^{-1}(B_2')$ is as follows. If $\l \neq \l_2$  it is given by

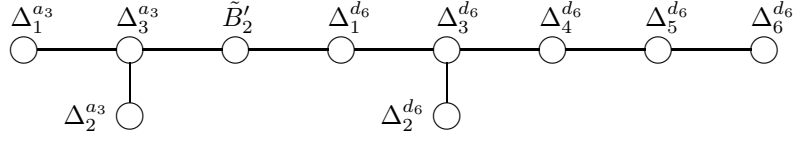
\centerline{\begin{picture}(300,57)
\put(10, 35){\circle{10}}\put(15,35){\line(1,0){30}}\put(5,44){$\D_1^{a_3}$}
\put(50, 35){\circle{10}}\put(55,35){\line(1,0){30}}\put(45,44){$\D_3^{a_3}$}
\put(50, 10){\circle{10}}\put(50,15){\line(0,1){15}}\put(25,7){$\D_2^{a_3}$}
\put(90, 35){\circle{10}}\put(95,35){\line(1,0){30}}\put(85,44){$\Bti_2'$}
\put(130, 35){\circle{10}}\put(125,44){$\D_1^{d_6}$}
\put(170, 10){\circle{10}}\put(145,7){$\D_2^{d_6}$}
\put(135,35){\line(1,0){30}}
\put(170,30){\line(0,-1){15}}
\put(170, 35){\circle{10}}\put(166,44){$\D_3^{d_6}$}
\put(175, 35){\line(1,0){30}}
\put(210, 35){\circle{10}}\put(215, 35){\line(1,0){30}}\put(205,44){$\D_4^{d_6}$}
\put(250, 35){\circle{10}}\put(255, 35){\line(1,0){30}}\put(245,44){$\D_5^{d_6}$}
\put(290, 35){\circle{10}}\put(285 ,44){$\D_6^{d_6}$}
\end{picture}}

\noindent  and if $\l = \l_2$ it is given by 

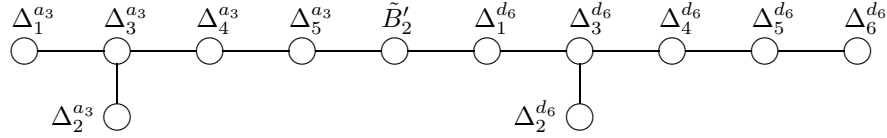
\centerline{\begin{picture}(335,57)
\put(10, 35){\circle{10}}\put(15,35){\line(1,0){25}}\put(5,44){$\D_1^{a_3}$}
\put(45, 35){\circle{10}}\put(50,35){\line(1,0){25}}\put(40,44){$\D_3^{a_3}$}
\put(45, 10){\circle{10}}\put(45,15){\line(0,1){15}}\put(20,7){$\D_2^{a_3}$}
\put(80, 35){\circle{10}}\put(85,35){\line(1,0){25}}\put(75,44){$\D_4^{a_3}$}
\put(115, 35){\circle{10}}\put(120,35){\line(1,0){25}}\put(110,44){$\D_5^{a_3}$}
\put(150, 35){\circle{10}}\put(155,35){\line(1,0){25}}\put(145,44){$\Bti_2'$}
\put(185, 35){\circle{10}}\put(180,44){$\D_1^{d_6}$}
\put(220, 10){\circle{10}}\put(195,7){$\D_2^{d_6}$}
\put(190,35){\line(1,0){25}}
\put(220,30){\line(0,-1){15}}
\put(220, 35){\circle{10}}\put(216,44){$\D_3^{d_6}$}
\put(225, 35){\line(1,0){25}}
\put(255, 35){\circle{10}}\put(260, 35){\line(1,0){25}}\put(250,44){$\D_4^{d_6}$}
\put(290, 35){\circle{10}}\put(295, 35){\line(1,0){25}}\put(285,44){$\D_5^{d_6}$}
\put(325, 35){\circle{10}}\put(320 ,44){$\D_6^{d_6}$}
\end{picture}}
\noindent Moreover:
\begin{itemize}
\itemth{a} The singular fiber $(\pht_\l^{a_3})^{-1}(p_{a_3}(d_6))$ is of type $\Eti_7$.

\itemth{b} The singular fiber $(\pht_\l^{d_6})^{-1}(p_{d_6}(a_3))$ is of type $\Dti_7$  
if $\l \neq \l_2$ and of type $\Dti_9$ if $\l=\l_2=-1/3$.
\end{itemize}
\end{lemma}

\bigskip

\begin{proof}
Only the statement~(b) has not been proved. First, $\pi_\l^{-1}(\ado_3)$ and $\Bti_2'$ 
are contained in $(\pht_\l^{d_6})^{-1}(p_{d_6}(a_3))$. Moreover, 
it follows from Proposition~\ref{prop:section-double-cover} that the curves 
$(\D_k^{d_6})_{1 \le k \le 4}$ are sent, through $\pht_\l^{d_6}$, to a single point 
of $\PM^1(\CM)$. Since $\Bti_2'$ meets $\D_1^{d_6}$, this point is necessarily 
$p_{d_6}(a_3)$. So 
$$(\pht_\l^{d_6})^{-1}(p_{d_6}(a_3))=\pi_\l^{-1}(\ado_3) \cup \Bti_2' \cup 
\Bigl(\bigcup_{k=1}^4 \D_k^{d_6}\Bigr),$$
and the result follows from the description of the intersection graph.
\end{proof}

\bigskip

\subsection{The elliptic fibration ${\boldsymbol{\pht_\l^{d_6}}}$} 
Since $d_6=[1:0:0]$, the maps $p_{d_6} : \PM^2(\CM) \setminus \{d_6\} \longto \PM^1(\CM)$ 
and $\ph_{d_6} : X_\l^{31} \setminus \{\ddo_6\} \longto \PM^1(\CM)$ 
are easily described by 
$$p_{d_6}([z_1:z_2:z_4])=[z_2:z_4]\qquad\text{and}\qquad\ph_{d_6}([z_1:z_2:z_4:t])=[z_2:z_4].$$
Since $\ddo_6$ is a $D_6$-singularity of $X_\l^{31}$, the reduced fiber $(\pih_\l^{d_6})^{-1}(d_6)$ 
is isomorphic to $\PM^1(\CM)$ and contains two singular points of $\Xhat_\l^{31}$: one, which 
we denote by $a$, is an $A_1$ singularity and the other, which we denote by $b$, 
is a $D_4$-singularity. A {\sc Magma} computation shows that 
$$\phh_{d_6}(a)=[1:-4\l(\l+1)]\qquad \text{and}\qquad \phh_{d_6}(b)=[0:1]=\ph_{d_6}(a_3).$$
The singular fiber above $[0:1]$ has been described in Lemma~\ref{lem:g31-graph}(b) 
so we concentrate now on the fiber above $[1:-4\l(\l+1)]$. 

We denote by $\D_\l$ the closure of $\ph_{d_6}^{-1}([1:-4\l(\l+1)])$ in $X_\l^{31}$: 
if we denote by $s_\l(z_1,z_2)$ the quadratic form 
$$s_\l(z_1,z_2)=
z_1^2 + (-71\,\l^2 - 52\,\l - 8)\,z_1z_2 + (8\,\l^4 + 28\,\l^3 + 36\,\l^2 + 20\,\l + 4)\,z_2^2,$$
we have
\eqna
\D_\l&=&\{[z_1:z_2:z_4:t] \in X_\l^{31}~|~z_4=-4\l(\l+1)z_2\} \\
&\simeq& \{[z_1:z_2:t] \in \PM(1,1,3)~|~t^2=z_2r_\l(z_1,z_2,-4\l(\l+1)z_2)\} \\
&=& \{[z_1:z_2:t] \in \PM(1,1,3)~|~t^2=-16\,\l^3(2\l+1)^3\,z_2^4
s_\l(z_1,z_2)\}.
\endeqna
Note the following fact:

\bigskip

\begin{lemma}\label{lem:delta-l}
If $\l \neq 0$, $-1/2$, $-8/17$, then the closed subvariety $\D_\l$ meets the singular 
locus of $X_\l^{31}$ at only one point (the point $\ddo_6$). 
\end{lemma}

\bigskip

\begin{proof}
This is a {\sc Magma} computation.
\end{proof}

\bigskip

Let $\Delt_\l$ denote the strict transform of $\D_\l$ in $\Xti_\l^{31}$, 
recall that $\D_6^{d_6}=(\pih_\l^{d_6})^{-1}(a)$ 
and let $\EC_\l$ denote the fiber $(\pht_\l^{d_6})^{-1}([1:-4\l(\l+1)])$. 
Then it follows from Lemma~\ref{lem:delta-l} that:

\bigskip

\begin{cor}\label{coro:delta-l}
If $\l \neq 0$, $-1/2$, $-8/17$, then $\EC_\l=\Delt_\l \cup \D_6^{d_6}$.
\end{cor}

\bigskip

However, it must be noticed that $\D_\l$ is not necessarily irreducible. Indeed, 
$$s_\l=(z_1 - \frac{71\,\l^2 + 52\,\l + 8}{2}\,z_2)^2 -\l\Bigl(\frac{17\,\l+8}{4}\Bigr)^3 z_2^2.$$
So $\D_\l$ is irreducible if and only if $\l \neq 0,-8/17$ (we retrieve 
the same special value as in Lemma~\ref{lem:delta-l}). 
We deduce from this the following result:

\bigskip

\begin{cor}\label{coro:e-lambda}
If $\l \neq 0$, $-1/2$, $-8/17$, then $\EC_\l$ is a singular fiber of type $\Irm_2$
\end{cor}

\bigskip

\begin{proof}
The hypothesis implies that $\EC_\l$ contains two irreducible components, namely $\Delt_\l$ 
and $\D_6^{d_6}$. It then follows from the classification of singular fibers that $\EC_\l$ 
is of type $\Irm_2$ or $\Irm\Irm\Irm$. 

Now, let $\Delh_\l$ denote the strict transform of $\D_\l$ in $\Xhat_\l^{31}$. 
From the equation of $\D_\l$, we see that $d_6$ is an $A_3$ singularity of $\D_\l$ so that, 
after blowing-up, $a$ is an $A_1$ singularity of $\Delh_\l$. So, after blowing-up $a$, 
we see that $\Delt_\l$ meets $\D_6^{d_6}$ in two different points, so that $\EC_\l$ is of type $\Ati_1$.
\end{proof}

\bigskip

\begin{pro}\label{prop:fibres-singulieres-31}
Let $\l \in \CM^\times$. Then the singular fibers of $\pht_\l^{d_6}$ 
are given by Table~\ref{table:k3-g31-fibration}.
\end{pro}

\bigskip

\begin{proof}
Assume first that $\l \neq -1/2$, $-8/17$. Then a {\sc Magma} computation 
shows that, if $x$ and $y$ are two different singular points of $X_\l^{31}\setminus\{\ddo_6\}$, 
then $\ph_{d_6}(x) \neq \ph_{d_6}(y)$. Then the result follows from Lemmas~\ref{lem:g31-graph}(b) 
and~\ref{lem:delta-l} and the same argument based on Euler characteristic in the proof 
of Corollary~\ref{coro:singular-g30} to distinguish between the different 
possibilities.

\medskip

The case where $\l=-1/2$ will be treated in~\S\ref{sub:g31-quartic}. So it remains to check 
the case where $\l=-8/17$. The numerical facts in what follows can be 
checked with {\sc Magma}. 
Whenever $\l=-8/17$, then $\D_{-8/17}$ is not irreducible and contains 
one of the $A_2$ singularities of $X_{-8/17}^{31}$ (let us call it $\ado_2$), 
the singularity $\ddo_6$ and no other singular points of $X_{-8/17}^{31}$. 
and splits into two irreducible components which we call $\D_{-8/17}^1$ and $\D_{-8/17}^2$. 
Their intersection contains only the points $\ddo_6$ and $\ado_2$. One can check that 
they are both smooth at $\ado_3$ and that the tangent line of $\D_{-8/17}^1$ at $\ado_2$ 
is different than the one of $\D_{-8/17}^2$. Therefore, $\EC_{-8/17}$ is the union 
of five irreducible components $\D_6^{d_6}$, $\Delt_{-8/17}^1$, $\Delt_{-8/17}^2$, $\D_1^{a_2}$ 
and $\D_2^{a_2}$ and the last four form an $A_4$ configuration whose extremal curves  
are $\Delt_{-8/17}^1$ and $\Delt_{-8/17}^2$. Since these extremal curves both meet $\D_6^{d_6}$, 
the only possibility for the singular fiber $(\pht_\l^{d_6})^{-1}(\ph_{d_6}(\ado_2))$ 
is to be of type $\Ati_4$. The other singular fibers are obtained as in the previous case, 
using again Euler characteristic to remove ambiguities.
\end{proof}

\bigskip

Recall from Proposition~\ref{prop:section-g31}(b) that the elliptic fibration 
$\pht_\l^{d_6}$ admits a section whose image is $\D_5^{d_6}$:

\bigskip

\begin{pro}\label{prop:g31-mordell}
Let $\l \in \CM^\times$. Then the Mordell-Weil group $\mweil(\pht_\l^{d_6})$ 
is given by Table~\ref{table:k3-g31-fibration}.
\end{pro}

\bigskip

\begin{proof}
In all cases, the rank of the Mordell-Weil group 
is equal to $0$. The torsion is given by~\cite{shimada}.
\end{proof}

\bigskip

We summarize all the datas collected in this section and the previous one 
in Table~\ref{table:k3-g31-fibration}. Observe that in all the cases except when the Mordell-Weil group has torsion, 
the Picard group of the K3 surface is $U+(Dynkin\, diagram \, of \, the \, singular \, fibers)$, i.e. in the generic cas is 
$U+D_7+3 A_2+3 A_1$. In the case when the Mordell-Weil group 
is $\mathbb{Z}/2\mathbb{Z}$ then one has to add the $2$-torsion section to get the whole Picard group. 

\bigskip

\begin{table}\label{table4}
$$\begin{array}{!{\vline width 2pt} c !{\vline width 1.3pt} c|c|c|c|c!{\vline width 2pt}}
\hlinewd{2pt}
\l & \ZC_\sing(f_{3,\l}) & \text{singularities of $X_\l^{31}$} & 
\rhob(\Xti_\l^{31}) & \text{singular fibers of $\pht_\l^{d_6}$} & \mweil(\pht_\l^{d_6}) \petitespace
\\
\hlinewd{2pt}
\neq \l_k,-8/17& \vide & D_6+A_3+3\,A_2+2\,A_1 & 
\ge 18^{(\dagger)} & \Dti_7+3\,\Ati_2+3\,\Ati_1 & 0^{(\ddagger)}\petitespace \\ \hline
-8/17 & \vide & D_6+A_3+3\,A_2+2\,A_1 & 19 & \Dti_7+\Ati_4+2\,\Ati_2+2\,\Ati_1 & 0\petitespace\\ \hline 
\l_{145}=-8/25 & 960\,A_1 & D_6+D_5+A_3+2\,A_2 & 19 & \Dti_7+\Dti_5+2\,\Ati_2+\Ati_1 & 0\petitespace\\ \hline
\l_{245}=-81/175 & 480\,A_1 & E_6+D_6+A_3+A_2+A_1 & 19 & \Eti_6+\Dti_7+\Ati_2+2\,\Ati_1 & 0\petitespace\\ \hline
\l_1=1 & 1920\,A_1 & D_6+A_5+A_3+A_2+2\,A_1 & 19 & \Dti_7+\Ati_5+\Ati_2+3\,\Ati_1 & \ZM/2\ZM \petitespace\\ \hline
\l_2=-1/3 & 1440\,A_2 & D_6+D_5+3\,A_2+A_1 & 19 & \Dti_9+3\,\Ati_2+2\,\Ati_1 & 0\petitespace\\ \hline
\l_3=-1/2 & 640\,A_3 & D_6+2\,A_3 +2\,A_2+2\,A_1 & 19 & \Dti_7 + \Dti_5 + 2\,\Ati_2+\Ati_1 & 0\petitespace\\ 
\hlinewd{2pt}
\end{array}$$
\caption{Some numerical data for the family of K3 surfaces 
$(\Xti_\l^{31})_{\l \in \CM^\times}$}\label{table:k3-g31-fibration}
\vskip-0.4cm
$^{(\dagger)}$ With equality for $\l$ generic \\
$^{(\ddagger)}$ Only for $\l$ generic
\end{table}

\bigskip

\bigskip

\subsection{Three particular cases} 
We study here the cases where $\l \in \{-8/17,1,-1/2\}$, which are 
all particular in their own way.

\bigskip

\subsubsection{The case $\l=-8/17$}\label{sub:g31-8-17}
We assume here, and only here, that $\l=-8/17$. 
As shown in Proposition~\ref{prop:fibres-singulieres-31}, 
the elliptic fibration $\ph_{-8/17}^{d_6}$ of the K3 surface $X_{-8/17}^{31}$ 
has the property that $\D_\l$ contains a singular point of $X_{-8/17}^{31}$ different 
from $d_6$ and the corresponding singular fiber $\EC_{-8/17}$ is of type $\Ati_4$. 
This has the following consequence for its Picard number, which makes $\Xti_{-8/17}^{31}$ 
a special member of the family obtained from minimal resolutions of quotients by $G_{31}'$ 
of the smooth family of icosics 
$(\ZC(f_{3,\l}))_{\l \in \CM^\times\setminus\{\l_{145},\l_{245},\l_1,\l_2,\l_3\}}$:

\bigskip

\begin{pro}\label{prop:picard-8-17}
$\rhob(\Xti_{-8/17}^{31})=19$. 
\end{pro}

\bigskip

\begin{proof}
For proving that $\rhob(X_{-8/17}^{31}) \ge 19$, we shall use the elliptic 
fibration $\pht_\l^{d_6}$. Indeed, this fibration admits a section, so 
$\rhob(X_{-8/17}^{31}) \ge 2 + m'$, where $m'$ is the rank of the subgroup 
of $\Pic(X_{-8/17}^{31})$ generated by irreducible components of the singular fibers 
(here, $2$ comes from the section and a general smooth fiber of $\pht_\l^{d_6}$). 
It follows from Table~\ref{table:k3-g31-fibration} that $m'=17$, so
$$\rhob(\Xti_{-8/17}^{31}) \ge 19.$$
Now, proving that $\rhob(\Xti_{-8/17}^{31})=19$ is done as in the proof 
of Corollary~\ref{coro:picard-g31}, thanks to {\sc Magma} computations 
and the Artin-Tate Conjecture.
\end{proof}

\bigskip

\subsubsection{The case $\l=\l_1=1$}\label{sub:g31-cubic}
We assume here, and only here, that $\l=\l_1=1$. We set
$$q_1=z_1z_2 - 1/108\,z_2^2 + 1/54\,z_2z_4 - 1/108\,z_4^2$$
$$c_1=z_1^2z_2 - 54\,z_1z_2^2 + 1/8\,z_1^2z_4 - 9\,z_1z_2z_4 - 1/4\,z_1z_4^2 + 1/8\,z_4^3.
\leqno{\text{and}}$$
Then $q_1$ and $c_1$ are irreducible and 
$$r_1=-864 q_1c_1.$$
So, if we denote by $Q_1=\ZC(q_1)$ and $C_1=\ZC(c_1)$, 
then $Q_1$ is a smooth conic while $C_1$ is a cuspidal cubic. 
Then
\equat\label{eq:r1}
R_1= \Bba_2 \cup Q_1 \cup C_1.
\endequat
The singular points of $R_1$ are given by
$$d_6=[1 : 0 : 0],\quad 
a_3=[1 : 0 : 1],\quad
a_5=[1:3:21],\quad a_2 =[1:1/27:-1/3],$$
$$a_1^+=[1 : {\SS{\frac{231\,\sqrt{33} + 1327}{2}}}, {\SS{\frac{165\,\sqrt{33} + 949}{2}}}]
\qquad\text{and}\qquad 
a_1^-=[1 : {\SS{\frac{-231\,\sqrt{33} + 1327}{2}}}, {\SS{\frac{-165\,\sqrt{33} + 949}{2}}}].$$
It is easily checked that 
$$d_6 \in \Bba_2 \cap Q_1 \cap C_1,\quad a_3 \in (\Bba_2 \cap C_1)\setminus Q_1,$$
$$a_2 \in C_1 \setminus (\Bba_2 \cup Q_1),\quad 
a_5,a_1^\pm \in (Q_1 \cap C_1) \setminus \Bba_2.$$
We denote by $Q_1'$ the preimage of $Q_1$ in $X_1^{31}$, endowed with its 
reduced structure (so that $Q_1' \simeq Q_1$) and we denote by $\Qti_1'$ 
the strict transform of $Q_1'$ in $\Xti_1^{31}$. Since $Q_1'$ is a smooth 
rational curve, we get that 
\equat\label{eq:q1}
\Qti_1' \simeq Q_1' \simeq Q_1.
\endequat
Since the smooth conic $Q_1$ goes through the point $d_6$, we get that
$\Qti_1'$ is a section of the elliptic fibration $\pht_1^{d_6}$. 
By Table~\ref{table:k3-g31-fibration}, we get:

\bigskip

\begin{pro}\label{prop:q1-section}
The smooth rational curve $\Qti_1'$ is a section of the elliptic fibration $\pht_1^{d_6}$. 
It is $2$-torsion and generates the Mordell-Weil group $\mweil(\pht_1^{d_6})$. 
\end{pro}

\bigskip

\subsubsection{The case $\l=\l_3=-1/2$}\label{sub:g31-quartic}
We assume here, and only here, that $\l=-1/2$. We set
\eqna
r_{-1/2}^\circ & = & 27\,z_1^3z_2 + 947/16\,z_1^2 z_2^2 + 54\,z_1z_2^3 
- 113/2\,z_1^2z_2z_4 \\ && - 171/2\,z_1z_2^2z_4 - z_1^2 z_4^2 + 59/2\,z_1z_2 z_4^2 
+ 2\,z_1z_4^3 + z_2z_4^3 - z_4^4
\endeqna
Then $r_{-1/2}^\circ$ is irreducible and
\equat\label{eq:facteur-rdemi}
r_{-1/2}=(z_2-z_4)r_{-1/2}^\circ.
\endequat
Let $L$ denote the line in $\PM^2(\CM)$ defined by $\ZC(z_2-z_4)$ and let 
$R_{-1/2}^\circ = \ZC(r_{-1/2}^\circ) \subset R_{-1/2}$. Then
\equat
R_{-1/2}= \Bba_2 \cup L \cup R_{-1/2}^\circ.
\endequat
The singular points of $R_{-1/2}$ are
$$d_6=[1 : 0 : 0],\quad 
a_3=[1 : 0 : 1],\quad
a_3^L=[0 : 1 : 1],$$
$$a_2^+=[1 : {\SS{\frac{52\,\sqrt{13} - 184}{27}}} : {\SS{\frac{13\,\sqrt{13} - 37}{6}}}],\quad
a_2^-=[1 : {\SS{\frac{-52\,\sqrt{13} - 184}{27}}} : {\SS{\frac{-13\,\sqrt{13} - 37}{6}}}].$$
$$a_1^L=[1 : -16 : -16]\qquad\text{and}\qquad
a_1=[1 : 1/32 : 7/8],$$
With this notation, $d_6$ is a $D_6$ singularity, $a_3$ and $a_3^L$ are $A_3$ singularities, 
$a_2^+$ and $a_2^-$ are $A_2$ singularities and $a_1$ and $a_1^L$ are $A_1$ singularities 
of $R_{-1/2}$. Note that, as sets,
\equat\label{eq:inter-rdemi}
\Bba_2 \cap L = \{d_6\}, \quad \Bba_2 \cap R_{-1/2}^\circ = \{d_6,a_3\}
\quad\text{and}\quad L \cap R_{-1/2}^\circ = \{d_6,a_3^L,a_1^L\}.
\endequat

Let $L'$ denote the preimage of $L$ in $X_{-1/2}^{31}$ endowed with 
its reduced structure. Then
$$L'=\{[z_1:z_2:z_3:t] \in \PM(1,1,1,3)~|~t=z_2-z_4=0\} \simeq \PM^1(\CM).$$
We denote by $\Lti'$ its strict transform in $\Xti_{-1/2}^{31}$. 
Then $(\pht_{-1/2}^{d_6})^{-1}(\ph_\l(a_3))$ contains $\D_6^{d_6}$, $\Lti'$ 
and the exceptional divisors above the singularities $\ado_1^L$ 
and $\ado_3^L$: the smooth rational curve $\Lti'$ meets $\D_6^{d_6}$, 
the exceptional divisor above $\ado_1^L$ and at least one of the 
exceptional divisors above $\ado_3^L$, so the only possibility 
is that $(\pht_{-1/2}^{d_6})^{-1}(\ph_\l(a_3))$ is a singular fiber of 
type $\Dti_5$. 

The other singular fibers are now determined easily and fit 
with the data in Table~\ref{table:k3-g31-fibration}. Note also 
that $\Lti'$ provides another section of all the fibrations 
$\pht_{-1/2}^{a_1}$, $\pht_{-1/2}^{a_2^\pm}$ and $\pht_{-1/2}^{a_3}$. 

\bigskip

\begin{centerline}{\bf\large Appendix: Morphisms to $\PM^1(\CM)$}\end{centerline}

\bigskip

We describe here two basic contructions of morphisms to $\PM^1(\CM)$ 
which are used in the body of the article for constructing 
elliptic fibrations on our K3 surfaces.

\bigskip
\renewcommand\thesection{\Alph{section}}
\setcounter{section}{0}

\section{Weighted projective space}\label{appendix-2}

\medskip

\boitegrise{{\bf Notation.} {\it We fix two natural numbers $k$ and $l$ 
such that $\gcd(k,l)=1$, we set $m=k+l$ and we denote here by $p$ the 
point $[0:0:1]$ of $\PM(k,l,m)$. 
It is an $A_{m-1}$-singularity of $\PM(k,l,m)$. We denote 
by $\pi : \hat{\PM}(k,l,m) \to \PM(k,l,m)$ the minimal 
resolution of the singularity $p$.}}{0.75\textwidth}

\medskip

Note that we have only resolved the singularity $p$, so that 
$\pi^{-1}(\PM(k,l,m) \setminus \{p\})$ may still have two singular points 
(above $[1:0:0]$ and $[0:1:0]$). Let 
$$\fonction{\ph_{k,l}}{\PM(k,l,m) \setminus \{p\}}{\PM^1(\CM)}{[x:y:z]}{[x^l:y^k].}$$
Then there exists a unique morphism of varieties
$$\phh_{k,l} : \hat{\PM}(k,l,m) \longto \PM^1(\CM)$$
making the diagram
\equat\label{eq:map-p1}
\diagram
\hat{\PM}(k,l,m) \setminus \pi^{-1}(p) 
\ddto_{\DS{\pi}}^{\rotatebox{90}{$\sim$}} \ar@{^{(}->}[rr]&& 
\hat{\PM}(k,l,m) \ddto^{\DS{\phh_{k,l}}} \\
&& \\
\PM(k,l,m) \setminus\{p\} \rrto^{\DS{\ph_{k,l}}} && \PM^1(\CM)
\enddiagram
\endequat
commutative.

\medskip

\begin{proof}
The uniqueness is trivial, so let us prove the existence. It is sufficient 
to work in the affine chart $U_z=\{[x:y:z] \in \PM(k,l,m)~|~z \neq 0\}$ 
of $\PM(k,l,m)$. We denote by $\Uti_z$ its minimal resolution of singularities. 
Through the variables $a=x^{m}$, $b=xy$ and $c=y^{m}$ (and setting $z=1$), 
we have
$$U_z =\{(a,b,c) \in \AM^{\!3}(\CM)~|~b^{m}=ac\}$$
and $p$ corresponds to the point $0$ of $U_z$ while the restriction of $\ph_{k,l}$ 
to $U_z \setminus \{0\}$ is given by
$$\ph_{k,l}(a,b,c)=
\begin{cases}
[a : b^k] & \text{if $(a,b) \neq (0,0)$,}\\
[b^l : c] & \text{if $(b,c) \neq (0,0)$.}
\end{cases}$$
A model for the minimal resolution of $U_z$ is given 
by\footnote{For $0 \le j \le m-1$, let $J_j$ denote the ideal 
of the algebra 
$\CM[U_z]=\CM[A,B,C]/\langle B^{m}-AC \rangle$ generated by $A$ and $B^j$. 
Then $\Uti_z$ 
is the blowing-up of 
$J_0J_1\cdots J_{m-1}=\langle (A^{m-j} B^{j(j-1)/2})_{1 \le j \le m} \rangle$: 
the variable $u_j$ corresponds to the generator $A^{m-j} B^{j(j-1)/2}$.}
\begin{multline*}
\Uti_z = \{((a,b,c),[u_1:u_2:\dots:u_m]) \in U_z \times \PM^{m-1}(\CM)~|~\\
\begin{cases}
\forall~2 \le j \le m,~au_j=b^{j-1} u_{j-1},\\
\forall~1 \le j \le m-1,~cu_j=b^{m-j} u_{j+1},\\
\forall~1 \le j < j' \le m,~u_ju_{j'} = b^{j'-j-1} u_{j+1}u_{j'-1}
\end{cases} \}.
\end{multline*}
Note that the last equation is automatically fulfilled if $j'=j+1$. 

We then define $\phh_{k,l} : \Uti_z \longto \PM^1(\CM)$ by
\equat\label{eq:extension}
\phh_{k,l}((a,b,c),[u_1\!:\!u_2\!:\!\cdots\!:\!u_m])\!=\!
\begin{cases}
[u_j\!:\!b^{k-j}u_{j+1}] & \!\!\text{if $u_j \neq 0$ and $j \le k$,}\\
[b^{j-1-k} u_{j-1}\!:\!u_j] & \!\!\text{if $u_j \neq 0$ and $j \ge k+1$.} \\
\end{cases}
\endequat
An immediate computation from the equations of $U_z$ shows that 
$\phh_{k,l}$ is well-defined and satisfies the required property.
\end{proof}

\bigskip

Let $\D_1$,\dots, $\D_{m-1}$ be the smooth projective 
lines in the exceptional divisor $\pi^{-1}(p)$ and we assume that 
they are numbered so that, in the open subset $\Uti_z$ described in 
the proof of~(\ref{eq:map-p1}), 
$$\D_j=\{0\} \times \{[u_1:\cdots:u_m] \in \PM^{m-1}(\CM)~|~
\forall~r \in \{1,2,\dots,m\} \setminus \{j,j+1\},~u_r=0\}.$$
Let
$$\D_x=\{[x:y:z] \in \PM(k,l,m)~|~x=0\} \hphantom{A}\text{and}\hphantom{A}
\D_y=\{[x:y:z] \in \PM(k,l,m)~|~y=0\}.$$
Then $\D_x$ and $\D_y$ are smooth rational curves and $\D_x \cap \D_y = \{p\}$. 
Let $\Delt_x$ and $\Delt_y$ denote the respective 
strict transforms of $\D_x$ and $\D_y$ 
in $\hat{\PM}(k,l,m)$. 

\bigskip

\begin{pro}\label{prop:graphe-eclatement}
The fiber $\phh_{k,l}^{-1}([1:0])$ (resp. $\phh_{k,l}^{-1}([0:1])$) is the union 
of the smooth rational curves $\Delt_y$, $\D_1$,\dots, $\D_{k-1}$ 
(resp. $\D_{k+1}$,\dots, $\D_{m-1}$, $\Delt_x$). The intersection 
graphs are given respectively by
$$\text
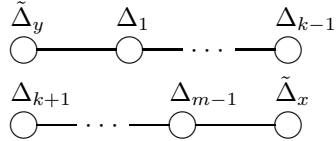
{\centerline{\begin{picture}(110,27)
\put(  5, 10){\circle{10}}\put(0,19){$\Delt_y$}
\put( 10, 10){\line(1,0){30}}
\put( 45, 10){\circle{10}}\put(40,19){$\D_1$}
\put( 50, 10){\line(1,0){15}}
\put(105, 10){\circle{10}}\put(100,19){$\D_{k-1}$}
\put(100, 10){\line(-1,0){15}}
\put( 68,7.3){$\cdot \cdot \cdot$}
\end{picture}}}$$
$$\text{\centerline{\begin{picture}(110,5)
\put(  5, 10){\circle{10}}\put(0,19){$\D_{k+1}$}
\put( 10, 10){\line(1,0){15}}
\put( 65, 10){\circle{10}}\put(60,19){$\D_{m-1}$}
\put( 60, 10){\line(-1,0){15}}
\put(105, 10){\circle{10}}\put(100,19){$\Delt_x$}
\put( 70, 10){\line(1,0){30}}
\put( 28,7.3){$\cdot \cdot \cdot$}
\end{picture}}}\leqno{\text{and}}$$
\end{pro}

\bigskip

\begin{proof}
One only needs to determine the intersections of the fibers 
$\phh_{k,l}^{-1}([1:0])$ and $\phh_{k,l}^{-1}([0:1])$ with the open set $\Uti_z$ of $\hat{\PM}(k,l,m)$. 
But this can be done from the explicit model and formula~\eqref{eq:extension} 
given in the proof of~\eqref{eq:map-p1}.
\end{proof}

\bigskip

\begin{remark}\label{rem:poids-section}
It follows from from~\eqref{eq:extension} that the restriction of $\phh_{k,l}$ 
to the smooth rational curve $\D_k$ is an isomorphism: this provides 
a section $\PM^1(\CM) \longto \hat{\PM}(k,l,m)$ to the morphism 
$\phh_{k,l} : \hat{\PM}(k,l,m) \longto \PM^1(\CM)$.\finl
\end{remark}

\bigskip

\begin{remark}\label{rem:p1-fibre-generale}
Let $p \in \PM^1(\CM) \setminus \{[1:0],[0:1]\}$. Then $\overline{\ph_{k,l}^{-1}(p)}$ 
is a smooth rational curve. Indeed, write $p=[1:\a]$ with $\a \neq 0$ (and assume 
that $k \le l$, the other case being similar). Then 
$$\overline{\ph_{k,l}^{-1}(p)}=\{[x:y:z] \in \PM(k,l,m)~|~y^k=\a x^l\}.$$
Working first in the affine chart $U_z$, we get
$$\overline{\ph_{k,l}^{-1}(p)} \cap U_z = \{(a,b,c) \in \AM^{\! 3}(\CM)~|~b^{m}=ac,~c=\a b^l~\text{and}~
a=\a^{-1} b^k\} \simeq \AM^{\! 1}(\CM).$$
Working now in the affine chart $U_x=\{[x:y:z] \in \PM(k,l,m)~|~x \neq 0\}$, 
we have
$$U_x=\{(u_0,u_1,\dots,u_k) \in \AM^{\! k+1}(\CM)~|~\forall 1 \le j \le j' \le k,
~v_jv_{j'} = v_{j-1} v_{j'+1}\}$$
(here, the variable $v_j$ stands for $y^j z^{k-j}$). 
Therefore, 
$$\overline{\ph_{k,l}^{-1}(p)} \cap U_x = \{(u_0,u_1,\dots,u_k) \in U_x~|~u_0=\a\} \simeq \AM^{\! 1}(\CM),$$
the isomorphism $\AM^{\! 1}(\CM) \longiso \overline{\ph_{k,l}^{-1}(p)} \cap U_x$ being given by 
$u \mapsto \a(1,u,u^2,\dots,u^{k-1})$. Since $\overline{\ph_{k,l}^{-1}(p)}$ is contained in 
$U_x \cup U_z$, we conclude that $\overline{\ph_{k,l}^{-1}(p)}$ 
is a smooth rational curve.\finl
\end{remark}

\bigskip

\section{Double cover of $\PM^2(\CM)$}\label{appendix-3}

\medskip

\boitegrise{{\bf Hypothesis and notation.} {\it We fix a non-zero natural number $m$ and a 
square-free homogeneous 
polynomial $F \in \CM[a,b,c]$ of degree $2m$, where $a$, $b$ and $c$ are of degree $1$. We denote by 
$\XC$ the surface 
$$\XC=\{[a:b:c:t] \in \PM(1,1,1,m)~|~t^2=F(a,b,c)\}$$
and by $\xi : \XC \to \PM^2(\CM)$, $[a:b:c:t] \mapsto [a:b:c]$. Let $\s$ denote the involutive 
automorphism of $\PM(1,1,1,m)$ defined by $\s([a:b:c:t])=[a:b:c:-t]$.}}{0.75\textwidth}

\bigskip

Then $\s$ stabilizes $\XC$ and $\xi$ is the double cover of $\PM^2(\CM)$ associated with $\s$. 
We denote by $\RC \subset \PM^2(\CM)$ its branch locus
$$\RC=\{[a:b:c] \in \PM^2(\CM)~|~F(a,b,c)=0\}.$$
If $x=[a:b:c] \in \RC$, we denote by $\xdo=[a:b:c:0]$ its unique preimage in $\XC$. 
We also define $p_x : \PM^2(\CM) \setminus \{x\} \longto \PM^1(\CM)$ to be 
the projection from $x$ and let $\b_x : \hat{\PM}_x^2(\CM) \longto \PM^2(\CM)$ denote the blow-up 
of $\PM^2(\CM)$ at $x$. Then the map 
$p_x \circ \b_x : \hat{\PM}_x^2(\CM)\setminus\b_x^{-1}(x) \longto \PM^1(\CM)$ extends uniquely 
to a morphism 
$$\phat_x : \hat{\PM}_x^2(\CM) \longto \PM^1(\CM),$$
which admits a section $\shat_x : \PM^1(\CM) \longto \hat{\PM}_x^2(\CM)$ whose 
image is $\b_x^{-1}(x) \simeq \PM^1(\CM)$. 
Finally, we denote by $\pi_x : \XCh_x \longto \XC$ the blow-up of $\XC$ at $\xdo$. 

\bigskip

\begin{pro}\label{prop:blowing-double-cover}
Assume that $x$ is a singular point of the branch locus $\RC$. Then the morphism 
$\xi : \XC \longto \PM^2(\CM)$ lifts uniquely to a morphism $\xih_x : \XCh_x \longto \hat{\PM}_x^2(\CM)$ 
making the diagram
$$\diagram
\XCh_x \rrto^{\DS{\xih_x}} \ddto_{\DS{\pi_x}} && \hat{\PM}_x^2(\CM) \ddto^{\DS{\b_x}} \\
&& \\
\XC \rrto^{\DS{\xi}} && \PM^2(\CM)
\enddiagram$$
commutative.
\end{pro}

\bigskip

\begin{remark}
The reader can easily check that, if $x$ is not a singular point of $R$, then 
the conclusion of proposition fails.\finl
\end{remark}

\bigskip

\begin{proof}
The uniqueness is clear, we only need to show the existence. 
By a linear change on the coordinates $a$, $b$, $c$, we may assume that $x=[0:0:1]$. 
It is sufficient to work in the open subsets $\UC$ and $U$ of $\XC$ and 
$\PM^2(\CM)$ defined by $c \neq 0$ (and we denote by $\UCh_x$ and $\Uhat_x$ their 
respective blow-up at $\xdo$ and $x$). 
We set $F_c(a,b)=F(a,b,1)$. Then, since $x$ is a singular point of $R$, this means 
that $F_c(0,0)=0$ and that the homogeneous component of degree $1$ of $F_c$ is zero. 
We can then write uniquely 
$$F_c(a,b)=a^2 \l(a,b) + ab \mu(b) + b^2\nu(b)$$
with $\l \in \CM[a,b]$ and $\mu$, $\nu \in \CM[b]$. Therefore:
\begin{multline*}
\UCh_x=\{((a,b,t),[A:B:T]) \in \AM^3(\CM) \times \PM^2(\CM)~|~
(a,b,t) \in [A:B:T]~\\ \text{and}~T^2=A^2\l(a,b)+AB \mu(b) + B^2 \nu(b)\}
\end{multline*}
$$\Uhat_x=\{((a,b),[A:B]) \in \AM^2(\CM) \times \PM^1(\CM)~|~(a,b) \in [A:B]\}.\leqno{\text{and}}$$
Then the map $\xih_x : \UCh_x \longto \Uhat_x$ defined by
$$\xih_x((a,b,t),[A:B:T])=((a,b),[A:B])$$
is well-defined and satisfies the requirements of the proposition.
\end{proof}

\bigskip

Proposition~\ref{prop:blowing-double-cover} allows to define a morphism 
$$\phh_x = \phat_x \circ \xih_x : \XCh_x \longto \PM^1(\CM).$$
We now investigate the question of sections of this morphism, whenever 
$x$ is an ADE singularity of the branch locus $\RC$ 
(by~\cite[Proposition~\ref{prop:double cover}]{bonnafe-sarti-1}, this implies that 
$\xdo$ is a simple singularity of the surface $\XC$ of the same type as $x$). 
First, note that, if $s_x : \PM^1(\CM) \longto \XCh_x$ is a section of $\phh_x$, then 
$\xih_x \circ s_x$ is a section of $\phat_x$. 
Therefore, the question amounts to study sections of the morphism 
$\xih_x : \pi_x^{-1}(\xdo) \longto \b_x^{-1}(x) \simeq \PM^1(\CM)$. 
Here, we endow $\pi_x^{-1}(\xdo)$ with its reduced structure. 
A first answer is given in the next proposition:

\bigskip

\begin{pro}\label{prop:section-double-cover}
Let $x$ be an ADE singularity of $\RC$. Then:
\begin{itemize}
\itemth{a} If $x$ is an $A_1$ singularity, then $\pi_x^{-1}(\xdo) \simeq \PM^1(\CM)$ 
and the morphism $\xih_x : \pi_x^{-1}(\xdo) \longto \b_x^{-1}(x)$ is a double 
cover admitting no section. 

\itemth{b} If $x$ is an $A_m$ singularity with $m \ge 2$, then $\pi_x^{-1}(\xdo)$ is the union 
of two smooth rational curves ($\simeq \PM^1(\CM)$) intersecting transversally at one point, 
and both smooth rational curves map isomorphically to $\b_x^{-1}(x)$. 
This gives two sections of $\xih_x$, each one being obtained from the other by 
composing with the involution $\s$.

\itemth{c} If $x$ is a DE singularity, then $\pi_x^{-1}(\xdo) \simeq \PM^1(\CM)$ mapping 
isomorphically on $\b_x^{-1}(x)$. This gives one section of $\xih_x$.
\end{itemize}
\end{pro}

\bigskip

\begin{proof}
As in the proof of the previous Proposition~\ref{prop:blowing-double-cover}, 
we may assume that $x=[0:0:1]$ and we keep the notation introduced in this above proof. 
In particular, 
$$\pi_x^{-1}(\xdo) \simeq \{[A:B:T] \in \PM^2(\CM)~|~T^2=\a(0,0)A^2+\b(0)AB+\g(0)B^2\}$$
and $\xih_x([A:B:T])=[A:B]$. Let us examine the different cases.

\medskip

(a) If $x$ is an $A_1$ singularity, then a linear change of coordinates in $a$, $b$ 
allows to assume that $\a(0,0)=\g(0)=0$ and $\b(0)=1$. Then 
$$\pi_x^{-1}(\xdo) \simeq \{[A:B:T] \in \PM^2(\CM)~|~T^2=AB\}$$
and the result follows.

\medskip

(b) If $x$ is an $A_m$ singularity with $m \ge 2$, 
then a linear change of coordinates in $a$, $b$ 
allows to assume that $\a(0,0)=1$ and $\b(0)=\g(0)=0$. 
Then 
$$\pi_x^{-1}(\xdo) \simeq \{[A:B:T] \in \PM^2(\CM)~|~T^2=A^2\} = \D_+ \cup \D_-,$$
where $\D_\pm=\{[A:B:T] \in \PM^2(\CM)~|~A=\pm T\}$. The result follows.

\medskip

(c) If $x$ is a DE singularity, then $\a(0,0)=\b(0)=\g(0)=0$, so 
$$\pi_x^{-1}(\xdo) \simeq \{[A:B:T] \in \PM^2(\CM)~|~T=0\} \simeq \PM^1(\CM),$$
so the result follows.
\end{proof}

\bigskip

Let us go on with the case where $x$ is an ADE singularity of $\RC$. 
We denote by $\pit_x : \XCt_x \longto \XC$ the resolution of $\XC$ 
{\it only at the point $\xdo$}. It factorizes through 
$\XCt_x \longto \XCh_x \longto \XC$. Let $m$ denote the Milnor number 
of $\xdo$. Then $\pit^{-1}_x(\xdo)$ is the union of $m$ smooth rational 
curves whose intersection graph is denoted by $\G_x$. If $x$ is not of type $A_1$, 
we denote by $\G_x^\#$ the graph obtained from $\G_x$ by removing 
the smooth rational curves which are mapped isomorphically to $\PM^1(\CM)$ under $\phh_x$. 
According to the discussion 
of Proposition~\ref{prop:section-double-cover}, easy computations 
give the following consequences about the behaviour of $\phh_x$ 
and the action of $\s$ on the corresponding 
intersection graph (here, type $D_2$ means type $A_1 \times A_1$ and type $D_3$ 
coincides with type $A_3$):

\bigskip

\begin{cor}\label{coro:section-double-cover}
All the smooth rational curves belonging to the same connected 
component of $\G_x^\#$ are mapped to the same point of $\PM^1(\CM)$ under 
$\phh_x$. If two smooth rational curves do $\G_x^\#$ do not belong to the 
same connected component, then they are mapped to different points of 
$\PM^1(\CM)$ under $\phh_x$. Moreover:
\begin{itemize}
\itemth{a} If $x$ is an $A_1$ singularity, then $\phh_x : \D_1^x \to \PM^1(\CM)$ 
is a double cover corresponding to the quotient by the action of $\s$.

\itemth{b} If $x$ is an $A_m$ singularity with $m \ge 2$, then $\G_x^\#$ 
is of type $A_{m-2}$ and $\s$ acts on $\G_x$ by the unique 
non-trivial involutive automorphism.

\itemth{c} If $x$ is a $D_m$ singularity with $m \ge 4$, then $\G_x^\#$ 
is of type $D_{m-2} \times A_1$. 
Moreover, $\s$ acts on the intersection graph as the identity if $m$ is even 
and as the unique non-trivial involutive automorphism if $m$ is odd.

\itemth{d} If $x$ is an $E_6$ singularity, then $\G_x^\#$ is of type $A_5$ 
and $\s$ acts on $\G_x$ as the unique non-trivial involutive automorphism.

\itemth{e} If $x$ is an $E_7$ (resp. $E_8$) singularity, then $\G_x^\#$ is of type $D_6$ 
(resp. $E_7$) and $\s$ acts on $\G_x$ as the identity.
\end{itemize}
\end{cor}

\bibliographystyle{amsplain}
\bibliography{Biblio}

\end{document}